\documentclass[10pt]{amsart}

\usepackage{multirow}
\usepackage{amsmath,fixltx2e}
\usepackage{amsthm}
\usepackage{amssymb}
\usepackage{amsfonts}
\usepackage{mathtools}
\usepackage{MnSymbol}

\newtheorem{theorem}[equation]{Theorem}
\newtheorem{lemma}[equation]{Lemma}
\newtheorem{proposition}[equation]{Proposition}
\newtheorem{corollary}[equation]{Corollary}

\theoremstyle{definition}
\newtheorem{definition}[equation]{Definition}

\theoremstyle{remark}
\newtheorem{remark}[equation]{Remark}

\numberwithin{equation}{section}

\newcommand{\cH}{\mathcal{H}}

\usepackage[utf8]{inputenc}
\usepackage{graphicx}
\usepackage[bookmarks=false]{hyperref}
\usepackage{tikz-cd}
\allowdisplaybreaks

\title[An Equivariant Main Conjecture and Applications]{An Unconditional Equivariant Main Conjecture \linebreak in Iwasawa Theory and Applications}

\author{Rusiru Gambheera}
\author{Cristian D. Popescu}
\address{Dept. of Mathematics, University of California at San Diego, San Diego, CA 92093, USA}
\email{rgambhee@ucsd.edu}
\email{cpopescu@ucsd.edu}



\thanks{}

\date{}

\begin{document}
\begin{abstract} We prove a stronger version of the keystone result of Dasgupta-Kakde \cite{Dasgupta-Kakde} on the $\Bbb Z[G(H/F)]^-$--Fitting ideals of certain Selmer modules $Sel_S^T(H)^-$ associated to an abelian, CM extension $H/F$ of a totally real number field $F$ and use this to compute the
$\Bbb Z_p[[G(H_\infty/F)]]^-$--Fitting ideal of the Iwasawa module analogues $Sel_S^T(H_\infty)_p^-$ of these Selmer modules, where $H_\infty$ is the cyclotomic $\Bbb Z_p$--extension of $H$, for an odd prime $p$. Our main Iwasawa theoretic result states that the $\Bbb Z_p[[G(H_\infty/F]]^-$--module $Sel_S^T(H_\infty)_p^-$ is of projective dimension $1$ (unlike the $\Bbb Z[G(H/F)]^-$--module $Sel_S^T(H)^-$ which could have infinite projective dimension), is quadratically presented, and that its Fitting ideal is principal, generated by an equivariant $p$--adic $L$--function $\Theta_S^T(H_\infty/F)$. Further, we establish a perfect duality pairing between $Sel_S^T(H_\infty)_p^-$ and a certain $\Bbb Z_p[[G(H_\infty/F)]]^-$--module $\mathcal M_S^T(H_\infty)^-$, essentially introduced by Greither and the second author in \cite{Greither-Popescu}. As a consequence, we recover the Equivariant Main Conjecture for the Tate module $T_p(\mathcal M_S^T(H_\infty))^-$, proved in loc.cit. under the hypothesis that the classical Iwasawa $\mu$--invariant associated to $H$ and $p$ vanishes. As a further consequence,
we give an unconditional proof of the refined Coates--Sinnott Conjecture, proved in loc.cit. under the same $\mu=0$ hypothesis, and also proved unconditionally but with different methods by Johnston--Nickel in \cite{Johnston-Nickel-EMC}, regarding the $\Bbb Z[G(H/F)]$--Fitting ideals of the higher Quillen $K$--groups $K_{2n-2}(\mathcal O_{H,S})$, for all $n\geq 2$. Finally, we combine the techniques developed in the process with the method of ``Taylor--Wiles primes'' (introduced by Wiles in \cite{Wiles-Brumer} and refined by Greither in \cite{Greither-Brumer}) to strengthen further the keystone result in \cite{Dasgupta-Kakde} and prove, as a consequence, a conjecture of Burns--Kurihara--Sano on Fitting ideals of Selmer groups of CM number fields.

\end{abstract}

\date{\today}

\maketitle

\section{Introduction}
Let $H/F$ be a finite, abelian, CM extension of a totally real number field $F$ and let $S$ an $T$ be finite, disjoint sets of places in $F$, such that $S$ contains the set $S_\infty(F)$ 
of archimedean places in $F$.  Let $G:={\rm Gal}(H/F)$. From the data $(K/k, S, T)$, one can construct a meromorphic $G$--equivariant Artin $L$--function 
$$\Theta_{S, H/F}^T: \Bbb C\to \Bbb C[G],$$
which is holomorphic on $\Bbb C\setminus\{1\}$, and a Selmer $\Bbb Z[G]$--module $Sel_S^T(H)$, as defined by Burns--Kurihara--Sano in \cite{BKS}. The Selmer $\Bbb Z[G]$--module is finitely generated and contains subtle information regarding the $T$--ray-class group $Cl^T(H)$ of $H$. In particular, there is a canonical surjective morphism and an isomorphism
$$Sel_{S}^T(H)^-\twoheadrightarrow Sel_{S_\infty}^T(H)^-\simeq Cl^T(H)^{-, \vee},$$
where $M^\vee$ stands for the Pontrjagin dual and $M^-$ stands for the $(-1)$--eigenspace with respect to the action of the unique complex conjugation automorphism $j$ of $H$ on $M\otimes\Bbb Z[1/2]$, for any $\Bbb Z[G]$--module $M$.
\\

Let $p>2$ be a prime and let $S_{ram}(H/F)$, $S_\infty(F)$, $S_p(F)$ denote the ramification locus of $H/F$, the set of infinite primes, respectively $p$--adic primes in $F$. We will work under the following hypothesis on the set $T$:
\\

\noindent $\mathbf{{\rm\bf Hyp}_p^T(H/F)}$: $T\not\subseteq S_{ram}(H/F)$ and there are no non--trivial $p$--power roots of unity in $H$ which are congruent to $1$ modulo all the primes in $T\setminus S_{ram}(H/F)$.
\\

The keystone result in \cite{Dasgupta-Kakde} (Theorem  3.3 in loc.cit.) can be stated as follows.
\begin{theorem}[Dasgupta--Kakde, \cite{Dasgupta-Kakde}]\label{keystone-DK} If $S$ and $T$ satisfy the following properties
\begin{equation}\label{DK-hyp} S=S_\infty(F)\cup (S_{ram}(H/F)\cap S_p(F)), \quad S_{ram}(H/F)\subseteq S\cup T, \quad {\rm Hyp}_p^T(H/F),\end{equation}
then the $\Bbb Z_p[G]^-$--module  $Sel_S^T(H)_p^-:=(Sel_S^T(H)\otimes\Bbb Z_p)^-$ is quadratically presented and we have 
$${\rm Fitt}_{\Bbb Z_p[G]^-} Sel_S^T(H)_p^-\simeq (\Theta_{S,H/F}^T(0)).$$
\end{theorem}
This is a powerful result which implies the $p$--primary part of a Conjecture of Kurihara, stating that 
\begin{equation}\label{Kurihara-conj}{\rm Fitt}_{\Bbb Z_p[G]^-} (Cl^T(H)_p^{-, \vee})={\rm SKu}^T(H/F)_p^-\end{equation}
where the right--side is the Sinnott--Kurihara ideal, explicitly described in terms of special values of equivariant $L$--functions (see Theorem 3.5 in loc.cit.). Further, this implies the $p$--primary part of the Strong Brumer-Stark Conjecture, stating that {\it if, in addition, $S_{ram}(H/F)\subseteq S$,} then 
\begin{equation}\label{strong-BS}\Theta_{S,H/F}^T(0)\in {\rm Fitt}_{\Bbb Z_p[G]^-} (Cl^T(H)_p^{-, \vee}).\end{equation}
(See Corollary 3.8. in loc.cit.) The methods and techniques developed in \cite{Dasgupta-Kakde} are not Iwasawa theoretic in nature. They rely on the Ritter--Weiss theory \cite{Ritter-Weiss} of generalized Tate sequences,  a theory of Hilbert modular forms with group--ring values, introduced by Wiles in \cite{Wiles} and developed further in \cite{Dasgupta-Kakde} and \cite{Silliman}, as well as a generalization of Ribet's method \cite{Ribet} for constructing large unramified extensions of CM number fields.
\\

Prior to the publication of \cite{Dasgupta-Kakde}, the $p$--primary part of the Strong--Brumer Stark Conjecture \eqref{strong-BS}, for all primes $p>2$,  had been tackled via Iwasawa theoretic methods and proved under the additional assumptions that $S_p(F)\subseteq S$ and that the Iwasawa $\mu$--invariant $\mu_{H,p}$ associated to the cyclotomic $\Bbb Z_p$--extension $H_\infty/H$ vanishes, as conjectured by Iwasawa. (See \cite{Greither-Popescu}.) The hypothesis $\mu_{H,p}=0$ was necessary so that a certain Iwasawa module $\mathcal M_S^T(H_\infty)^-$ defined essentially  in \cite{Greither-Popescu} (see \S7 below for a definition in a more general setting) is $p$--divisible, (automatically) of finite $\Bbb Z_p$--corank. This implies that one has a canonical  isomorphism of  $\Bbb Z_p[[G(H_\infty/F)]]^-$--modules 
\begin{equation}\label{GP-iso}\mathcal M_S^T(H_\infty)_p^{-,\vee}\simeq T_p(\mathcal M_S^T(H_\infty)^{-})^\ast\end{equation}
between the Pontrjagin dual and, respectively, the $\Bbb Z_p$--dual of the Tate module of $\mathcal M_S^T(H_\infty)^-$. The results in \cite{Greither-Popescu} imply that the module $\mathcal M_S^T(H_\infty)^{-,\vee}$ is $G(H_\infty/F_\infty)$--cohomologically trivial, in general (no $\mu_{H,p}=0$ assumption necessary.) This fact, combined with its $\Bbb Z_p$--freeness (assuming $\mu_{H,p}=0$) shows that  the module $\mathcal M_S^T(H_\infty)^{-,\vee}$ is $\Bbb Z_p[G(H_\infty/F_\infty)]$--projective, and therefore of projective dimension $1$ over the equivariant Iwasawa algebra $\Bbb Z_p[[G(H_\infty/F)]]^-$. This result allows for the use of homological algebra techniques  developed in \cite{Greither-Popescu-Picard} to compute the 
Fitting ideal of this module over the equivariant Iwasawa algebra. When combining these techniques with links between  $T_p(\mathcal M_S^T(H_\infty)^{-})$ and the classical Iwasawa modules $\mathfrak X_S$, as well as the non-equivariant  main conjecture 
proved by Wiles in \cite{Wiles} for $\mathfrak X_S$, the authors of \cite{Greither-Popescu} obtain the equality
\begin{equation}\label{EMC-GP}{\rm Fitt}_{\Bbb Z_p[[G(H_\infty/F]]^-}\mathcal M_S^T(H_\infty)^{-,\vee}=(\Theta_S^T(H_\infty/F)),\end{equation}
where $\Theta_S^T(H_\infty/H)=\{\Theta_{S, H_n/F}^T(0)\}_n$ is the equivariant $p$--adic $L$--function, obtained by taking the projective limit in the Iwasawa algebra of the special values $\Theta_{S,H_n/F}^T(0)$ at the finite levels $H_n/F$ of the tower $H_\infty/F$. 

Equality \eqref{EMC-GP} is the Equivariant Main Conjecture (EMC) for the module  $\mathcal M_S^T(H_\infty)_p^{-,\vee}$, proved conditionally in \cite{Greither-Popescu}, under the assumption that $\mu_{H,p}=0$. Further, it is proved in \cite{Greither-Popescu} that, via Iwasawa codescent, \eqref{EMC-GP} implies the $p$--primary parts of the Strong Brumer--Stark Conjecture \eqref{strong-BS} and Refined Coates--Sinnott Conjecture, stating that
\begin{equation}\label{strong-CS}{\rm Fitt}_{\Bbb Z_p[G]}(K_{2n-1}(\mathcal O_{H,S})^{tors}_p)\Theta_S^T(1-n)=e_n\cdot {\rm Fitt}_{\Bbb Z_p[G]}(K_{2n-2}(\mathcal O_{H,S})_p), \qquad \text{ for all } n\geq 2,\end{equation}
if $S_p(F)\subseteq S$, which is a serious (imprimitivity) assumption for \eqref{strong-BS}, but not significant for \eqref{strong-CS}. Above, $K_i(O_{H,S})$ denotes the $i$--th Quillen $K$--group of the ring $\mathcal O_{H,S}$ of $S$--integers in $H$ and $e_n=1/2(1+(-1)^nj)$.

\begin{remark} In \cite{Greither-Popescu}, the authors prove an EMC for the Iwasawa module $T_p(\mathcal M_S^T(H_\infty)^-)$, but taking into account isomorphism \eqref{GP-iso} and the finite projective dimension of these modules (which makes the Fitting ideals of $T_p(\mathcal M_S^T(H_\infty)^-)$ and its $\Bbb Z_p$--dual equal), this is equivalent to an EMC for the module $\mathcal  M_S^T(H_\infty)^{-,\vee}$.
\end{remark}
\medskip

The main goal of this paper is to improve upon Theorem \ref{keystone-DK} of Dasgupta-Kakde {\it by relaxing hypotheses \eqref{DK-hyp} on $S$ and $T$,}  and use the improved version to prove unconditionally an Equivariant Main Conjecture in the spirit of \cite{Greither-Popescu} for the Selmer modules associated to the  cyclotomic $\Bbb Z_p$--extensions $H_\infty/H$:
$$Sel_S^T(H_\infty)_p^-:=\varprojlim_n Sel_S^T(H_n)_p^-,$$
where the projective limit is taken with respect to certain canonical restriction maps, giving it a natural $\Bbb Z_p[[G(H_\infty/F)]]^-$--module structure. Below, we let $\mathcal G:=G(H_\infty/F)$.
\smallskip

The first strengthening of Theorem \ref{keystone-DK} we obtain in this paper (see Theorem \ref{keystone-improved} below) is the following.
\begin{theorem}\label{keystone-improved-intro}
If $S$ and $T$ satisfy the following properties
\begin{equation}\label{improved-hyp} S_{\infty}(F)\cup S_p(F)\subseteq S,\quad  S_{ram}(H/F)\subseteq S\cup T,\qquad  {\rm Hyp}_p^T(H/F),\end{equation} 
then the $\Bbb Z_p[G]^-$--module $Sel_S^T(H)_p^-$ is quadratically presented and we have
$${\rm Fitt}_{\mathbb{Z}_p[G]^-}(Sel_S^T(H)_p^-)=(\Theta_{S, H/F}^T(0)).$$
\end{theorem}

Observe that the main difference between the old hypotheses \eqref{DK-hyp} and the new ones \eqref{improved-hyp} is allowing non $p$--adic primes (ramified or not) in $S$. Moving non $p$--adic ramified primes from $T$ to $S$ in the theorem above is a technically difficult 
task. This is achieved in several stages in \S6 below, through Theorems \ref{partial}, \ref{full} and \ref{keystone-improved}. The key ingredient is a construction of canonical maps between the quadratic presentations of the Selmer modules
$Sel_S^T(H')_p^-$ and   $Sel_S^T(H)_p^-$ for two abelian CM extensions $H'/F$ and $H/F$, with $H\subseteq H'$, achieved in \S5 and especially Theorem \ref{QP-compatibility}, which allows us to prove further that, under certain hypotheses, the Selmer groups in question behave well with respect to taking 
Galois--coinvariants (see the key Corollary \ref{covariance} below for the precise statement.) The results and techniques developed in \S5 are of interest in their own right and will be used further in our upcoming paper on an Iwasawa theory for Ritter--Weiss modules.
\smallskip

It is also worth noting that the methods of \cite{Dasgupta-Kakde} cannot be used (at least directly) to prove Theorem \ref{keystone-improved-intro}. For instance, when the set $S$ ($\Sigma$ in the notation of \cite{Dasgupta-Kakde}) contains non $p-$adic ramified primes, the proof of Proposition 9.5 in \cite{Dasgupta-Kakde}, which is central, will no longer work. 
\smallskip

The main result of this paper is the following unconditional Equivariant Main Conjecture for the Iwasawa module $Sel_S^T(H_\infty)_p^-$ in the spirit of \cite{Greither-Popescu}. (See Theorem \ref{EMC} below.)
\begin{theorem}[EMC] \label{EMC-intro}  Suppose that $S$ and $T$ satisfy hypotheses \eqref{improved-hyp}. Then, the following hold.
\begin{enumerate}
\item We have an equality of $\Bbb Z_p[[\mathcal G]]^-$--ideals
$${\rm Fitt}_{\mathbb{Z}_p[[\mathcal{G}]]^-}(Sel_S^T(H_{\infty})_p^-)=(\Theta_S^T(H_\infty/F)).$$
\item The $\Bbb Z_p[[\mathcal G]]^-$--module $Sel_S^T(H_{\infty})_p^-$ is quadratically presented and of projective dimension $1$.
\end{enumerate}
\end{theorem}

\noindent Part (1) of the above Theorem is obtained by applying Theorem \ref{keystone-improved-intro} to all the extensions $H_n/F$ and, via an application of a theorem of Greither and Kurihara \cite{Greither-Kurihara} (see Proposition \ref{limit-of-Fitt} below), taking a projective limit as $n\to\infty$. Part (2) follows from part (1), after proving the important fact that 
$\Theta_S^T(H_\infty/F)$ is a non zero--divisor in $\Bbb Z_p[[\mathcal G]]^-$ (see Proposition \ref{Theta-non-zero-divisor} below) and applying a slight improvement (see Proposition \ref{Cornacchia-Greither} below) of a homological algebra result due to Cornacchia and Greither \cite{Cornacchia-Greither} 
\\

The last section of this paper (\S7) is concerned with applications of the above results, with a special focus on linking the Selmer modules $Sel_S^T(H_\infty/H)_p^-$ and the Iwasawa modules introduced in \cite{Greither-Popescu} (called abstract Picard $1$--motives in loc.cit.) and providing unconditional proofs for the results proved in loc.cit. under the hypothesis that $\mu_{H,p}=0$. More precisely, in \S7.1 we work under the (weaker) hypotheses that $S\cap T=\emptyset$ and $S_p(F)\cup S_\infty(F)\subseteq S$ and construct a $\Bbb Z_p[[\mathcal G]]$--module $\mathcal M_S^T(H_\infty)$, which is $\Bbb Z_p$--torsion and had been essentially constructed in \cite{Greither-Popescu} under the stronger hypotheses that $S_{ram}(H_\infty/F)\subseteq S$ and $\mu_{H,p}=0$. Via an explicit perfect duality pairing (see \S7.1), we establish an isomorphism of $\Bbb Z_p[[\mathcal G]]^-$--modules 
\begin{equation}\label{iso-sel-mst}Sel_S^T(H_\infty)_p^-\simeq \mathcal M_S^T(H_\infty)^{\vee, -},\end{equation}
(See Proposititon \ref{link-sel-mst} below.) Under the additional hypthesis that $\mu_{H,p}=0$ (making the $\Bbb Z_p$--module $\mathcal M_S^T(H_\infty)^{\vee, -}$ divisible of finite corank), the isomorphism above becomes
\begin{equation}\label{iso-sel-tpmst} Sel_S^T(H_\infty)_p^-\simeq T_p(\mathcal M_S^T(H_\infty)^{-})^\ast,\end{equation}
which brings forth the module $T_p(\mathcal M_S^T(H_\infty)^-)$ (the $p$--adic realization of the abstract Picard $1$--motive $\mathcal M_S^T(H_\infty)^-$) studied in \cite{Greither-Popescu}. When combining our Theorem \ref{EMC-intro} with isomorphism \eqref{iso-sel-mst} (respectively, \eqref{iso-sel-tpmst}) we obtain an unconditional proof (respectively, a proof) of the main result in \cite{Greither-Popescu}, under weaker hypotheses on $S$ and $T$ than those in loc.cit. (See Corollary \ref{true-Greither-Popescu}, respectively Theorem \ref{Greither-Popescu-Main} below and keep in mind that in loc.cit. one has $S_{ram}(H_\infty/F)\subseteq S.$) 

Further, in \S7.3, we work under the hypothesis that $S_{ram}(H_\infty/F)\subseteq S$ and combine the technical machinery developed in \S6 of \cite{Greither-Popescu} with isomorphism \eqref{iso-sel-mst} to establish links between the Selmer modules $Sel_S^T(H_\infty)_p^-$ and the \'etale cohomology groups ${\rm H}^i_{et}(\mathcal O_{K,S}, \Bbb Z_p(n))$, for $i=1,2$. (See Proposition \ref{link-etale-selmer} and exact sequences \eqref{long-seq-sel-twist} and \eqref{coinvariant-seq}.) When combined with our Theorem \ref{EMC-intro} and the strategy developed in \S6 of \cite{Greither-Popescu}, these links lead to an unconditional proof of the Refined Coates--Sinnott Conjecture \eqref{strong-CS}. (See Theorem \ref{CS-K} below.) Equalities \eqref{strong-CS} had been established in \S6 of \cite{Greither-Popescu}, under the hypothesis $\mu_{H,p}=0$.

Finally, in \S7.4, we combine the techniques developed in the process of proving Theorem \ref{keystone-improved-intro} with the method of ``Taylor--Wiles primes'' (in the refined version of Greither \cite{Greither-Brumer} to strengthen further Theorem 1.9, by replacing the hypothesis $S_p(F)\subseteq S$ with the weaker $S_p(F)\cap S_{ram}(H/F)\subseteq S$, and prove a conjecture of Burns--Kurihara--Sano, as a consequence. Namely, we obtain the following. (See Theorem \ref{BKS-Conj} below.)

\begin{theorem}\label{BKS-Conj-intro}The following hold.
\begin{enumerate}
\item If $p$ is an odd prime and $S$ and $T$ satisfy hypotheses
$$ S_{\infty}(F)\cup (S_p(F)\cap S_{ram}(H/F))\subseteq S,\quad  S_{ram}(H/F)\subseteq S\cup T,\qquad  {\rm Hyp}_p^T(H/F),$$
then the $\Bbb Z_p[G]^-$--module $Sel_S^T(H)_p^-$ is quadratically presented and we have
$${\rm Fitt}_{\mathbb{Z}_p[G]^-}(Sel_S^T(H)_p^-)=(\Theta_{S, H/F}^T(0)).$$
\item {\bf (The Burns--Kurihara--Sano Conjecture.)}  If $S$ and $T$ satisfy the above hypotheses for all odd primes $p$ (e.g. if $S_{ram}(H/F)\subseteq S$ and $H_T^\times$ has no coprime--to--$2$ torsion), then we have
$${\rm Fitt}_{\mathbb{Z}[G]^-}(Sel_S^T(H)^-)=(\Theta_{S, H/F}^T(0)).$$
\end{enumerate}
\end{theorem}

We conclude this introduction by making the reader aware that an unconditional proof of the Refined Coates--Sinnott Conjecture \eqref{strong-CS} was also given in the recent work of Johnston and Nickel \cite{Johnston-Nickel-EMC}, as a consequence of their proof in loc.cit. of an unconditional Equivariant Main Conjecture whose flavour is different from that of our Theorem \ref{EMC-intro}. Our work and methods are independent and different from those in loc.cit.

\section{Notations}
 \subsection{Class--groups and units} Let $H/F$ be a Galois extension of number fields of Galois group $G$. Let  $\mathcal O_H$ be the ring of algebraic integers in $H$. For any finite place $w$ of $H$ (i.e. any maximal ideal of $\mathcal O_H$), we let ${\rm ord}_w:H^\times\to\Bbb Z$ denote the valuation associated to $w$, normalized so that ${\rm ord}_w(H^\times)=\Bbb Z$. If $S$ and $T$ are finite, disjoint sets of places of $F$, such that the set of infinite places $S_\infty(F)$ is contained in $S$, we let $S_H$ and $T_H$ be the sets of places in $H$ sitting above places in $S$ and $T$, respectively. For simplicity of notation and if no confusion arises, we use $S$, $T$ and $S_\infty$ to mean $S_H$, $T_H$ and $S_\infty(H)$, respectively.
Let
$$\mathcal O_{H,S,T}^{\times}:= \{x\in H^\times ; {\rm ord}_w(x)= 0, \text{ for all }w\not\in S_H, \quad {\rm ord}_{w}(x-1)>0, \text{ for all } w\in T_H\}.$$
This is a subgroup of finite index of the group of $S$--units $\mathcal O_{H,S}^\times$ of $H$. We let
$$H_T^{\times}:= \{x\in H^{\times} ; {\rm ord}_{w}(x-1)>0,\ \text{for all } w\in T_H\}.$$
In what follows,  $Cl^T(H)$ denotes the ray--class group of $H$ of conductor $(\prod_{w\in T_H}w)$ and $Cl^T_S(H)$ denotes the quotient $Cl^T(H)$ by the subgroup generated by the classes of ideals in $S_H\setminus S_\infty(H)$. We let $Y_S(H)$ and $Y_{\overline S}(H)$ denote the free $\Bbb Z$--modules of divisors (with integral coefficients) of $H$, supported at primes in $S_H\setminus S_\infty(H)$ and outside of $S_H$, respectively.

From the definitions, we have an obvious canonical exact sequence of $\Bbb Z[G]$--modules
\begin{equation}\label{class-group-sequence} \begin{tikzcd}
0\arrow{r} &\mathcal O_{H, S,T}^\times\arrow{r} & H_T^{\times}\arrow{r}\arrow{r}{{\text div}_{\overline{S\cup T}}} & Y_{\overline{S\cup T}}(H)\arrow{r} & Cl^T_S(H)\arrow{r} &  0.
\end{tikzcd}
\end{equation}
Above, the map ${\text div}_{\overline{S\cup T}}(\ast):=\sum_{w\notin S_H\cup T_H}{\rm ord}_w(\ast)\cdot w$ is the usual ($S\cup T$--depleted) divisor map and the right--most non--zero map is the divisor--class map.\\

\subsection{$L$--functions} Assume further that $G$ is abelian. For a place $v$ of $F$, we let $G_v$ and $I_v$ denote its decomposition and inertia groups in $G$, respectively and fix $\sigma_v\in G_v$ a Frobenius element for $v$.
The idempotent associated to the trivial character of $I_v$ in $\Bbb Q[I_v]$ is given by
$$e_v:=\frac{1}{|I_v|} N_{I_v}:=\frac{1}{|I_v|}\sum_{\sigma\in I_v}\sigma$$
Then $e_v\sigma_v^{-1}\in \Bbb Q[G]$ is independent of the choice of $\sigma_v$.
As in \cite{Dasgupta-Kakde}, the $\Bbb C[G]$--valued ($G$--equivariant) $L$--function associated to $(H/F, S, T)$ of complex variable $s$ is given by the meromorphic continuation to the entire complex plane of the following holomorphic function
$$\Theta_{S, H/F}^T(s):=\prod_{v\not\in S}(1-e_v\sigma_v^{-1}\cdot Nv^{-s})^{-1}\cdot\prod_{v\in T}(1-e_v\sigma_v^{-1}\cdot Nv^{1-s}), \qquad \text{Re}(s)>0.$$
The resulting meromorphic continuation (also denoted by $\Theta_{S, H/F}^T(s)$) is holomorphic on $\Bbb C\setminus\{1\}$. It is easily seen that if $H'/F$ is abelian, finite, of Galois group $G'$ and $H\subseteq H'$, then
\begin{equation}\label{projection-Theta}\Theta_{S, H/F}^T(s)=\pi_{G'\to G}(\Theta_{S, H'/F}^T(s)),\qquad \text{ for all }s\in\Bbb C\setminus\{1\},\end{equation}
where $\pi_{G'\to G}:\Bbb C[G']\to\Bbb C[G]$ is the usual $\Bbb C$--algebra morphism induced by Galois restriction.\\

Classical results proved independently by  Klingen--Siegel (\cite{Klingen} and \cite{Siegel}) and Shintani \cite{Shintani} show that
\begin{equation}\label{Klingen-Siegel}\Theta_{S, H/F}^T(1-n)\in\Bbb Q[G], \qquad \text{ for all } n\in\Bbb Z_{\geq 1}.\end{equation}
Further, if $T$ is chosen so that $H_T^\times$ is $\Bbb Z$--torsion free (so, necessarily $T\ne\emptyset$) and $S$ contains the ramification locus $S_{ram}(H/F)$ of $H/F$ then, another classical result, due independently to  Pi. Cassou--Nogu\`es \cite{Nogues} and Deligne--Ribet \cite{Deligne-Ribet}, shows that
\begin{equation}\label{Deligne-Ribet}\Theta^T_{S,H/F}(1-n)\in\Bbb Z[G], \qquad \text{ for all }n\in\Bbb Z_{\geq 1}.\end{equation}

Let $p>2$ be a prime and $S_p(F)$ be the set of $p$--adic primes in $F$. A $p$--adic improvement on \eqref{Deligne-Ribet} for $n=1$, due essentially to Kurihara, is the following.
\begin{lemma}\label{Deligne-Ribet-Kurihara} Assume that the non--empty, disjoint sets of primes $S$ and $T$ in $F$ satisfy the following set of hypotheses, for a prime $p>2$.
$${\rm Hyp}(H/F)_p:\, \begin{cases}
S_{ram}(H/F)\subseteq S\cup T;\\
(S_{ram}(H/F)\cap S_p(F))\cup S_\infty(F)\subseteq S;\\
\text{$T\not\subseteq S_{ram}(H/F)$ and $(H_{T\setminus S_{ram}(H/F)})^\times$ has no $p$--torsion.}
\end{cases}$$
Then, we have  $\Theta_S^T(H/F):=\Theta_{S,H/F}^T(0)\in\Bbb Z_p[G]$
\end{lemma}
\begin{proof} See the proof of Lemma 3.4 and Remark 3.6 in \cite{Dasgupta-Kakde}. Note that the third hypothesis above is satisfied if $T\not\subseteq S_{ram}(H/F)\cup S_p(F).$ This is an equivalence if the order $p$ roots of unity are in $H$. \end{proof}

Now, assume further that $H$ is CM and $F$ is totally real, and let $j$ be the complex conjugation automorphism of $H$, viewed as an element in $G$. We call a $\Bbb C$--valued character $\chi$ of $G$ even if $\chi(j)=+1$ and odd if $\chi(j)=-1$.   The following is a well--known consequence of the functional equation of $\Theta_{S,H/F}^T(s)$ and the fact that
$\Theta_{S,H/F}^T(s)\in\Bbb C[G]^\times$, for ${\rm Re}(s)>1$.
\begin{lemma}\label{vanishing-of-Theta} With notations as above, the following hold.
\begin{enumerate}
\item If $\chi$ is an even character of $G$, then $\chi(\Theta_S^T(H/F))=0$.
\item If $\chi$ is an odd character of $G$, then $\chi(\Theta_S^T(H/F))=0$ if and only if there exists a prime $v\in S\setminus S_\infty(F)$, such that $\chi(G_v)=\{1\}$.
\end{enumerate}
\end{lemma}
\begin{proof} See \cite{Tate}. \end{proof}
\vskip .2in

\subsection{Modules over commutative rings.}  If $G$ is a (pro)finite abelian group and $M$ is a $\Bbb Z[G]$--module, we endow its $\Bbb Z$--dual and Pontrjagin dual
$$M^\ast:={\rm Hom}_{\Bbb Z}(M, \Bbb Z),\qquad M^\vee:={\rm Hom}_{\Bbb Z}(M, \Bbb Q/\Bbb Z)$$
with the {\it covariant} $G$--actions $gf(x)=f(gx)$, for $g\in G$. Further, if $G$ contains a canonical element $j$ of order $2$ (usually a compex conjugation automorphism), then we let
$$M^-:=\frac{1}{2}(1-j)\cdot (M\otimes_{\Bbb Z}\Bbb Z[1/2]).$$
This is the $(-1)$--eigenspace of $M\otimes_{\Bbb Z}\Bbb Z[1/2]$ under the action of $j$. It has a natural structure of $\Bbb Z[1/2][G]^-$--module. If $p$ is an odd prime, then we let
$$M_p:=M\otimes_{\Bbb Z}\Bbb Z_p, \qquad M_p^-:=\frac{1}{2}(1-j)M_p=M^-\otimes_{\Bbb Z[1/2]}\Bbb Z_p,$$
viewed as modules over $\Bbb Z_p[G]$ and $\Bbb Z_p[G]^-$, respectively. Obviously, $(\ast\to \ast^-)$ and $(\ast\to \ast^-_p)$ are exact functors. Tacitly, we identify the rings $\Bbb Z_p[G]/(1+j)$ and $\Bbb Z_p[G]^-$, respectively  $\Bbb Z[1/2][G]/(1+j)$ and $\Bbb Z[1/2][G]^-$, via the isomorphism given by multiplication with $\frac{1}{2}(1-j)$.

Slightly abusively, if $N$ is a $\Bbb Z_p[[G]]$--module, we use the same notations as above
$$N^\ast:={\rm Hom}_{\Bbb Z_p}(N, \Bbb Z_p),\qquad N^\vee:={\rm Hom}_{\Bbb Z_p}(N, \Bbb Q_p/\Bbb Z_p)$$
to denote its $\Bbb Z_p$--dual and Pontrjagin dual, always endowed with the covariant $G$--action.
\\

If $\phi: R\to R'$ is a morphism of commutative rings and $M$ is an $R$--module, we let
$$M_{R'}:=M\otimes_R R'.$$
We remind the reader that if the $R$--module $M$ is finitely presented, then its $0$--th $R$--Fitting ideal ${\rm Fitt}_R(M)$ is well defined and we have an equality
\begin{equation}\label{base-change-Fitt}\phi({\rm Fitt}_R(M))R'={\rm Fitt}_{R'}(M\otimes_R R').\end{equation}
Also, if $M$ and $N$ are finitely presented $R$--modules and $\psi: M\twoheadrightarrow N$ is a surjective $R$--linear map, then we have an inclusion of $R$--ideals
\begin{equation}\label{inclusion-Fitt}{\rm Fitt}_R(M)\subseteq {\rm Fitt}_R(N).\end{equation}
\begin{definition}A finitely presented $R$--module $M$ is called quadratically presented if there is a $k\geq 1$ and an exact sequence of $R$--modules
$$R^k\longrightarrow R^k\longrightarrow M\longrightarrow 0. $$
Once a quadratic presentation of $M$ is given, we refer to $k$ as its rank.
\end{definition}
 The following results on Fitting ideals will be very useful in our considerations.
\begin{lemma}[Johnston--Nickel]\label{Fitt-ses}
If we have a  short exact sequence
$$0\longrightarrow A\longrightarrow B\longrightarrow C\longrightarrow 0$$
of finitely presented $R$--modules and $C$ is quadratically presented, then
$${\rm Fitt}_R(B)={\rm Fitt}_R(A)\cdot{\rm Fitt}_R(C).$$
\end{lemma}
\begin{proof}
See Lemma 2.13 in \cite{Johnston-Nickel}.
\end{proof}
\begin{lemma}[Burns--Greither]\label{Burns-Greither}
Assume that $R$ is a finitely generated $\Bbb Z_p$--algebra, which is $\Bbb Z_p$--free and relatively Gorenstein over $\Bbb Z_p$, for some prime $p$. Then, for any exact sequence of finite $R$--modules
$$0\to A\to P\to P'\to A'\to 0,$$
such that ${\rm pd}_R P\leq 1$ and ${\rm pd}_R P'\leq 1$, we have
$${\rm Fitt}_R(A^\vee)\cdot{\rm Fitt}_R(P')={\rm Fitt}_R(A')\cdot{\rm Fitt}_R(P).$$
\end{lemma}
\begin{proof} See Lemma 5, p.179 in \cite{Burns-Greither}. See also Proposition 7.3 in \cite{Greither-Popescu}.
\end{proof}
For additional basic properties of Fitting ideals needed throughout, the reader can consult \cite{Greither-Popescu} and \cite{Dasgupta-Kakde}.

\section{Selmer modules and class--groups}\label{Selmer-section}

\subsection{Selmer modules and class--groups in finite extensions}
Let $H/F$ be a finite abelian extension of number fields of Galois group $G$, with $F$ totally real and $H$ a CM field. Let $S$ and $T$ be finite, disjoint sets of places of $F$, such that $S_{\infty}\subseteq S$.
We recall the definition of the Selmer module defined by Burns-Kurihara-Sano in \cite{BKS} and studied further by Burns in \cite{Burns}. The reader should also consult \S3.1 in \cite{Dasgupta-Kakde}.  By taking $\Bbb Z$--duals in exact sequence \eqref{class-group-sequence}, we obtain a canonical, injective morphism of $\Bbb Z[G]$--modules
$$Y_{\overline{S\cup T}}(H)^\ast \hookrightarrow{} (H_T^{\times})^\ast.$$
If one identifies $\prod_{w\notin S_H\cup T_H}\mathbb{Z}\simeq Y_{\overline{S\cup T}}(H)^\ast$ in the obvious manner, the above injection sends the tuple $(x_w)_w\in\prod_{w}\Bbb Z$ to the homomorphism $\ast\to \sum_w x_w {\rm ord}_w(\ast)$.

\begin{definition} The Selmer $\Bbb Z[G]$--module for the data $(H/F, S, T)$ is given by
$${\rm Sel}_S^T(H):=(H_T^{\times})^\ast/Y_{\overline{S\cup T}}(H)^\ast \simeq (H_T^{\times})^\ast/\prod_{w\notin S_H\cup T_H}\mathbb{Z},$$
where the isomorphism is given by the identification described above.
\end{definition}

\begin{remark}\label{larger-S-remark} Assume that $S'$ is a finite set of primes in $F$, containing $S$ and disjoint from $T$. Then, since $Y_{\overline{S\cup T}}=Y_{S'\setminus S}\oplus Y_{\overline{S'\cup T}}$, then we obtain a short exact sequence of
$\Bbb Z[G]$--modules
\begin{equation}\label{enlarge-S-Selmer-sequence}
0\to Y_{S'\setminus S}(H)^\ast\to Sel_{S'}^T(H)\to Sel_S^T(H)\to 0.
\end{equation}
\end{remark}
The following is a more general and explicit version of Lemma 3.1 in [3], providing a link between class--groups and Selmer groups.
\begin{lemma}\label{class group=Selmer Module} For $(H/F, S, T)$ as above, the following hold.
 \begin{enumerate}
\item[(1)] We have a canonical exact sequence of $\Bbb Z[G]$--modules
\[ \begin{tikzcd}
0\arrow{r} & Cl^T_S(H)^{\vee}\arrow{r}{i} & Sel_{S}^T(H)\arrow{r}{j}  & (O_{H,T,S}^{\times})^\ast \arrow{r}  &  0,
\end{tikzcd}
\]
where $i$ sends $\phi\in Cl^{T}_S(H)^\vee$ to the coset in $Sel_{S}^T(H)$ of $[a\mapsto \sum_{w\notin T\cup S}\widetilde{\phi([w])}ord_w(a)]$ (viewed as an element in $(H_T^\times)^\ast$),  $\Tilde{x}\in\mathbb{Q}$ denotes any lift of $x\in\mathbb{Q}/\mathbb{Z}$ to $\Bbb Q$, $j$ is the restriction map, and $[w]$ is the ideal--class of $w$ in $Cl^T_S(H)$.

\item[(2)] In particular, if $S=S_\infty$, then $i$ induces an isomorphism of $\Bbb Z[G]^-$--modules
$$i: Cl^{T}(H)^{\vee,-}\simeq {\rm Sel}_{S_\infty}^T(H)^-,$$
explicitely described in (1) above.
\end{enumerate}
\end{lemma}
\begin{proof}
We apply the exact functor $\ast\to {\rm Hom}_{\mathbb{Z}}(\ast,\mathbb{Q})$ to exact sequence \eqref{class-group-sequence}. Since $Cl^T_S(H)$ is finite, we obtain the following exact sequence of $\Bbb Z[G]$--modules.

\[ \begin{tikzcd}
0\arrow{r} & {\rm Hom}_{\Bbb Z} (Y_{\overline{S\cup T}}, \Bbb Q)\arrow{r}  & {\rm Hom}_{\Bbb Z}(H_T^{\times}, \Bbb Q)\arrow{r} & {\rm Hom}_{\Bbb Z}(O_{H,T,S}^{\times}, \Bbb Q)\arrow{r}  &  0
\end{tikzcd}
\]
Now, by taking Pontryagin duals in \eqref{class-group-sequence}, we obtain an exact sequence of $\Bbb Z[G]$--modules.

\[ \begin{tikzcd}
0\arrow{r} & Cl^T(H)^{\vee}\arrow{r} & (Y_{\overline{S\cup T}})^\vee \arrow{r}  & (H_T^{\times})^\vee\arrow{r} & (O_{H,T,S}^{\times})^\vee\arrow{r}  &  0
\end{tikzcd}
\]
Now consider the following commutative diagram of $\Bbb Z[G]$--modules with exact rows.

\[ \begin{tikzcd}
& 0\arrow[swap]{d} & 0\arrow[swap]{d} & Cl^T_S(H)^{\vee} \arrow[hook]{d}\\%
0\arrow{r} & (Y_{\overline{S\cup T}})^\ast \arrow{r}\arrow[swap]{d} &{\rm Hom}_{\Bbb Z} (Y_{\overline{S\cup T}}, \Bbb Q)  \arrow{r}\arrow{d} & (Y_{\overline{S\cup T}})^\vee \arrow{r}\arrow{d}& 0\\%
0\arrow{r}& (H_T^{\times})^\ast \arrow{r}\arrow[two heads]{d} & {\rm Hom}_{\mathbb{Z}}(H_T^{\times},\mathbb{Q})\arrow{r}\arrow[two heads]{d} & (H_T^{\times})^\vee \arrow[two heads]{d} & \\%
& Sel_{S}^T(H)  & {\rm Hom}_{\mathbb{Z}}(O_{H,T,S}^{\times},\mathbb{Q}) & (O_{H,T,S}^{\times})^\vee
\end{tikzcd}
\]
Now, the snake lemma leads to the following exact sequence of $\Bbb Z[G]$--modules.

\[ \begin{tikzcd}
0\arrow{r} & Cl^T_S(H)^{\vee}\arrow{r}{i} & Sel_{S}^T(H)\arrow{r}{j}  & {\rm Hom}_{\mathbb{Z}}(O_{H,T,S}^{\times},\mathbb{Q})\arrow{r}{f} & (O_{H,T,S}^\times)^\vee
\end{tikzcd}\]
Clearly, $\ker f=(O_{H,T,S}^{\times})^\ast$. So we have the following short exact sequence.
\[ \begin{tikzcd}
0\arrow{r} & Cl^T_S(H)^{\vee}\arrow{r}{i} & Sel_{S}^T(H)\arrow{r} {j} & (O_{H,T,S}^{\times})^\ast \arrow{r}  &  0
\end{tikzcd}
\]
where $j$ is the restriction map. Now let us find the isomorphism $i$ explicitly. If $\phi\in Cl^T_S(H)^{\vee}$ a diagram chase corresponding to the above snake lemma set--up gives the following.

\[ \begin{tikzcd}
& & \phi \arrow[mapsto]{d}\\%
& (\widetilde{\phi([w])})_w \arrow[swap, mapsto]{r}\arrow[mapsto]{d} & (\phi([w]))_w \\%
(a\mapsto \sum_{w\notin T\cup S_\infty}\widetilde{\phi([w])}ord_w(a))\arrow[mapsto]{r}\arrow[mapsto]{d} & (a\mapsto \sum_{w\notin T\cup S_\infty}\widetilde{\phi([w])}ord_w(a))\\%
{[a\mapsto \sum_{w\notin T\cup\ S_{\infty}}\widetilde{\phi([w])}ord_w(a)]}
\end{tikzcd}
\]
This completes the proof of part (1).

(2) Now, since $H$ is a CM field, the group $(O_{H,T,S_{\infty}}^{\times})^-$ is finite. Hence, $(O_{H,T,S_{\infty}}^{\times})^{\ast,\, -}=0$. Therefore, by taking the minus part of the short exact sequence in part (1), we have the desired isomorphism.
This completes the proof of part (2).
\end{proof}

\begin{remark}\label{large-S-remark} Assume that $S'$ is a finite set of places in $F$, disjoint from $T$, containing $S$ and sufficiently large so that $Cl_S^T(H)$ is generated by ideals in $S'_H\setminus S_H$. Then, we have the following obvious analogue of exact sequence \eqref{class-group-sequence}.
$$ \begin{tikzcd}
0\arrow{r} &\mathcal O_{H, S,T}^\times\arrow{r} & \mathcal O_{H, S', T}^\times\arrow{r}\arrow{r}{{\text div}_{S'\setminus S}} & Y_{S'\setminus S}(H)\arrow{r} & Cl^T_S(H)\arrow{r} &  0.
\end{tikzcd}
$$
When combined with the proof of part (1) of the above Lemma, this leads to an isomorphism
\begin{equation}\label{large-S-iso}
Sel_S^T(H)\simeq\frac{ (O_{H, S', T}^\times)^\ast}{(Y_{S'\setminus S}(H))^\ast},\end{equation}
induced by the restriction map $(H_T^\times)^\ast\to (O_{H, S', T}^\times)^\ast.$
\end{remark}
\medskip

Let $H\subseteq H'$, where $H'$ is a CM abelian extension of $F$, of Galois group $G'$. The natural morphisms  $H_{T}^{\times}\hookrightarrow {H'_{T}}^{\times}$, $Y_{\overline{S\cup T}}(H)\hookrightarrow Y_{\overline{S\cup T}}(H')$,   and $Cl^T_S (H) \xrightarrow{} Cl^T_S(H')$ induce $G'$--equivariant restriction maps
$$res : Sel_S^T(H')\xrightarrow{} Sel_S^T(H), \qquad \widetilde{res} : Cl^T_S(H)^{\vee} \xrightarrow{ } Cl^T_S(H')^{\vee}.$$

\begin{lemma}\label{commutative-square}
We have the following commutative diagram of $G'$-modules
\[ \begin{tikzcd}
Cl^{T}_S(H')^{\vee} \arrow{r}{i'} \arrow{d}{\widetilde{res}} & Sel_{S}^T(H') \arrow{d}{res} \\%
Cl^{T}_S(H)^{\vee} \arrow{r}{i} &  Sel_{S}^T(H)
\end{tikzcd}
\]
where the horizontal maps are the injective maps $i$ given in Lemma \ref{class group=Selmer Module} (1).
\end{lemma}
\begin{proof} Let $\phi\in Cl^{T}_S(H')^{\vee}$. Then, if $v$ and $w$ denote primes in $H$ and $H'$, respectively, Lemma \ref{class group=Selmer Module} gives the following.
 $$res(i'(\phi))=res([a\mapsto \sum_{w\notin T\cup\ S_{\infty}}\widetilde{\phi([w])}ord_w(a)])=[a\mapsto \sum_{v\notin T\cup\ S_{\infty}}\sum_{w|v}\widetilde{\phi([w])}ord_w(a)].$$
On the other hand, $i(\widetilde{res}(\phi))=[a\mapsto \sum_{v\notin T\cup\ S_{\infty}}\widetilde{\phi([v])}ord_v(a)].$
However, if $e$ is the ramification index of $v$ in $H'/H$, since $ord_v(\ast)=e\cdot ord_w(\ast)$, we have
$$\sum_{w|v}\widetilde{\phi([w])}ord_w(a)=\sum_{w|v}\widetilde{\phi([w^e])}ord_v(a)=\widetilde{\phi(\prod_{w/v}[w^e])} ord_v(a)=\widetilde{ \phi([v])}ord_v(a),$$
making the diagram above commutative.
\end{proof}\
We conclude this section with a lemma on the surjectivity of the Selmer restriction.
\begin{lemma}\label{surjective-res-Selmer}
With notations as above, assume that $(H'_T)^{\times}$ has no $p$--torsion, for some prime $p$. Then the restriction map
$res : Sel_S^T(H')_p\xrightarrow[]{} Sel_S^T(H)_p$
is surjective.
\end{lemma}
\begin{proof} Let $S'$ be a finite set of primes in $F$, containing $S$, disjoint from $T$ and large enough so that $Cl_S^T(H)$ and $Cl_S^T(H')$ are generated by classes of $S'$--supported ideals. Then, Remark \ref{large-S-remark} above shows that it suffices to show that
the restriction map
$$res: (O_{H', S', T}^\times)^\ast_p\to  (O_{H, S', T}^\times)^\ast_p$$
is surjective. By taking $\Bbb Z$--duals in the exact sequence of finitely generated $\Bbb Z$--modules
$$0\xrightarrow[]{} O_{H,S',T}^{\times}\xrightarrow[]{} O_{H',S',T}^{\times} \xrightarrow[]{} O_{H',S',T}^{\times}/O_{H,S',T}^{\times}\xrightarrow[]{} 0,$$
and then tensoring with $\Bbb Z_p$,
this amounts to showing that $O_{H',S',T}^{\times}/O_{H,S',T}^{\times}$ has no $p$--torsion, which we show below.

Let $x\in O_{H',S',T}^{\times}$ such that $x^p\in O_{H,S',T}^{\times}$. Then, we have $1=(x^p)^{g-1}=(x^{g-1})^p$, for all $g\in Gal(H'/H)$. But, since $(H'_T)^{\times}$ is $p$--torsion free, so is $O_{H',S',T}^{\times}$. Therefore, we have $x^{g-1}=1$ for all $g\in Gal(H'/H)$. Hence, $x\in O_{H,S',T}^{\times}$. Consequently,  $O_{H',S',T}^{\times}/O_{H,S',T}^{\times}$ has no $p$--torsion, which concludes the proof.
\end{proof}

\begin{corollary}\label{maps-Fitt-Selmer}
Under the assumptions of Lemma \ref{surjective-res-Selmer}, if $\pi:\Bbb Z_p[G']\to\Bbb Z_p[G]$ is the $\Bbb Z_p$--algebra morphism induced by Galois restriction, then we have
$$\pi({\rm Fitt}_{\Bbb Z_p[G']} Sel_S^T(H')_p)\subseteq {\rm Fitt}_{\Bbb Z_p[G]}Sel_S^T(H)_p.$$
\end{corollary}
\begin{proof}
Note that the surjective restriction map $res : Sel_S^T(H')_p\twoheadrightarrow Sel_S^T(H)_p$  (see Lemma \ref{surjective-res-Selmer}) induces a surjective morphism of $\Bbb Z_p[G]$--modules
 $$\widehat{res} : (Sel_S^T(H')_p)_{G(H'/H)}\twoheadrightarrow Sel_S^T(H)_p.$$
 Consequently, if we apply properties \eqref{base-change-Fitt} and \eqref{inclusion-Fitt} of Fitting ideals, we obtain
 $$\pi({\rm Fitt}_{\Bbb Z_p[G']} Sel_S^T(H')_p)={\rm Fitt}_{\Bbb Z_p[G]} (Sel_S^T(H')_p)_{G(H'/H)}\subseteq {\rm Fitt}_{\Bbb Z_p[G]}Sel_S^T(H)_p,$$
 which concludes the proof.
\end{proof}

\subsection{Selmer modules and class--groups in Iwasawa towers} Let $H/F$, $S$ and $T$ be as in the previous section. We denote by $H_{\infty}$ the cyclotomic $\mathbb{Z}_p$--extension of $H$. Let $H_n$ be its $n$--th layer,  characterized by the usual equality $[H_n:H]=p^n$. Then, $H_n$ and $H_\infty$ are also CM fields and the complex conjugation automorphim of $H_n$ is the restriction of that of $H_\infty$. Let $G_n:=Gal(H_n/F)$, for all $n\in \mathbb{N}$ and let $\mathcal{G}:=Gal(H_\infty/F).$

We consider the $G_{n+1}$--equivariant restriction maps corresponding to $H_{n+1}/H_n$
$$res_n : Sel_S^T(H_{n+1})\xrightarrow{} Sel_S^T(H_n), \qquad \widetilde{res_n} : Cl^T_S(H_{n+1})^{\vee} \xrightarrow{ } Cl^T_S(H_{n})^{\vee},$$
and define $p$--primary Selmer groups at the top of the Iwasawa tower $H_\infty/H$ as follows.

\begin{definition} We define the $\Bbb Z_p[[\mathcal G]]$--module
$$Sel_S^T(H_{\infty})_p=\underset{n}{\varprojlim}\text{ }Sel_S^T(H_n)_p,$$
where the projective limit is taken with respect to the maps $res_n$ defined above.
\end{definition}
The following Lemma shows that  $Sel_S^T(H_{\infty})_p$ only depends on the data $(H_{\infty}/F ,S,T,p)$, not on the intermediate CM field $H$.
\begin{lemma}\label{independence}
Let $H, H'$ be two CM, abelian extensions of the totally real number field $F$, such that  $H_\infty=H'_\infty$. Then, the following hold.
\begin{enumerate}
\item We have an isomorphism of $\Bbb Z_p[[\mathcal G]]$--modules
$$\varprojlim_n\, Sel_S^T(H_n)_p\cong\varprojlim_n\, Sel_S^T(H_n')_p,$$
compatible (in the obvious way) with the restriction maps.
\item We have an equality of $\Bbb Z_p[[\mathcal G]]$--ideals
$$\varprojlim_n\, {\rm Fitt}_{\Bbb Z_p[G_n]}\, Sel_S^T(H_n)_p\cong\varprojlim_n\, {\rm Fitt}_{\Bbb Z_p[G'_n]} Sel_S^T(H_n')_p,$$
where $G'_n:=G(H'_n/F)$, for all $n$, and the transition maps at the level of Fitting ideals are induced by the Galois restriction maps, as in Corollary \ref{maps-Fitt-Selmer}.
\end{enumerate}
\end{lemma}
\begin{proof} Obvious, left to the reader.
\end{proof}

\begin{definition}
We let $A_n^T:=Cl^T(H_n)_p$, for all $n$, and define the $\Bbb Z_p[[\mathcal G]]$--module
$$A_{\infty}^T:=\underset{n}{\varinjlim}\text{ }A_n^T,$$
where the injective limit is taken with respect to the natural maps at the level of class--groups.
\end{definition}
\noindent Obviously, we have an isomorphism of  $\mathbb{Z}_p[[\mathcal{G}]]$--modules
$$A_{\infty}^{T,\vee}:=\underset{n}{\varprojlim}\text{ }A_n^{T,\vee},$$
where the projective limit is taken with respect to the maps $\widetilde{res_n}$.\\\\

In order to simplify notation, in what follows we let $Y_{n, R}:=Y_R(H_n)$, for any $n$ and any finite set of primes $R$ in $F$ which contains $S_\infty$. We let
$$Y_{\infty, R}=\varinjlim_n Y_{n, R},$$
where the injective limit is taken with respect to the usual injective maps at divisor level. Since the primes above $p$ in $F$ (whose set is denoted by $S_p$) are totally ramified
in $H_\infty/H$, for $n\gg 0$, we have an isomorphism of $\Bbb Z$--modules
$$ Y_{\infty, S_p}\simeq \bigoplus_{w\in S_p(H_\infty)}\Bbb Z[1/p].$$
In particular $Y_{\infty, S_p}$ is $p$--divisible and, consequently,  we have
\begin{equation}\label{divisor-dual-vanishing}\varprojlim_n Y_{n, S_p}^\ast\simeq Y_{\infty, S_p}^\ast=0, \qquad \varprojlim_n Y_{n, S}^\ast\simeq Y_{\infty, S}^\ast= Y_{\infty, S\setminus S_p}^\ast,\end{equation}
for any finite set $S$ of primes in $F$.

\begin{proposition}\label{SES-infinity}
We have the following short exact sequence of $\mathbb{Z}_p[[\mathcal{G}]]^-$--modules.
\[ \begin{tikzcd}
0\arrow{r} & (Y_{\infty,S\setminus S_p}^\ast \otimes_{\mathbb{Z}}\mathbb{Z}_p)^-\arrow{r} & Sel_S^T(H_{\infty})_p^-\arrow{r}{\lambda} & A_{\infty}^{T,\vee,-}\arrow{r} & 0
\end{tikzcd}
\]
\end{proposition}
\begin{proof}
The (vertical) restriction maps lead to the following commutative diagrams with exact rows (see Remark \ref{larger-S-remark} and note that $(Y_{n, S}\otimes_{\Bbb Z}\Bbb Z_p)^-=(Y_{n, S\setminus S_\infty}\otimes_{\Bbb Z}\Bbb Z_p)^-$, for all $n$) in the category of $\Bbb Z_p[[\mathcal G]]$--modules, for all $n\geq 0$.
\[\begin{tikzcd}
0\arrow{r} & (Y_{n+1,S}^\ast\otimes_{\mathbb{Z}}\mathbb{Z}_p)^-\arrow{r}\arrow{d} & Sel_S^T(H_{n+1})^-_p\arrow{r}\arrow{d} & Sel_{S_{\infty}}^T(H_{n+1})^-_p\arrow{r}\arrow{d} & 0\\%
0\arrow{r} & (Y_{n,S}^\ast\otimes_{\mathbb{Z}}\mathbb{Z}_p)^-\arrow{r} & Sel_S^T(H_n)^-_p\arrow{r} & Sel_{S_{\infty}}^T(H_n)^-_p\arrow{r} & 0.
\end{tikzcd}
\]
Note that all the modules above are finitely generated $\Bbb Z_p$--modules, therefore compact in the $p$--adic topology, and that all the morphisms involved are continuous. Therefore, when taking the projective limit with respect to the restriction maps, one gets a short exact sequence of $\Bbb Z_p[[\mathcal G]]^-$--modules. When combined with \eqref{divisor-dual-vanishing} and the isomorphism in Lemma \ref{class group=Selmer Module}(2), that exact sequence is identical to the one in the statement of the theorem.
\end{proof}

\begin{corollary}\label{finitely-generated}
The  $\mathbb{Z}_p[[\mathcal{G}]]^-$--module $Sel_S^T(H_{\infty})_p^-$ is finitely generated and torsion. In particular, its Fitting ideal ${\rm Fitt}_{\Bbb Z_p[[\mathcal G]]^-} Sel_S^T(H_\infty)_p^-$ is well defined.
\end{corollary}
\begin{proof}
By the above proposition, we need to show that both $A_{\infty}^{T,\vee,-}$ and $(Y_{\infty,S\setminus S_p}^\ast \otimes_{\mathbb{Z}}\mathbb{Z}_p)^-$ are finitely generated and torsion $\mathbb{Z}_p[[\mathcal{G}]]^-$--module.

By \cite{Greither-Popescu} (see exact sequence (11) and proof of Lemma 2.8), there are non--negative integers $\delta$ and $\delta_T$ and  an exact sequence of $\mathbb{Z}_p[[\mathcal{G}]]^-$--modules:
$$0\xrightarrow[]{} (\mathbb{Q}_p/\mathbb{Z}_p)^{\delta}(1)\xrightarrow[]{}(\mathbb{Q}_p/\mathbb{Z}_p)^{\delta_T}(1)\xrightarrow[]{}A_{\infty}^{T,-}\xrightarrow[]{}A_{\infty}^{\emptyset,-}\xrightarrow[]{}0.$$
When taking Pontryagin duals, we obtain an exact sequence of $\Bbb Z_p[[\mathcal G]]^-$--modules
$$0\xrightarrow[]{}A_{\infty}^{\emptyset,\vee,-}\xrightarrow[]{} A_{\infty}^{T,\vee,-}\xrightarrow[]{}\mathbb{Z}_p^{\delta'}(1)\xrightarrow[]{}\mathbb{Z}_p^{\delta}(1)\xrightarrow[]{}0.$$
However, Kummer duality gives an isomorphism of $\Bbb Z_p[[\mathcal G]]^-$--modules
$$A_{\infty}^{\emptyset,\vee,-}\cong X_{S_p}^+(-1),$$
where $X_{S_p}=Gal(M_{S_p}^{ab}/H_{\infty})$ and  $M_{S_p}^{ab}$ is the maximal abelian pro-$p$ extension of $H_{\infty}$, unramified outside $S_p$. (See \cite{Washington}, Proposition 13.32.)
Therefore, since $X_{S_p}^+$ is finitely generated and torsion as a $\Bbb Z_p[[\mathcal G]]^-$--module (see \cite{Washington}, \S13.5), the last last displayed exact sequence and isomorphism show that $A_{\infty}^{T,\vee,-}$  is also a finitely generated, torsion $\Bbb Z_p[[\mathcal G]]^-$--module. \\

Now, if we let $\mathcal G_v$ denote the decompostion group of $v\in S\setminus (S_\infty\cup S_p)$ in $H_\infty/F$, it is easily seen that we have an isomorphism of $\Bbb Z_p[[\mathcal G]]^-$--modules
$$(Y_{\infty,S\setminus S_p}\otimes_{\mathbb{Z}}\mathbb{Z}_p)^-\simeq \bigoplus _{v\in S\setminus (S_\infty\cup S_p)}\Bbb Z_p[[\mathcal G/\mathcal G_v]]^-.$$
Since the groups $\mathcal G/\mathcal G_v$ are finite, for all $v$ as above, the module $(Y_{\infty,S\setminus S_p}\otimes_{\mathbb{Z}}\mathbb{Z}_p)^-$ is finitely generated over $\Bbb Z_p$. Consequently, its $\Bbb Z_p$--dual is also finitely generated over $\Bbb Z_p$, and therefore it is finitely generated and torsion over $\Bbb Z_p[[\mathcal G]]^-$.
\end{proof}

\begin{proposition}\label{limit-of-Fitt} If $H_{{\infty}, T}^\times$ has no $p$--torsion, we have an equality of $\Bbb Z_p[[\mathcal G]]^-$--ideals
$$\varprojlim_n\,  {\rm Fitt}_{\Bbb Z_p[G_n]^-} Sel_S^T(H_n)_p^- = {\rm Fitt}_{\Bbb Z_p[[\mathcal G]]^-} Sel_S^T(H_\infty)_p^-.$$
\end{proposition}
\begin{proof}
First, we know from \cite{Greither-Popescu} that there is a CM number field $H'$, such that
$$F\subseteq H'\subseteq H_\infty, \quad H'_\infty=H_\infty, \quad H'\cap F_\infty=F.$$
If we  let $G':=G(H'/F)$, we have $\mathcal G\simeq G'\times \Gamma$, where $\Gamma= G(F_\infty/F).$ Therefore, we have the following isomorphisms of topological $\Bbb Z_p$--algebras
$$\Bbb Z_p[[\mathcal G]]^-\simeq\Bbb Z_p[G']^-[[\Gamma]]\simeq\Bbb Z_p[[G']]^-[[T]].$$
Now, Lemma \ref{independence} allows us to replace $H$ by $H'$ in the equality we are trying to prove. If we let $X_n:=Sel_S^T(H'_n)_p^-$ and $X:=\varprojlim_n X_n=Sel_S^T(H_\infty)_p^-$
and recall that the transition maps $X_{n+1}\to X_n$ are surjective (see Lemma \ref{surjective-res-Selmer}) and that $X$ is finitely generated and torsion over $R=\Bbb Z_p[[\mathcal G]]^-$ (see Corollary \ref{finitely-generated}), then we can apply Theorem 2.1 in \cite{Greither-Kurihara}
to conclude that
$${\rm Fitt}_R(X)=\varprojlim_n\, {\rm Fitt}_{R_n}(X_n),$$
where $R_n:=\Bbb Z_p[G_n]^-$, for all $n$. This concludes the proof.
\end{proof}

On of the main results of this paper (the equivariant Iwasawa main conjecture) describes the Fitting ideal  ${\rm Fitt}_{\Bbb Z_p[[\mathcal G]]^-} Sel_S^T(H_\infty)_p^-$ in terms of certain equivariant $p$--adic $L$--functions, which will be defined and studied at the end of the next section.

\section{Some algebraic results}
In this section, we will prove some algebraic results needed for the rest of the paper. We also define the relevant equivariant $p$--adic $L$--functions and show that they are not zero divisors in the appropriate Iwasawa algebras.

 \begin{lemma}\label{finite-rank}
 Let $(A_n)_n$ be a projective system of $R$--modules, where $R$ is a PID. Suppose that the modules $A_n$ are $R$--free, of finite and bounded rank. Then, the $R$--module
$$A_\infty:=\underset{n}{\varprojlim}\text{ }A_n$$
is free of finite rank.
 \end{lemma}
\begin{proof}
For each $n$, we have the canonical projection map $\pi_n: A_{\infty} \xrightarrow[]{} A_n$. For all $n$, let $B_n:=Image(\pi_n)$. Clearly, $(B_n)_n$ forms a projective system of free $R$--modules of finite, bounded rank (as $R$--submodules of the $A_n$'s), whose transition maps are surjective  and
$$A_{\infty}=\underset{n}{\varprojlim}\text{ }B_n.$$
 Hence, we may assume from the very beginning that the transition maps $A_{n+1} \xrightarrow[]{} A_n$ are all surjective.
 Consequently, we have $rank_R(A_n)\leq rank_R(A_{n+1})$, for all $n$. Since the increasing sequence  $(rank_R(A_n))_n$ is bounded,  it is stationary for $n\gg 0$. Therefore, the transition maps give isomorphisms $A_{n+1}\cong A_n$, for $n\gg 0$.
Therefore, $A_{\infty} \cong A_k$, for $k\gg 0$. Hence, $A_{\infty}$ is a free $R$-module of finite rank, as desired.
 \end{proof}


Next, we fix a prime $p>2$ and take the cyclotomic $\Bbb Z_p$--extension $H_\infty/H$. We keep the notations of the previous subsection. In particular, $G_n:=G(H_n/F)$ and $\mathcal G:=G(H_\infty/F)$.
Note that if $G=G'\times G_p$, where $G_p$ is the Sylow $p$--subgroup of $G$, then we have group isomorphisms
$G_n\cong G'\times G_{p,n}$ where $G_{p,n}$ is the Sylow $p$--subgroup of $G_n$, for all $n$. We assume that the finite sets of primes $S$ and $T$ in $F$ satisfy the hypotehses in Lemma \ref{Deligne-Ribet-Kurihara}.

We have group isomorphisms
$$\mathcal{G}\cong G'\times G_{p,\infty}, \quad G_{p, \infty}\cong  t(G_{p, \infty})\times \Gamma, $$
where $G_{p,\infty}={\varprojlim}_n\text{ }G_{p,n}$, $t(G_{p,\infty})$ is the  torsion subgroup of $G_{p,\infty}$ (a finite group isomorphic via Galois restriction to asubgroup of $G_p$),  and $\Gamma\simeq \mathbb{Z}_p$. Consequently, we have continuous isomorphisms of  $\Bbb Z_p$--algebras
$$R_n:=\Bbb Z_p[G_n]\simeq \Bbb Z_p[G'][G_{p, n}], \quad R_\infty:=\Bbb Z_p[[\mathcal G]]\simeq \Bbb Z_p[G'][[G_{p,\infty}]]\simeq \Bbb Z_p[G'][t(G_{p,\infty})][[\Gamma]].$$

For a finite, abelian group $A$, we denote by $\widehat A$ its group of $\overline{\Bbb Q_p}$--valued characters. We let $[\widehat A]$ denote the set of equivalence classes of such characters  with respect to the equivalence relation
$\chi\sim \chi'$ if there exists $\sigma\in G_{\Bbb Q_p}$ such that $\chi'=\sigma\circ\chi$. For characters $\chi\in\widehat{G'}$, we let
$$R_n^\chi:=\Bbb Z_p[\chi][G_{p, n}], \qquad R_\infty^\chi:=\Bbb Z_p[\chi][[G_{p, \infty}]].$$
Then, we have topological ring isomorphisms given by the usual character evaluation maps
\begin{equation}\label{semilocal}R_n\cong \bigoplus_{[\chi]\in[\widehat{G'}]}R^{\chi}_n, \qquad R_\infty \cong \bigoplus_{[\chi]\in[\widehat{G'}]}R^{\chi}_{\infty}, \qquad R_\infty^\chi\cong \varprojlim_n R_n^\chi.\end{equation}

Now, we fix a character $\chi\in\widehat{G'}$ and let $\mathcal O_\chi=:\mathbb{Z}_p[\chi]$. This is a finite, unramified extension of $\mathbb{Z}_p$. Let us also fix a finite set of non--archimedean places $\Bbb S$ in $F$ and, for every $v\in \Bbb S$ let $G_{n, v}$ and $\mathcal G_v$ denote the decomposition groups of $v$ inside $G_n$, for all $n$, and $\mathcal G$, respectively. Character evaluation maps $\iota_n^\chi$ (injective) and  projections $\pi_n^\chi$ (surjective) give the following morphisms of $\mathcal O_\chi$--algebras
\[\begin{tikzcd}
f_n^{\chi, \Bbb S} :R_n^\chi \arrow{r}{\iota_n^\chi}& \bigoplus\limits_{\substack{\psi_n|_{G'}=\chi\\\psi_n\in \widehat{G_n}}} \mathcal{O}_\chi[\psi_n] \arrow{r}{\pi_n^{\chi, \Bbb S}}&\bigoplus\limits_{\substack{\psi_n|_{G'}=\chi\\ \psi_n(G_{n,v})\neq 1\\ \forall v\in \Bbb S}}\mathcal{O}_\chi[\psi_n]:=R_n^{\chi, \Bbb S}.
\end{tikzcd}\]

\begin{lemma} \label{whole-ring-at-infinity}
With notations as sbove, we let $R_\infty^{\chi, \Bbb S}:={\varprojlim}_n R_n^{\chi, \Bbb S}$, where the projective limit is taken with respect to the canonical transition maps. Then, the following
$$f_\infty^{\chi, \Bbb S}:=\varprojlim_n f_n^{\chi, \Bbb S}: R_{\infty}^{\chi}\longrightarrow R_\infty^{\chi, \Bbb S}$$
is an injective morphism of $\mathcal O_\chi$--algebras.
\end{lemma}
\begin{proof}
Since $\iota_n^\chi$ is injective, we can view $R_n^\chi$ as a subring of $\bigoplus\limits_{\substack{\psi|_{G'}=\chi\\ \psi\in \widehat{G_n}}} \mathcal{O}_\chi[\psi_n]$. Then
$$ker(f_n^{\chi, \Bbb S})\subseteq ker(\pi_n^{\chi, \Bbb S})=\bigoplus\limits_{\substack{\psi_n|_{G'}=\chi\\ \psi_n(G_{n,v}) = 1\\ \text{for some } v\in \Bbb S}}\mathcal{O}_\chi[\psi_n]. $$
Now, observe that the characters $\psi_n$ on the right factor through $G_n/G_{n,v_0}$, so via the surjectve group morphism  $\mathcal G/\mathcal G_{v_0}\to G_n/G_{n, v_0}$, they are characters of $\mathcal G/\mathcal G_{v_0}$, for some $v_0\in\Bbb S$.
Therefore, for all $n$, we have the following inequalities
\begin{equation}\label{bounded-rank-kernels} rank_{\mathcal O_\chi}(ker(f_n^{\chi, \Bbb S}))\leq rank_{\mathcal{O}_\chi}(\bigoplus\limits_{\substack{\psi_n|_{G'}=\chi\\ \psi_n(G_{n,v}) = 1\\ \text{for some } v\in \Bbb S}}\mathcal{O}_\chi[\psi_n])\leq rank_{\mathcal{O}_\chi}(\bigoplus_{v\in \Bbb S}\mathcal{O}_\chi[\mathcal G/\mathcal G_v])=K<\infty,\end{equation}
where the last inequality is due to the fact that $\mathcal G/\mathcal G_v$ is a finite group, for any non--archimedean prime $v$ of $F$. Note that the integer $K$ does not depend on $n$.
 Now,  observe that, for each $n$, we have the following commutative diagram of finitely generated $\mathcal O_\chi$--modules whose rows are exact.
\[\begin{tikzcd}
0 \arrow{r} & ker(f_{n+1}^{\chi, \Bbb S}) \arrow{r}\arrow{d} & R_{n+1}^\chi \arrow{r}{f_{n+1}^{\chi, \Bbb S}}\arrow{d} & R_{n+1}^{\chi, \Bbb S}\ar{d} \\%
0 \arrow{r} & ker(f_n^{\chi, \Bbb S}) \arrow{r} & R_n^\chi \arrow{r}{f_n^{\chi, \Bbb S}} & R_n^{\chi, \Bbb S}
\end{tikzcd}
\]
The middle vertical maps are Galois restrictions and all the other vertical maps are induced by them. Since all the maps involved are continuos in the $p$--adic topology, these are commutative diagrams in the category of topological, compact $\mathcal O_\chi$--modules. Therefore, after taking a projective limit, we obtain an exact sequence of compact $\mathcal O_\chi$--modules
\[\begin{tikzcd}
0 \arrow{r} & \underset{n}{\varprojlim}\text{  }ker(f_n^{\chi, \Bbb S}) \arrow{r} & R_{\infty}^{\chi} \arrow{r}{f_\infty^{\chi, \Bbb S}} & R_\infty^{\chi, \Bbb S}.
\end{tikzcd}
\]
Now, observe that $\mathcal O_\chi$ is a PID and the projective system of $\mathcal{O}_\chi$--modules $(ker(f_n^{\chi, \Bbb S}))_n$ satisfies the conditions of Lemma \ref{finite-rank}. So, we can
conclude that $ker(f_\infty^{\chi, \Bbb S})={\varprojlim}_n\, ker(f_n^{\chi, \Bbb S})$ is a free $\mathcal{O}_\chi$--module of finite rank.
However, since the maps $f_n^{\chi, \Bbb S}$ are ring morphisms, so is $f_\infty^{\chi, \Bbb S}$  and therefore $ker(f_\infty^{\chi, \Bbb S})$ is in fact an ideal in the ring
$$R^{\chi}_{\infty}=\mathcal O_\chi[[G_{p,\infty}]]\cong \mathcal O_\chi[t(G_{p,\infty})][[T]],$$
where the (non-canonical) isomorphism sends a fixed topological generator of $\Gamma$ to $(1+T)$, as usual.
It is easily seen that the only ideal in the ring $\mathcal O_\chi[G_p][[T]]$ which has finite $\mathcal O_\chi$--rank is the zero ideal. Therefore $ker(f_\infty^{\chi, \Bbb S})=0$, which concludes the proof.
\end{proof}
\vskip .2in

Now, we take two finite, disjoint sets $S$ and $T$ of primes in $F$, satisfying the hypotheses in Lemma \ref{Deligne-Ribet-Kurihara} for $H_n/F$, for all $n$. Note that, in particular $S_p(F)\subseteq S$, where $S_p(F)$ is the set of all $p$--adic primes in $F$.  Lemma \ref{Deligne-Ribet-Kurihara} and Lemma \ref{vanishing-of-Theta}(1) imply that $\Theta_S^T(H_n/F)\in \Bbb Z_p[G_n]^-$, for all $n$.
Further, equality \eqref{projection-Theta} implies that we have a well--defined element
$$\Theta_S^T(H_\infty/F):=\left(\Theta_S^T(H_n/F)\right)_n \in \Bbb Z_p[[\mathcal G]]^{-}.$$

\begin{definition}\label{p-adic-L-definition}
The element $\Theta_S^T(H_\infty/F)\in  \Bbb Z_p[[\mathcal G]]^{-}$ defined above is called the $S$--incomplete, $T$--smooth, $G$--equivariant $p$--adic $L$--function associated to the data $(H/F, S, T, p)$.
\end{definition}

\begin{proposition}\label{Theta-non-zero-divisor}
 Under the above hypotheses, the $p$--adic $L$--function $\Theta_S^T(H_\infty/F)$ is not a zero--divisor in $\Bbb Z_p[[\mathcal G]]^-$.
\end{proposition}
\begin{proof} Let $\Bbb S:=S\setminus S_\infty$. Note that $\Bbb S$ is a non-empty set of finite primes in $F$, as it contains $S_p(F)$. Lemma \ref{vanishing-of-Theta}(2) implies that, for all odd characters $\chi\in\widehat{G'}$, the elements
$f_n^{\chi, \Bbb S}(\Theta_S^T(H_n/F)^\chi)$ are not zero--divisors in $R_n^{\chi, \Bbb S}$, for all $n$. This shows that
$$f_\infty^{\chi, \Bbb S}(\Theta_S^T(H_\infty/F)^\chi)=(f_n^{\chi, \Bbb S}(\Theta_S^T(H_n/F)^\chi))_n$$
is not a zero--divisor in $R_\infty^{\chi, \Bbb S}$. However, since $f_\infty^{\chi, \Bbb S}$ is an injective ring morphism (see Lemma \ref{whole-ring-at-infinity}),  this implies that
$\Theta_S^T(H_\infty/F)^\chi$ is not a zero--divior in $R_\infty^\chi$, for all odd characters $\chi$ of $G'$. Since, by definition,  we  have
$$\Theta_S^T(H_\infty/F)=(\Theta_S^T(H_\infty/F)^\chi)_{[\chi]\, \text{odd}} \quad \text{ in }\quad\bigoplus_{[\chi]\,\text{odd}}R_\infty^\chi=\Bbb Z_p[[\mathcal G]]^-,$$
this implies the desired result.
\end{proof}
\medskip

Next, we prove a Lemma regarding ideals in projective limits of topological, compact rings. Typical examples of such systems are $(R_n)_n$ and $(R_n^\chi)_n$, defined above, endowed with their $p$--adic topologies and the usual transition maps given by
Galois restriction.

\begin{lemma}\label{proj-lim-ideals} Assume that $(S_n)_n$ is a projective system of topological, compact rings whose transition maps $\pi_n^{n+1}:S_{n+1}\to S_n$ are continuous.
Let $S_\infty:=\varprojlim_n S_n$.
Assume that $I_n$ is an ideal in $S_n$ of generators $(x_n^1, \dots, x_n^k)$ such that
$$\pi_n^{n+1}(x_{n+1}^i)=x_n^i, \qquad\text{ for all $i=1, \dots, k$ and all $n\geq 0$}.$$
Then the $S_\infty$--ideal $I_\infty:=\varprojlim_n I_n$
is generated by $(x_\infty^1:=(x_n^1)_n, \dots, x_\infty^k:=(x_n^k)_n)$.
\end{lemma}
\begin{proof} Consider the projective system of exact sequences of compact $S_n$--modules
$$ 0\longrightarrow\ker(\rho_n)\longrightarrow S_n^k \overset{\rho_n}\longrightarrow I_n\longrightarrow 0,$$
where $\rho_n((\alpha_i)_i):=\sum_i\alpha_i\cdot x_n^i$ and the transition maps are given by the $\pi_n^{n+1}$'s. Since all the modules above are compact, the projective limit
$$0\longrightarrow\varprojlim_n \ker(\rho_n)\longrightarrow S_\infty^k \overset{\rho_\infty}\longrightarrow I_\infty\longrightarrow 0,$$
stays exact. Clearly, $\rho_\infty ((\alpha_i)_i):=\sum_i\alpha_i\cdot x_\infty^i$, which concludes the proof.
\end{proof}
\medskip

Finally, we give a couple of technical results on Fitting ideals, needed in the later sections.
\begin{lemma}\label{Fitt-Y_S-dual}
Let $S$ be a finite set of non-archimedean primes in $F$. The following hold.
\begin{enumerate}
\item There is a canonical isomorphism of $\Bbb Z[G]$--modules
$$Y_S(H)^\sharp\simeq Y_S(H)^\ast,$$
where $Y_S(H)^\sharp$ is $Y_S(H)$ endowed with the $\Bbb Z[G]$--action induced by the $\Bbb Z$--linear involution
$\sharp:\Bbb Z[G]\simeq \Bbb Z[G]$ sending $g\mapsto g^{-1}$, for all $g\in G$.
\item We have equalities of $\Bbb Z[G]$--ideals
$${\rm Fitt}_{\Bbb Z[G]}Y_S(H)^\ast = {\rm Fitt}_{\Bbb Z[G]}Y_S(H)=\prod_{v\in S}I_{G_v}=\prod_{v\in S} (1-\sigma_v, 1-\tau\,\vert\, \tau\in I_v),$$
where $G_v$ and $I_v$ are the decomposition, respectively inertia groups of $v$ in $G$, $I_{G_v}$ is the ideal in $\Bbb Z[G]$ generated by $\{\sigma-1\mid \sigma\in G_v\}$
and $\sigma_v$ is any choice of Frobenius at $v$ in $G_v$.
\end{enumerate}
\end{lemma}
\begin{proof} Since we have a direct sum decomposition $Y_S(H)=\oplus_{v\in S} Y_v(H)$ of $\Bbb Z[G]$--modules, it suffices to prove the result above for a singleton $S=\{v\}$.
In this case, the map
$$Y_v(H)^\sharp\overset\sim\rightarrow Y_v(H)^\ast, \quad w\to w^\ast, \quad\text{ for all $w$ in $H$ with  }w|v, \quad\text{ where } w^\ast(w')=\begin{cases} 1, & \text{ if } w'=w,\\ 0, &\text{ if }w'\ne w,\end{cases}  $$
is a $\Bbb Z[G]$--module isomorphism, as the reader can easily check. This proves part (1).

Part (2) follows immediately from part (1), the obvious $\Bbb Z[G]$--module isomorphisms
$$Y_v(H)\simeq \Bbb Z[G/G_v]\simeq \Bbb Z[G]/I_{G_v}, $$
and equalities  $\sharp(I_{G_v})=I_{G_v}$ and $I_{G_v}=(\sigma_v-1, 1-\tau\mid \tau\in I_v)$, which can be  easily checked.
\end{proof}

The next result is essentially due to Cornacchia--Greither \cite{Cornacchia-Greither}.

\begin{proposition}\label{Cornacchia-Greither}
Let $R$ be a commutative, semi--local ring (i.e. a finite direct sum of local rings) and let $M$ be a finitely presented $R$--module. The following are equivalent.
\begin{enumerate}
\item There exists an $n\geq 1$ and a short exact sequence of $R$--modules
$$0\longrightarrow R^n\overset{\phi}\longrightarrow R^n\longrightarrow M\longrightarrow 0.$$
\item The Fitting ideal ${\rm Fitt}_R(M)$ is principal, generated by a non--zero divisor in $R$.
\end{enumerate}
\end{proposition}
\begin{proof}
This is essentially an extension to the general semi--local case of the equivalence $(2)\Leftrightarrow (3)$ in Proposition 4.2 of \cite{Cornacchia-Greither}, which only deals with a special type of a local ring. We will give a brief sketch of the proof.

If (1) is true, then (2) is clearly true because ${\rm Fitt}_R(M)$ is generated by ${\rm det}(\phi)$, which is a non--zero divisor in $R$ because $\phi$ is injective.

In order to prove that $(2)\Rightarrow (1)$, observe that since $R=\bigoplus_i R_i$ with $R_i$ local, then $M=\bigoplus_i M_i$, with $M_i:=M\otimes_R R_i$ and ${\rm Fitt}_{R_i}(M_i)=\pi_{R\to R_i}({\rm Fitt}_R(M))$ is a principal ideal, generated by a non--zero divisor in $R_i$, for all $i$. The proof of the implication $(3)\Rightarrow (2)$ in Proposition 4.2 in \cite{Cornacchia-Greither}, which only uses the fact that the ring in loc.cit. is local and not the special form of that ring, applies here and produces short exact sequences of $R_i$--modules
$$0\longrightarrow R_i^{n_i}\overset{\phi_i}\longrightarrow R^{n_i}\overset{\pi_i}\longrightarrow  M_i\longrightarrow  0,$$
for all $i$. Now, we take  $n:=\max\{n_i\vert i\}$ and obtain an exact sequence of $R$--modules
$$0\longrightarrow\ R^n\overset{\phi}\longrightarrow R^n\overset{\pi}\longrightarrow  M\longrightarrow 0,$$
where $\phi\vert_{R_i^n}=\phi_i\oplus {\rm id}_{n-n_i}$ and  $\pi\vert_{R_i^n}=\pi_i\oplus {\rm 0}_{n-n_i}$, for all $i$.
\end{proof}

\begin{remark}\label{semi-local-remark}
Examples of semilocal rings $R$ relevant in our context are the following.
$$\Bbb Z_p[G_n], \quad \Bbb Z_p[G_n]^-, \quad \Bbb Z_p[[\mathcal G]], \quad \Bbb Z_p[[\mathcal G]]^-$$
Indeed, in the direct sum decompositions \eqref{semilocal} the rings $R_n^\chi$ and $R_\infty^\chi$ are local, of respective maximal ideals
$$\frak m_n^\chi:=pR_n^\chi+I_{G_{n, p}}, \qquad \frak m_\infty^\chi:=pR_\infty^\chi+I_{G_{\infty, p}},$$ 
where  $I_{G_{n, p}}$ and $I_{G_{\infty, p}}$ are the augmentation ideals
in $R_n^\chi$ and $R_\infty^\chi$, respectively.
\end{remark}

\section{Quadratic Presentations of Selmer Modules}

In this section we discuss quadratic presentations of Selmer and related (Ritter--Weiss) modules in finite, Galois extensiosn $H/F$ as above.  The first subsection is a brief summary of results obtained in \cite{Ritter-Weiss} and \cite{Dasgupta-Kakde}. The second subsection is dedicated to studying transition maps between these quadratic presentations in towers of extensions of a given base field $F$, with an eye on eventually passing to a projective limit in an Iwasawa tower $H_\infty/F$.

\subsection{Quadratic presentations of Selmer modules in finite extensions}

In this section, we will recall some notations, constructions  and results from \cite{Ritter-Weiss} and  Appendix A of \cite{Dasgupta-Kakde}.
We build on notations developed in \S\ref{Selmer-section}. In particular, $H/F$ is a finite, abelian, CM extension of a totally real number field, of Galois group $G$, and $S$ and $T$ are two finite, disjoint sets of places in $F$, such that
$S_{ram}(H/F)\subseteq S\cup T$ and $S_\infty(F)\subseteq S$.

In what follows, if $w$ is a prime in $H$, then $H_w$ denotes the completion of $H$ in its $w$--adic topology. If $w$ is finite, then $\mathcal O_w$ denotes the ring of integeres in $H_w$, and $U_w$  denotes the subgroup of principal units inside $\mathcal O_w^\times$. Also, $J_H$ and $C_H$ denote the topological groups of id\`eles and id\`ele classes of $H$, respectively. We denote by $G_w$ the decomposition group of $w$ in $H/F$. We let $\Delta G$ and $\Delta G_w$ the augmentation ideals of $\Bbb Z[G]$ and $\Bbb Z[G_w]$, respectively.

\begin{definition}\label{larger-S'}  We call a finite set  $S'$ of places in $F$, a $(H/F, S, T)$--large set if.
\begin{itemize}
  \item $S\subseteq S'$ and $S'\cap T=\phi$
  \item  $Cl_{S'}^T(H)=1$
  \item $\bigcup_{w\in S'(H)} G_w=G$.
\end{itemize}
\end{definition}
\medskip

\noindent For the rest of this section, we fix an $(H/F, S, T)$--large
 set $S'$.\\

Let us  fix a place $v$ of $F$ and a place $w$ of $H$ which sits above $v$. Following Ritter and Weiss \cite{Ritter-Weiss}, we define a $G_w$--module $V_w$, as an extension of $\Delta G_w$ by $H_w^\times$
\begin{equation}\label{local-exact-sequence}
0\xrightarrow[]{} H_w^{\times}\xrightarrow[]{s}V_w\xrightarrow[]{}\Delta G_w\xrightarrow[]{}0, \end{equation}
whose extension class $\alpha_w$ corresponds via the canonical isomorphisms (see \cite{Ritter-Weiss})
$$Ext_{G_w}^1(\Delta G_w,H_w^{\times})=H^1(G_w,Hom(\Delta G_w,H_w^{\times}))\overset{\delta_w'}{\simeq}  H^2(G_w,H_w^{\times}), \qquad \alpha_w\to u_{H_w/F_v}, $$
to the local fundamental class $u_{H_w/F_v}\in  H^2(G_w,H_w^{\times})$.
Further, the composite injective map $\mathcal O_w^\times \subseteq H_w^\times\overset{s}{\longrightarrow} V_w$ gives rise to an exact sequence of $G_w$--modules
\begin{equation}\label{local-exact-sequence-units}0\xrightarrow[]{} O_w^{\times}\xrightarrow[]{}V_w \xrightarrow[]{}W_w\xrightarrow[]{}0,\end{equation}
where $W_w\simeq V_w/s(\mathcal O_w^\times)$ .
Moreover, since, $H_w^{\times}/O_w^{\times}\simeq \Bbb Z$ (with $G_w$ acting trivially on $\Bbb Z$),  we have the following exact sequence of $G_w$--modules.
\begin{equation}\label{local-exact-sequence-W}0\xrightarrow[]{}\mathbb{Z}\xrightarrow[]{i}W_w\xrightarrow[]{j}\Delta G_w\xrightarrow[]{}0.\end{equation}
A similar argument gives rise to a global analogue of the exact sequence \eqref{local-exact-sequence}
\begin{equation}\label{global-exact-sequence}0\xrightarrow[]{}C_H\xrightarrow[]{}D\xrightarrow[]{}\Delta G\xrightarrow[]{}0,\end{equation}
whose extension class $\alpha$ corresponds via the canonical isomorphisms (see \cite{Ritter-Weiss})
$$Ext_G^1(\Delta G,C_H)=H^1(G,Hom(\Delta G,C_H))\overset{\delta'}{\simeq} H^2(G,C_H), \qquad \alpha\to u_{H/F}$$
to the global fundamental class $u_{H/F}\in H^2(G,C_H).$

\begin{definition} Assume that $R$ is a set of finite places in $F$, and that for every $v\in R$ we fix a place $w$ in $H$ above $v$ and a $G_w$--module $M_w$.  Then, we define the following $G$--module
\[\widetilde{\prod_{v\in R}} M_w:=\prod_{v\in R} Ind_{G_w}^G M_w,
\]
where $Ind_{G_w}^G M_w:=M_w\otimes_{\Bbb Z[G_w]}\Bbb Z[G]$ is the usual induced module.
\end{definition}

Now, we can define the following $G$--modules
$$J:=\prod_{v\notin S\cup T}^{\sim} O_w^{\times}\times \prod_{v\in S}^{\sim} H_w^{\times}\times \prod_{v\in T}^{\sim} U_w, \qquad
J':=\prod_{v\notin S'\cup T}^{\sim} O_w^{\times}\times \prod_{v\in S'}^{\sim} H_w^{\times}\times \prod_{v\in T}^{\sim} U_w$$
$$V:=\prod_{v\notin S'\cup T}^{\sim} O_w^{\times}\times \prod_{v\in S'}^{\sim} V_w\times \prod_{v\in T}^{\sim} U_w, \qquad W:=\prod_{v\in S'\setminus S}^{\sim} W_w\times \prod_{v\in S}^{\sim} \Delta G_w, \qquad W':=\prod_{v\in S'}^{\sim} \Delta G_w.$$
If we combine \eqref{local-exact-sequence}--\eqref{local-exact-sequence-W}, we obtain the short exact sequences of $G$--modules
\begin{equation}\label{exact-sequences-J}0\xrightarrow[]{} J \xrightarrow[]{} V\xrightarrow[]{} W \xrightarrow[]{} 0,\qquad
0\xrightarrow[]{} J' \xrightarrow[]{} V \xrightarrow[]{} W' \xrightarrow[]{} 0.\end{equation}
Now, Theorem 1 in \cite{Ritter-Weiss} gives us the  the following commutative diagram of $G$--modules
\[\begin{tikzcd}
& 0\arrow{d} & 0\arrow{d} & 0\arrow{d} \\%
& O_{H,S,T}^{\times}\arrow{d} & V^{\theta}\arrow{d} & W^{\theta}\arrow{d} \\%
0\arrow{r} & J\arrow{r}\arrow{d}{\theta_J} & V \arrow{r}\arrow{d}{\theta_V} & W \arrow{r}\arrow{d}{\theta_W} & 0 \\%
0\arrow{r} & C_H\arrow{r}\arrow[two heads]{d} & D \arrow{r}\arrow{d} & \Delta G\arrow{r}\arrow{d} & 0 \\%
& Cl_S^T(H) & 0 & 0
\end{tikzcd}
\]
where $\theta_J$ and $\theta_W$ are the obvious maps and the existence and uniqueness of $\theta_V$ is a direct consequence of the compatibility between the local and global fundamental classes.
The snake lemma gives an exact sequence of $G$--modules
\begin{equation}\label{sequence-units-class-S}
0\xrightarrow[]{} O_{H,S,T}^{\times}\xrightarrow[]{} V^{\theta}\xrightarrow[]{} W^{\theta}\xrightarrow[]{} Cl_S^T(H)\xrightarrow[]{} 0.\end{equation}
An identical line of arguments, involving $J'$ and $W'$ instead of $J$ and $W$, produces the following exact sequence of $G$--modules.
\begin{equation}\label{sequence-class-units-S'}0\xrightarrow[]{} O_{H,S',T}^{\times}\xrightarrow[]{} V^{\theta}\xrightarrow[]{} W'^{\theta}\xrightarrow[]{} 0.\end{equation}
Moreover, the definition of the maps $\theta_W$ and $\theta_{W'}$ combined with \eqref{local-exact-sequence-W} gives a short exact sequence of $G$--modules
\begin{equation}\label{sequence-W-W'}0\to Y_{S'\setminus S}\to W^\theta\to W'^{\theta}\to 0.\end{equation}\\

Now, we define two additional $G$--modules
\[B:= \prod_{v\in S'}^{\sim}\mathbb{Z}[G_w]\simeq \prod_{v\in S'}\Bbb Z[G], \qquad
Z:=\prod_{v\in S}^{\sim} \mathbb{Z}.\]
Then, Appendix A1 of \cite{Dasgupta-Kakde} gives us a commutative diagram of $G$- modules
\begin{equation}\label{diagram-B}
\begin{tikzcd}
& 0\arrow{d} & 0\arrow{d} & 0\arrow{d} \\%
0\arrow{r}& W^{\theta}\arrow{d}\arrow{r} & B^{\theta}\arrow{d}\arrow{r} & Z^{\theta}\arrow{r}\arrow{d}&0 \\%
0\arrow{r} & W\arrow{r}{\gamma}\arrow{d}{\theta_W} & B \arrow{r}\arrow{d}{\theta_B} & Z \arrow{r}\arrow{d}{\theta_Z} & 0 \\%
0\arrow{r} & \Delta G\arrow{r}\arrow{d} & \mathbb{Z}[G]\arrow{r}\arrow{d} & \mathbb{Z}\arrow{r}\arrow{d} & 0 \\%
& 0 & 0 & 0
\end{tikzcd}
\end{equation}
where $W^\theta$, $B^\theta$ and $Z^\theta$ are the kernels of $\theta_W$, $\theta_B$ and $\theta_Z$, respectively. (See loc.cit for the precise definitions of the maps involved in the diagram above.)
\begin{definition}\label{Ritter-Weiss-module-definition}
Define the Ritter-Weiss module associated to the data $(H/F, S, T)$ by
\[\nabla_S^T(H):=coker(V^{\theta}\xrightarrow[]{} W^{\theta} \xrightarrow[]{} B^{\theta}) \]
\end{definition}
If we combine \eqref{sequence-units-class-S} and \eqref{diagram-B}, we obtain the following exact seqeunces of $G$--modules.
\begin{equation}\label{nabla-sequences}
0\xrightarrow[]{} O_{H,S,T}^{\times}\xrightarrow[]{} V^{\theta} \xrightarrow[]{} B^{\theta}\xrightarrow[]{}\nabla_S^T(H)\xrightarrow[]{}0, \qquad 0\xrightarrow[]{} Cl_S^T(H)\xrightarrow[]{}\nabla_S^T(H)\xrightarrow[]{} Z^\theta\xrightarrow[]{} 0.
\end{equation}

\begin{lemma}[Appendix A, \cite{Dasgupta-Kakde}]\label{nabla-lemma} The following hold.
\begin{enumerate}\item $\nabla_S^T(H)$ does not depend on the auxiliary set of primes $S'$.
\item If $p>2$ is a prime and $S_p(F)\cap S_{ram}(H/F)\subseteq S$, then $(B^\theta)_p$ and $(V^\theta)_p$ are free $\Bbb Z_p[G]$--modules of constant local rank equal to $(|S'|-1)$
\end{enumerate}
\end{lemma}
\begin{proof} See loc.cit. \end{proof}
Now, one can define the transposed Ritter-Weiss module (in the sense of \cite{Jannsen}) as follows.
\begin{definition}\label{nabla-transposed-definition}
The transposed Ritter--Weiss module for the data $(H/F, S, T)$ is given by
$$\nabla_S^T(H)^{tr}:=coker((B^\theta)^\ast\to (V^\theta)^\ast)=coker((W^\theta)^\ast\to (V^\theta)^\ast),$$
where the $\Bbb Z$--duals in question are endowed with the contravariant $G$--action, as usual, and the equality above is a consequence of the obvious surjectivity
of $(B^\theta)^\ast\to(W^\theta)^\ast$.
\end{definition}

Now, by taking $\Bbb Z$--duals in \eqref{sequence-units-class-S}--\eqref{sequence-W-W'}, we obtain a commutative diagram of $G$--modules.
\begin{equation}\label{ritter-weiss-transpose}
    \begin{tikzcd}
& 0\arrow{d} & 0\arrow{d} \\%
0\arrow{r} & (W'^{\theta})^\ast\arrow{r}\arrow{d} & (V^{\theta})^\ast\arrow{r}\arrow[equal]{d} & (O_{H,S',T}^{\times})^\ast\arrow{r}\arrow{d} & 0 \\%
0\arrow{r} & (W^{\theta})^\ast\arrow{r}\arrow[two heads]{d} & (V^{\theta})^\ast\arrow{r} \arrow{d} & \nabla_S^T(H)^{tr}\arrow{r} & 0 \\%
& (Y_{S'\setminus S})^\ast & 0
\end{tikzcd}
\end{equation}
The snake lemma applied to the diagram above produces an exact sequence of $G$--modules
$$0\to (Y_{S'\setminus S})^\ast\to (O_{H, S', T}^\times)^\ast\to {\nabla_S^T(H)}^{tr}\to 0.$$
When combined with \eqref{large-S-iso}, this exact sequence gives an isomorphism of $G$--modules
\begin{equation}\label{selmer-nablatr-iso} Sel_S^T(H)\simeq {\nabla_S^T(H)}^{tr}.\end{equation}
(See Appendix A in \cite{Dasgupta-Kakde} for more details.)  Consequently, we obtain the following.

\begin{lemma}\label{quadratic-pres-lemma} For any $(H/F, S, T, S')$ as above, and any prime $p>2$, we have the following.
\begin{enumerate}
\item The associated Ritter--Weiss and Selmer modules admit $\Bbb Z[G]$--module presentations
$$V^\theta\to B^\theta\to\nabla_S^T(H)\to 0, \qquad (B^\theta)^\ast\to (V^\theta)^\ast\to Sel_S^T(H)\to 0.$$
\item If $S_p(F)\cap S_{ram}(H/F)\subseteq S$, then the resulting $\Bbb Z_p[G]$--module presentations
$$(V^\theta)_p\to (B^\theta)_p\to\nabla_S^T(H)_p\to 0, \qquad (B^\theta)^\ast_p\to (V^\theta)^\ast_p\to Sel_S^T(H)_p\to 0$$
are $\Bbb Z_p[G]$--free, quadratic, of rank $(|S'|-1)$.
\end{enumerate}
\end{lemma}

\begin{proof} See Appendix A in \cite{Dasgupta-Kakde}.\end{proof}

The quadratic presentations in Lemma \ref{quadratic-pres-lemma}(2) will play a very important role in the proof of the main results of this paper. We will study them further in the next section.
\medskip

\subsection{Transition maps between quadratic presentations}

In this section, we consider two finite, abelian, CM extensions $K_1$ and $K_2$ of a totally real number field $F$, such that $K_1\subseteq K_2$. We take sets of places $S$, $T$ in $F$ satisfying $S_{ram}(K_2/F)\subseteq S\cup T$ and $S_\infty(F)\subseteq S$, and an auxiliary set of places $S'$ in $F$, satisfying the required ``largeness'' conditions (see Definition \ref{larger-S'}) for both $(K_2/F, S, T)$ and $(K_1/F, S, T)$. 

Our main goal is to construct maps between the presentations provided by Lemma \ref{quadratic-pres-lemma} for the Selmer groups $Sel_S^T(K_2)$ and $Sel_S^T(K_1)$ and their $p$--primary parts, compatible with the restriction map 
$$res: Sel_S^T(K_2)\to Sel_S^T(K_1).$$
This section is rather technical, so the reader is advised to skip ahead to the final results (Theorems \ref{Z-presentations-compatible} and \ref{QP-compatibility}, and especially Corollary \ref{covariance}) and go over the technical details afterwards. \\

In what follows, we let $G_1:=Gal(K_1/F)$, $G_2:=Gal(K_2/F)$, $G_{2,1}:=G(K_2/K_1)$, and let  $\pi:G_2\twoheadrightarrow G_1$ be the Galois restriction morphism. All $\Bbb Z[G_1]$--modules are viewed as $\Bbb Z[G_2]$--modules via $\pi$.
Throughout, if $M$ is any of the modules $J, J', V, W, W'$ etc. introduced in the last section, we denote by $M_1$ and $M_2$ (and at times $M^1$ and $M^2$) the corresponding modules at levels $K_1$ and $K_2$, respectively.
\\

We start with a local construction. To that end, let us fix a place $u$ of $F$, a place $v$ of $K_1$ above $u$ and a place $w$ of $K_2$ above $v$. Let $F_u$, $K_{1,v}$ and $K_{2,w}$ be completions with respect to those places. Then the decomposition groups of $u$ in $K_1/F$ and $K_2/F$ are $G_v:=Gal(K_{1,v}/F_u)$ and $G_w:=Gal(K_{2,w}/F_u)$, respectively. Further, let $G_{w/v}:=G_w\cap G_{2,1}$. All $\Bbb Z[G_v]$--modules are viewed as $\Bbb Z[G_w]$--modules via the Galois restriction $\pi_v: G_w\twoheadrightarrow G_v$.
Define the $\mathbb{Z}[G_w]$--module morphism
$$\psi_v:\mathbb{Z}[G_v] \xrightarrow[]{} \mathbb{Z}[G_w], \qquad \psi_v(g)=\sum_{\pi_v(\tilde g)=g}\widetilde{g}.$$
At the level of augmentation ideals, we have
$\psi_v(\Delta G_v)\subseteq\Delta G_w$. The following gives maps between the short exact sequences \eqref{local-exact-sequence} at levels $K_{1,v}$ and $K_{2,w}$.

\begin{proposition}\label{compatibility}
There exists a unique $\mathbb{Z}[G_w]$--module morphism $f_v$ making the following diagram
\[\begin{tikzcd}
0\arrow{r} & K_{2,w}^{\times}\arrow{r}\arrow[leftarrow]{d}{i_v} & V_w^2 \arrow{r}\arrow[leftarrow]{d}{f_v} & \Delta G_w \arrow{r} \arrow[leftarrow]{d}{\psi_v} & 0 \\%
0\arrow{r} & K_{1,v}^{\times}\arrow{r} & V_v^1 \arrow{r} & \Delta G_v \arrow{r} & 0
\end{tikzcd}
\]
commutative, where $i_v$ is the usual inclusion and $\psi_v$ is the map defined above.
\end{proposition}
\begin{proof}
We have the following commutative diagram of $\Bbb Z[G_w]$--modules with exact rows.

\[\begin{tikzcd}
0\arrow{r} & K_{2,w}^{\times}\arrow{r}\arrow[leftarrow]{d}{} & V_w^2 \arrow{r}\arrow[leftarrow]{d}{} & \Delta G_w \arrow{r} \arrow[leftarrow]{d}{} & 0 \\%
0\arrow{r} &(K_{2,w}^{\times})^{G_{w/v}}\arrow{r}\arrow[leftarrow,equal]{d}{} & (V_w^2)^{G_{w/v}} \arrow{r}\arrow[leftarrow]{d}{} & (\Delta G_w)^{G_{w/v}} \arrow{r} \arrow[leftarrow]{d}{} & 0 \\%
0\arrow{r} & K_{1,v}^{\times}\arrow{r} & V_v^1 \arrow{r} & \Delta G_v \arrow{r} & 0
\end{tikzcd}
\]
Here, the maps between the second and first row are induced by inclusion. The maps between the third and second row are given in diagram (163) in \cite{Dasgupta-Kakde}. In particular (see loc.cit.), the map
$$\Delta G_v\to (\Delta G_w)^{G_{w/v}}\hookrightarrow\Delta G_w$$ is our map $\psi_v$. Consequently, the desired map $f_v$ is the composition of the middle vertical maps.
\end{proof}

\begin{remark}\label{V-W-local} For primes $v\in S'$, since $i_v(\mathcal O_v^\times)\subseteq \mathcal O_w^\times$, the maps $f_v: V_v^1\to V_w^2$ clearly induce maps $f_v': W_v^1\to W_w^2$ which make the following diagram commute
\[\begin{tikzcd}
V_w^2\arrow{r} & W_w^2 \arrow{r}& \Delta G_w\\%
V_v^1\arrow{r}\arrow{u}{f_v} & W_v^1 \arrow{r}\arrow{u}{f_v'}& \Delta G_v\arrow{u}{\psi_v}
\end{tikzcd}
\]
where the rows are given by the canonical morphisms between the given modules.
\end{remark} 
\begin{lemma}\label{induction-lemma}  Let $M_v$ be a $\Bbb Z[G_v]$--module, $M_w$ a $\Bbb Z[G_w]$--module, and $\xi:M_v\to M_w$ a $\Bbb Z[G_w]$--linear map. Then, we have a well defined  $\Bbb Z[G_2]$--linear morphism
$$\widetilde\xi: M_v\otimes_{\Bbb Z[G_v]}\Bbb Z[G_1]\to M_w\otimes_{\Bbb Z[G_w]}\Bbb Z[G_2],\qquad \widetilde\xi(m\otimes\sigma):=\xi(x)\otimes(\tilde\sigma\cdot\sum_{\tau}\tau),$$
where $\tilde\sigma\in G_2$ is a lift of $\sigma\in G_1$ and $\tau$ runs through a complete set of $G_{w/v}$--coset representatives of $G_{2,1}$.
\end{lemma}
\begin{proof} This is a consequence of the inclusion $\xi(M_v)\subseteq M_w^{G_{w,v}}$. We leave the details to the reader. \end{proof}

\begin{remark}\label{Psi-definition} Observe that, via the procedure described in the Lemma,  the local map $\psi_v:\Bbb Z[G_v]\to\Bbb Z[G_w]$ defined above induces the $\Bbb Z[G_2]$--linear global map
$$\widetilde{\psi_v}:Ind_{G_v}^{G_1}\Bbb Z[G_v]\cong\Bbb Z[G_1] \longrightarrow Ind_{G_w}^{G_2}\Bbb Z[G_w]\cong\Bbb Z[G_2],\qquad \widetilde{\psi_v}(g)=\sum_{\tilde g\twoheadrightarrow g}\tilde g,$$
where $g\in G_1$ and $G_2\twoheadrightarrow  G_1$ is the Galois restriction. Note that $\widetilde{\psi_v}$ does not depend on $v$ and $\widetilde{\psi_v}(\Delta G_1)\subseteq \Delta G_2$.
\end{remark}

Next, we construct (global) maps between the short exact sequences \eqref{sequence-class-units-S'} at levels $K_1$ and $K_2$.  
\begin{proposition}\label{global-compatibility}
We have a commutative diagram of $\Bbb Z[G_2]$--modules with exact rows
\[\begin{tikzcd}
0\arrow{r} & O_{K_2,S',T}^{\times}\arrow[leftarrow]{d}{\iota}\arrow{r} & V_2^{\theta}\arrow[leftarrow]{d}{\widetilde f}\arrow{r} & W_2^{'\theta}\arrow[leftarrow]{d}{\widetilde \psi}\arrow{r} & 0\\%
0\arrow{r} & O_{K_1,S',T}^{\times}\arrow{r} & V_1^{\theta}\arrow{r} & W_1^{'\theta}\arrow{r} & 0
\end{tikzcd}
\]
where the left vertical map is inclusion, and $\widetilde f$ and $\widetilde\psi$ are componentwise induced as in Lemma \ref{induction-lemma} by the maps $\{f_v\}_v$ and $\{\psi_v\}_v$ defined above.
\end{proposition}
\begin{proof}
From the definitions of the modules $W_i'$ and $V_i$, we have a $\Bbb Z[G_2]$--linear morphism
$$\widetilde{\psi}=\{\widetilde{\psi_v}\}_{v\in S'}:W'_1 \xrightarrow[]{} W'_2, \qquad \widetilde{f}=\{\widetilde{f_v}\}_v:V_1\xrightarrow[]{} V_2 .$$
Now, Proposition \ref{compatibility} leads to the following commutative diagram of $\Bbb Z[G_2]$--modules
\begin{equation}\label{JVW-diagram}\begin{tikzcd}
0\arrow{r} & J'_2\arrow[leftarrow]{d}{\iota}\arrow{r} & V_2\arrow[leftarrow]{d}{\widetilde f}\arrow{r} & W'_2\arrow[leftarrow]{d}{\widetilde\psi}\arrow{r} & 0\\%
0\arrow{r} & J'_1\arrow{r} & V_1\arrow{r} & W'_1\arrow{r} & 0
\end{tikzcd}
\end{equation}
where the left vertical arrow is the usual inclusion map. Now, observe that we have the following commutative diagrams of $\Bbb Z[G_2]$--modules
\begin{equation}\label{W-J-V-diagrams}\begin{tikzcd}
W'_1\arrow{d}{\theta_{W'_1}}\arrow{r}{\Tilde{\psi}} & W'_2\arrow{d}{\theta_{W'_2}} & \\%
\Delta G_1 \arrow{r}{\Tilde{\psi}} & \Delta G_2
\end{tikzcd}\qquad
\begin{tikzcd}
J_1'\arrow{d}{\theta_{J'_1}}\arrow{r}{\iota} & J'_2\arrow{d}{\theta_{J'_2}} & \\%
C_{K_1} \arrow{r} & C_{K_2}
\end{tikzcd}
\qquad
\begin{tikzcd}
V_1\arrow{d}{\theta_{V_1}}\arrow{r}{\widetilde f} & V_2\arrow{d}{\theta_{V_2}} & \\%
D_1 \arrow{r} & D_2
\end{tikzcd}
\end{equation}
where the first two diagrams are obviously commutative from the definitions and the third is essentially diagram (164) of \cite{Dasgupta-Kakde}. This implies right away that
$\widetilde{\psi}(W_1^{'\theta})\subseteq W_2^{'\theta}$ and $\widetilde{f}(V_1^\theta)\subseteq V_2^\theta$ and also that $\iota$ induces the inclusion map
$J_1^\theta=\mathcal O_{K_1, S', T}^\times \to \mathcal O_{K_2, S', T}^\times=J_2^\theta.$ Consequently, \eqref{JVW-diagram} and \eqref{W-J-V-diagrams} give the desired commutative diagram in the statement of the proposition.
\end{proof}
Next, we relate the transition maps above to the restriction maps between Selmer modules.  With the notation in \cite{Dasgupta-Kakde}, if $G$ is a finite, abelian group, $M$ is a $\Bbb Z[G]$--module, and $p$ is a prime number, we let
$$F(M):={\rm Hom}_{\Bbb Z} (M, \Bbb Z)= M^\ast,\qquad F_p(M):={\rm Hom}_{\Bbb Z_p} (M\otimes\Bbb Z_p, \Bbb Z_p)=(M_p)^\ast,$$
endowed, as usual, with the covariant  $G$--action.
\begin{proposition}\label{restriction-compatibility}
We have the following commutative diagram of $\Bbb Z[G_2]$--modules with exact rows
\[\begin{tikzcd}
0\arrow{r} & F(W_2^{\theta})\arrow{d}{F(\widetilde{f'})}\arrow{r} & F(V_2^{\theta})\arrow{d}{F(\tilde f)}\arrow{r} & Sel_S^T(K_2)\arrow{d}{res}\arrow{r} & 0\\%
0\arrow{r} & F(W_1^{\theta})\arrow{r} & F(V_1^{\theta})\arrow{r} & Sel_S^T(K_1)\arrow{r} & 0,
\end{tikzcd}
\]
where the rows are as in diagram \eqref{ritter-weiss-transpose} and $\widetilde{ f'}:=(\widetilde{{f'_v}})_{v\in{S'\setminus S}}\oplus (\widetilde{\psi_v})_{v\in S}$, with $f_v'$ defined in Remark \ref{V-W-local}.
\end{proposition}
\begin{proof}
Apply the functor $F$ to the diagram in Proposition \ref{global-compatibility}. Observe that, since $W_1^{\theta}$ and $W_2^{\theta}$ are $\mathbb{Z}$--free, the rows will stay exact and the right vertical map is the restriction map.

\begin{equation}\label{dual-W'V-compatibility}\begin{tikzcd}
0\arrow{r} & F(W_2^{'\theta})\arrow{d}{F(\tilde\psi)}\arrow{r} & F(V_2^{\theta})\arrow{d}{F(\tilde f)}\arrow{r} & F(O_{K_2,S',T}^{\times})\arrow{d}{F(\iota)}\arrow{r} & 0\\%
0\arrow{r} & F(W_1^{'\theta})\arrow{r} & F(V_1^{\theta})\arrow{r} & F(O_{K_1,S',T}^{\times})\arrow{r} & 0
\end{tikzcd}
\end{equation}
Consequently, we have the following diagram with commutative squares,
$$
    \begin{tikzcd}
F(V_2^{\theta})\arrow{d}{F(\tilde f)}\arrow{r} & F(O_{K_2,S',T}^{\times})\arrow{d}{F(\iota)}\arrow{r} & Sel_S^T(K_2)\arrow{d}{res} & \\%
F(V_1^{\theta})\arrow{r}& F(O_{K_1,S',T}^{\times})\arrow{r} & Sel_S^T(K_1)
\end{tikzcd}
$$
where the left square is commutative by \eqref{dual-W'V-compatibility} and the right square is commutative from the definition of the Selmer group and the fact that $F(\iota)$ is the restriction map. Therefore, we have just proved the commutativity of the right--most square of the diagram in the statement of the Proposition.

Next, note that Remark \ref{V-W-local} yields a diagram of $\Bbb Z[G_2]$--modules with commutative cells 
\begin{equation}\label{W-theta-compatible}
    \begin{tikzcd}
V_2\arrow{r} & W_2\arrow{r}{\theta_{W_2}} & \Delta G_2\\%
V_1\arrow{r}\arrow{u}{\widetilde{f}} & W_1\arrow{u}{\widetilde f'}\arrow{r}{\theta_{W_1}} &\Delta G_1\arrow{u}{\widetilde \psi}.
\end{tikzcd}
\end{equation}
On the other hand, if $*\in\{1, 2\}$, we also have the following diagram, where the horizonthal maps are canonical.
\[\begin{tikzcd}
V_\ast\arrow{d}{\theta_{V_\ast}}\arrow{r} &
W_\ast\arrow{d}{\theta_{W_\ast}} \\%
D_\ast\arrow{r} &\Delta G_\ast
\end{tikzcd}
\]
Hence, by taking $\theta$--kernels, the last two commutative diagrams give us the following commutative diagram.
\[\begin{tikzcd}
V_2^{\theta}\arrow{r} & W_2^{\theta} & \\%
V_1^{\theta}\arrow{r}\arrow{u}{\widetilde{f}} & W_1^{\theta}\arrow{u}{\widetilde{f'}}
\end{tikzcd}
\]
After applying the functor $F$ to the above diagram, we obtain the commutativity of the left square of the diagram in the Proposition. The exactness of the rows is clear from \eqref{ritter-weiss-transpose}.
\end{proof}
\medskip

Now, let us look at the maps $\gamma:W\xrightarrow[]{} B$ as in diagram \ref{diagram-B} at the levels $K_1$ and $K_2$, and then construct maps between the two levels. In order to do that, let us understand the modules $W$ more explicitly. Here we are using the description in \cite{Gruenberg-Weiss}. We focus on one level, say $K_1$. Let $\overline{G_v}=\langle F \rangle$ be the Galois group of the residual field extension corresponding to $K_{1,v}/F_u$,  where $F$ is the Frobenius automorphism. Then, we have 
$$W_v=\{(x,y)\in \Delta G_v\bigoplus \mathbb{Z}[\overline{G_v}]\, \vert\, \Bar{x}=(F-1)y\}$$ where $\Bar{x}$ is the image of $x$ via the projection $\Bbb Z[G_v]\to \mathbb{Z}[\overline{G_v}]$. Clearly, $W_v$ is a free $\mathbb{Z}$--module of basis 
$$\{w_g:=(g-1,\sum_{i=0}^{a(g)-1}F^i)\, \vert\, g\in G_v\},\qquad\text{ where } \Bar{g}=F^{a(g)} \text { and }0<a(g)\leq l_v:=|\overline{G_v}|.$$
With notation as in \eqref{local-exact-sequence-W}, we have $$i(1)=w_1, \qquad j(w_g)=g-1, \quad \text{for all }g\in G_v.$$ 
The $G_v$--action on these basis elements is given by $$g\cdot w_h=w_{gh}-w_g+a_{g,h}w_1, \qquad \text{where } a(g)+a(h)=a(gh)+l_va_{g,h},$$
for all $g, h\in G_v$. Observe that $a_{1,h}=1$ for each $h\in G_v.$\\\\
Next, we prove technical result about the quantities $a(g)$ and $a_{g, h}$ introduced above.
\begin{proposition}\label{formulas}
Let $\Tilde{h}\in G_w$ is a lift of $h\in G_v$ and let $g\in G_v$. Let $l$ be the residue class degree and $e$ be the ramification index of the extension $K_{2,w}/K_{1,v}$ Then, the following hold.\\\\
\begin{enumerate}
\item $a(\Tilde{h})=a(h)+k_{\Tilde{h}}l_v$,  for some $k_{\Tilde{h}}\in\{0,1,2,...\text{ }  l-1\}$.\\
\item $\sum_{\Tilde{g}\twoheadrightarrow g} a_{\Tilde{g},\Tilde{h}}=e(a_{g,h}+k_{\Tilde{h}})$
\end{enumerate}
\end{proposition}
\begin{proof}
(1) If $F_w$ and $F_v$ are the Frobenius morphisms at levels $K_2$ and $K_1$, respectively, we have, $\overline{\Tilde{h}}=F_w^{a(\Tilde{h})}$ and $\overline{h}=F_v^{a(h)}$. However, under Galois restriction, $F_w$ maps to $F_v$. Hence, we have $F_v^{a(\Tilde{h})}=F_v^{a(h)}$. Hence,  $a(\Tilde{h})=a(h)+k_{\Tilde{h}}l_v.$ for some $k_{\Tilde{h}}\in\mathbb{Z}$. But, since $0<a(\Tilde{h})\leq l_w$ and $l=l_w/l_v$, we have $0\leq k_{\Tilde{h}}<l$ as desired.\\\\
(2) If $\Tilde{g}\in G_w$ is a lift of $g$, by the definition of $a_{\Tilde{g},\Tilde{h}}$ we have
 $$a(\Tilde{g})+a(\Tilde{h})=a(\Tilde{g}\Tilde{h})+l_wa_{\Tilde{g},\Tilde{h}}.$$ By applying part (1), we get
$$a(g)+k_{\Tilde{g}}l_v+a(h)+k_{\Tilde{h}}l_v=a(gh)+k_{\Tilde{g}\Tilde{h}}l_v+l_wa_{\Tilde{g},\Tilde{h}}.$$
In the equality above, take the sum when $\Tilde{g}$ varies through all the lifts of $g$. Since, $k_{\Tilde{g}}$ is distributed evenly among the elements in $\{0,1,...\text{ }l-1\}$ when $\Tilde{g}$ varies through all the lifts of $g$, if we take the summation, the terms  $k_{\Tilde{g}}$ and $k_{\Tilde{g}\Tilde{h}}$  cancel out. Therefore, we have 
$$m\cdot a(g)+m\cdot (a(h)+k_{\Tilde{h}}l_v)=m\cdot a(gh)+l_w\sum_{\Tilde{g}\xrightarrow[]{}g}a_{\Tilde{g},\Tilde{h}},\qquad 
ml_v(a_{g,h}+k_{\Tilde{h}})=l_w\sum_{\Tilde{g}\xrightarrow[]{}g}a_{\Tilde{g},\Tilde{h}},$$
where $m=[K_{2,w}:K_{1,v}]$. Since $m=el$, we have the desired result.
\end{proof}
Next, we give an explicit description of the $\Bbb Z[G_w]$--linear map $f_v'$ defined in Remark \ref{V-W-local}.
\begin{proposition}
The map $f_v':W_v^1\xrightarrow[]{} W_w^2$ is given by $$f_v'(w_g)=\sum_{\Tilde{g}\xrightarrow[]{}g} w_{\Tilde{g}}-\sum_{\Tilde{g}\xrightarrow[]{}1} w_{\Tilde{g}}+ew_{\Tilde{1}}$$ where $e$ is the ramification index of the extension $K_{2,w}/K_{1,v}$.
\end{proposition}
\begin{proof}
From Propostion \ref{compatibility}, we obtain the following commutative diagram whose rows are the exact sequences \eqref{local-exact-sequence-W} at levels $K_1$ and $K_2$, respectively.
\begin{equation}\label{f'-W-compatibility}
    \begin{tikzcd}
0\arrow{r} & \mathbb{Z}\arrow{r}{i_2}\arrow[leftarrow]{d}{\times e} & W_w^2\arrow{r}{j_2}\arrow[leftarrow]{d}{f_v'} & \Delta G_w \arrow{r}\arrow[leftarrow]{d}{\psi_v} & 0 \\%
0\arrow{r} & \mathbb{Z}\arrow{r}{i_1} & W_v^1\arrow{r}{j_1} & \Delta G_v \arrow{r} & 0
\end{tikzcd}
\end{equation}
Observe that for any $g\in G_v$, by the commutativity of the right square, we have $$j_2(f_v'(w_g))=\psi_v(j_1(w_g))=\psi_v(g-1)=\sum_{\Tilde{g}\xrightarrow[]{}g}\Tilde{g}-\sum_{\Tilde{g}\xrightarrow[]{}1}\Tilde{g}, \qquad j_2(\sum_{\Tilde{g}\xrightarrow[]{}g}w_{\Tilde{g}}-\sum_{\Tilde{g}\xrightarrow[]{}1}w_{\Tilde{g}})=\sum_{\Tilde{g}\xrightarrow[]{}g}\Tilde{g}-\sum_{\Tilde{g}\xrightarrow[]{}1}\Tilde{g}.$$
Therefore, by the exactness of the upper row in the above diagram, we have $$f'_v(w_g)=\sum_{\Tilde{g}\xrightarrow[]{}g}w_{\Tilde{g}}-\sum_{\Tilde{g}\xrightarrow[]{}1}w_{\Tilde{g}}+e_g w_{\Tilde{1}}$$ for some $e_g\in\mathbb{Z}$. Here, $\Tilde{1}$ is the identity element in group $G_w$. 

We need to prove that $e_g=e$, for all $g\in G_v$.
By the $\Bbb Z[G_w]$--linearity of $f'_v$,  if $g,h\in G_v$ and if $\Tilde{h}$ is any lift of $h$ in $G_w$, we have $\Tilde{h}\cdot(f'_v(w_g))=f'_v(h\cdot w_g)$.
Further, observe that
\begin{align}
\Tilde{h}\cdot(f_v'(w_g))&=\sum_{\Tilde{g}\xrightarrow[]{}g}(w_{\Tilde{g}\Tilde{h}}-w_{\Tilde{h}}+a_{\Tilde{g},\Tilde{h}}w_{\Tilde{1}})-\sum_{\Tilde{g}\xrightarrow[]{}1}(w_{\Tilde{g}\Tilde{h}}-w_{\Tilde{h}}+a_{\Tilde{g},\Tilde{h}}w_{\Tilde{1}})+e_g w_{\Tilde{1}} \nonumber \\
&=\sum_{\Tilde{gh}\xrightarrow[]{}gh}w_{\Tilde{gh}}-\sum_{\Tilde{h}\xrightarrow[]{}h}w_{\Tilde{h}}+(\sum_{\Tilde{g}\xrightarrow[]{}g}a_{\Tilde{g},\Tilde{h}}-\sum_{\Tilde{g}\xrightarrow[]{}1}a_{\Tilde{g},\Tilde{h}}+e_g)w_{\Tilde{1}}. \nonumber
\end{align}
Via Proposition \ref{formulas}(2), we can simplify this further:
\begin{align}
 \Tilde{h}\cdot(f_v'(w_g))&=\sum_{\Tilde{gh}\xrightarrow[]{}gh}w_{\Tilde{gh}}-\sum_{\Tilde{h}\xrightarrow[]{}h}w_{\Tilde{h}}+e w_{\Tilde{1}}[(a_{g,h}+k_{\Tilde{h}})-(a_{1,h}+k_{\Tilde{h}})]+e_g w_{\Tilde{1}} \nonumber \\
 &=\sum_{\Tilde{gh}\xrightarrow[]{}gh}w_{\Tilde{gh}}-\sum_{\Tilde{h}\xrightarrow[]{}h}w_{\Tilde{h}}+(ea_{g,h}-e+e_g)w_{\Tilde{1}}. \nonumber
\end{align}
On the other hand, we have
\begin{align}
f_v'(h.w_g)&=f_v'(w_{gh}-w_h+a_{g,h}w_1) \nonumber \\
&=(\sum_{\Tilde{gh}\xrightarrow[]{}gh}w_{\Tilde{gh}}-\sum_{\Tilde{g}\xrightarrow[]{}1}w_{\Tilde{g}}+e_{gh}w_{\Tilde{1}})-(\sum_{\Tilde{h}\xrightarrow[]{}h}w_{\Tilde{h}}-\sum_{\Tilde{g}\xrightarrow[]{}1}w_{\Tilde{g}}+e_{h}w_{\Tilde{1}})+ea_{g,h}w_{\Tilde{1}} \nonumber \\
&=\sum_{\Tilde{gh}\xrightarrow[]{}gh}w_{\Tilde{gh}}-\sum_{\Tilde{h}\xrightarrow[]{}h}w_{\Tilde{h}}+(e_{gh}-e_h+ea_{g,h})w_{\Tilde{1}} \nonumber
\end{align}
Therefore, for all $g,h\in G_v$, we have $e_g-e=e_{gh}-e_h$. Now, by taking the summation when $h$ varies through all elements in $G_v$, we get $e_g=e$ as desired. This completes the proof.
\end{proof}
Now let us recall the description in \cite{Dasgupta-Kakde} and \cite{Gruenberg-Weiss} of the map $\gamma:W\xrightarrow[]{} B$ in diagram \eqref{diagram-B}. Its componentwise definition at the level $K_1$ is as follows.
\begin{itemize}
    \item For $v\in S,$ $\gamma_v$ is induced by the inclusion $\Delta G_v\subseteq\mathbb{Z}[G_w]$
    \item For $v\in S'\setminus S$, $\gamma_v$ is induced by the map $s$ given by
\begin{equation*}s(w_g)=\sum_{h\in G_v}(r(g)+1-a_{g^{-1},h})h, \qquad\text{ where }
 r(g) =  \left\{
        \begin{array}{ll}
            1 & if \quad g\in I_v \\
            0 & if \quad g\notin I_v
        \end{array}
    \right., \qquad\text{ for all }g\in G_v,
\end{equation*} 
and $I_v$ i s the ramification group. 
 
\end{itemize}
From the right square of the commutative diagram \eqref{f'-W-compatibility}, the map $f'$ is compatible with the maps $j$ and $\psi$. Next, we prove a similar result for the maps $s$ at levels $K_{1,v}$ and $K_{2,w}$.

\begin{proposition}\label{W-compatibility}
We have a commutative diagram of $\Bbb Z[G_w]$--modules
\[\begin{tikzcd}
W_w^2\arrow{r}{s_2}\arrow[leftarrow]{d}{f_v'} & \mathbb{Z}[G_w]\arrow[leftarrow]{d}{e\psi_v} \\%
W_v^1\arrow{r}{s_1} & \mathbb{Z}[G_v],
\end{tikzcd}
\]
where $e:=e(w/v)$ is the ramification index.
\end{proposition}
\begin{proof}
Observe that for any $g\in G_v$ we have
\begin{align}
    e\psi_v(s_1(w_g)) &= e\psi_v(\sum_{h\in G_v}(r(g)+1-a_{g^{-1},h})h) \nonumber \\
    &= e\sum_{h\in G_v}(r(g)+1-a_{g^{-1},h})\cdot \sum_{\Tilde{h}\xrightarrow[]{}h} \Tilde{h}  \nonumber  \\
    &=e\sum_{\Tilde{h}\in G_w}(r(g)+1-a_{g^{-1},h})\Tilde{h}. \nonumber
\end{align}
On the other hand, we have
\begin{align}
    s_2(f_v'(w_g)) &= s_2(\sum_{\Tilde{g}\xrightarrow[]{}g} w_{\Tilde{g}}-\sum_{\Tilde{g}\xrightarrow[]{}1} w_{\Tilde{g}}+ew_{\Tilde{1}}) \nonumber \\
    &= \sum_{\Tilde{g}\xrightarrow[]{}g} \sum_{\Tilde{h}\in G_w}(r(\Tilde{g})+1-a_{\Tilde{g}^{-1},\Tilde{h}})\Tilde{h}-\sum_{\Tilde{g}\xrightarrow[]{}1} \sum_{\Tilde{h}\in G_w}(r(\Tilde{g})+1-a_{\Tilde{g}^{-1},\Tilde{h}})\Tilde{h}+e\sum_{\Tilde{h}\in G_w}(r(\Tilde{1})+1-a_{\Tilde{1},\Tilde{h}})\Tilde{h}. \nonumber
\end{align}
Let $m=[K_{2,w}:K_{1,v}]$. Clearly, $r(\Tilde{1})=1$ and $\sum_{\Tilde{g}\xrightarrow[]{}g}r(\Tilde{g})=e\cdot r(g)$. Now, when using these equalities together with Proposition \ref{formulas}(2) and simplifying further, we obtain

\begin{align}
   s_2(f_v'(w_g)) &= \sum_{\Tilde{h}\in G_w}\Tilde{h}[(e\cdot r(g)+m-e(a_{g^{-1},h}+k_{\Tilde{h}}))-(e\cdot1+m-e(1+k_{\Tilde{h}}))+e]  \nonumber \\
   &=e\sum_{\Tilde{h}\in G_w}(r(g)+1-a_{g^{-1},h})\Tilde{h}. \nonumber
\end{align}
Therefore, we have $e\psi_v\circ s_1=s_2\circ f_v'$. This completes the proof.
\end{proof}
The compatibility of local maps in Proposition \ref{W-compatibility} and the discussion before that gives compatibility of global maps, which is stated precisely in the following.
\begin{proposition}\label{WB-compatibility}
We have a commutative diagram of $\Bbb Z[G_2]$--modules
\[\begin{tikzcd}
W_2\arrow{r}{\gamma_2} & B_2 \\%
W_1\arrow{r}{\gamma_1}\arrow{u}{\widetilde{f'}} & B_1\arrow{u}{\widetilde{\psi'}}
\end{tikzcd},
\]
where $\widetilde{f'}$ is defined in Proposition \ref{restriction-compatibility} and $\widetilde{\psi'}$ is induced by $\widetilde{\psi}=(\Tilde{\psi_v})_{v\in S'}$ (defined in Remark \ref{Psi-definition}.)
\end{proposition}
\medskip

Next, we focus on the map $\theta_B:B:= \prod_{v\in S'}\Bbb Z[G]\to \Bbb Z[G]$ from diagram \eqref{diagram-B}. In \cite{Dasgupta-Kakde}, this map is defined componentwise as follows.
\begin{itemize}
    \item For $v\in S$, $\theta_B$ is the identity on $\Bbb Z[G]$.
    \item For $v\in S'\setminus S$, $\theta_B(x)=(\sigma_v-1)x$, for all $x\in\Bbb Z[G]$, where $\sigma_v$ is the Frobenius automorphism. (Observe that in this case $v$ is unramified.)
\end{itemize}

\begin{remark}\label{B-theta-free} Note that since $\theta_B$ is surjective, its kernel $B^\theta$ is always a free $\Bbb Z[G]$--module of rank $|S'|-1$. 
Further, note that when restricted to the category of $\Bbb Z[G]$--modules,  the functor $F(\ast)$  is equivalent to  $F_G(\ast):={\rm Hom}_{\mathbb{Z}[G]}(*,\mathbb{Z}[G])$, and  the equivalence is given by $\phi\mapsto (m\mapsto \sum_{g\in G}\phi(gm)g^{-1})$. Hence, we have a $\Bbb Z[G]$--module isomorphism 
$$\beta_G : F(\mathbb{Z}[G])\cong F_G(\Bbb Z[G])\cong  \mathbb{Z}[G], \qquad \beta_G(\phi)=\sum_{g\in G}\phi(g)g^{-1}.$$
This implies right away that the $\Bbb Z[G]$--module $F(B)$ and $F(B^\theta)$ are free of ranks $|S'|$ and $|S'|-1$, respectively.
\end{remark}
\medskip

\begin{proposition}\label{B-compatibility}
We have the following commutative diagram.
\[\begin{tikzcd}
B_1\arrow{r}{\Tilde{\psi}'}\arrow[]{d}{\theta_{B_1}} & B_2\arrow[]{d}{\theta_{B_2}} \\%
\mathbb{Z}[G_1]\arrow{r}{\Tilde{\psi}} & \mathbb{Z}[G_2]
\end{tikzcd}
\]
\end{proposition}
\begin{proof}
The commutativity of the diagram above is proved componentwise, for every $v\in S'$.  It is obvious when $v\in S$. In the case $v\in S'\setminus S$ (therefore $v$ and $w$ are unramified), we have
$$\theta_{B_2}(\Tilde{\psi}'(g))=\theta_{B_2}(\sum_{\Tilde{g} \twoheadrightarrow g}\Tilde{g})=(\sigma_w-1)\sum_{\Tilde{g} \twoheadrightarrow g}\Tilde{g}=\sum_{\Tilde{h} \twoheadrightarrow \sigma_vg}\Tilde{h}-\sum_{\Tilde{g} \twoheadrightarrow g}\Tilde{g}, \qquad \text{ for all }g\in G_1.$$
On the other hand, we also have 
$$\Tilde{\psi}(\theta_{B_1}(g))=\Tilde{\psi}((\sigma_v-1)g)=\sum_{\Tilde{h} \twoheadrightarrow \sigma_vg}\Tilde{h}-\sum_{\Tilde{g} \twoheadrightarrow g}\Tilde{g}, \qquad\text{ for all }g\in G_1.$$
This concludes the proof.
\end{proof}
Next, we finally construct the desired maps between the presentations of the Selmer modules at levels $K_2$ and $K_1$, described in Lemma \ref{quadratic-pres-lemma}(1).
\begin{theorem}\label{Z-presentations-compatible}
We have a commutative diagram of $\Bbb Z[G_2]$--modules
\[\begin{tikzcd}
F(B_2^{\theta})\arrow{d}{\pi':=F(\Tilde{\psi'})}\arrow{r} & F(V_2^{\theta})\arrow{d}{F(\Tilde{f})}\arrow{r} & Sel_S^T(K_2)\arrow{d}{res}\arrow{r} & 0 \\%
F(B_1^{\theta})\arrow{r} & F(V_1^{\theta})\arrow{r} & Sel_S^T(K_1)\arrow{r} & 0,
\end{tikzcd}
\]
where the left vertical map $\pi'$ is induced by Galois restriction, in appropriately chosen bases for the free $\Bbb Z [G_i]$--modules $F(B_i^\theta)$ of rank $(|S'|-1)$, for $i=1,2$.
\end{theorem}
\begin{proof}
Propositions \ref{WB-compatibility} and \ref{B-compatibility}, the lower left square in \eqref{diagram-B} and the right square of \eqref{W-theta-compatible}  give the following commutative diagram of $\Bbb Z[G_2]$--modules, whose vertical maps are the $\theta_W$'s and the $\theta_B$'s.
\[
 \begin{tikzcd}[row sep=1.5em, column sep = 1.5em]
    W_1\arrow[rr] \arrow[dr, swap] \arrow[dd,swap] &&
    W_2 \arrow[dd] \arrow[dr] \\
    & B_1 \arrow[rr] &&
    B_2 \arrow[dd] \\
    \Delta G_1 \arrow[rr,] \arrow[dr] && \Delta G_2 \arrow[dr] \\
    & \mathbb{Z}[G_1] \arrow[rr] \arrow[uu,leftarrow]&&
    \mathbb{Z}[G_2]
    \end{tikzcd}
\]
Now, when taking the kernels of the vertical maps above, we get a commutative diagram of $\Bbb Z[G_2]$--modules
\[
 \begin{tikzcd}
 W_2^{\theta}\arrow{r}\arrow[leftarrow]{d} & B_2^{\theta}\arrow[leftarrow]{d} \\%
 W_1^{\theta}\arrow{r} & B_1^{\theta}
    \end{tikzcd}
\]
When applying the functor $F$ to the above diagram, and then combining it with the diagram in Proposition \ref{restriction-compatibility}, we obtain the desired commutative diagram in the statement of the Theorem. The rows are exact by \cite{Dasgupta-Kakde}.\\

Now, we proceed to proving that the map $\pi'$ is induced by Galois restriction, in appropriately chosen bases for the free $\Bbb Z[G_i]$--modules $F(B_i^\theta)$, with $i=1,2$.
We know that, in the standard bases for $B_i=\Bbb Z[G_i]^{|S'|}$, the map $\Tilde{\psi}':B_1\xrightarrow[]{} B_2$ is given by $(\Tilde{\psi_v})_{v\in S'}$, where  
$$\Tilde{\psi_v}:\mathbb{Z}[G_1]\xrightarrow[]{}\mathbb{Z}[G_2], \qquad \Tilde{\psi_v}(g)=\sum_{\Tilde{g} \twoheadrightarrow g}\Tilde{g}, \qquad \text{ for all }v\in S', g\in G_1.$$ 
With notations as in Remark \ref{B-theta-free}, we claim that we have the following commutative diagrams
\[
 \begin{tikzcd}
F(\mathbb{Z}[G_2])\arrow{r}{F(\Tilde{\psi_v})}\arrow{d}{\beta_{G_2}} & F(\mathbb{Z}[G_1])\arrow{d}{\beta_{G_1}}  \\%
\mathbb{Z}[G_2]\arrow{r}{\pi} & \mathbb{Z}[G_1],
    \end{tikzcd}
\]
for all $v\in S'$, where the lower horizonthal map is Galois restriction. Indeed, observe that we have
$$\pi(\beta_{G_2}(\phi))=\pi(\sum_{\Tilde{g}\in G_2}\phi(\Tilde{g})\Tilde{g}^{-1})=\sum_{g\in G_1}\sum_{\Tilde{g}\xrightarrow[]{}g}\phi(\Tilde{g})g^{-1}.$$
On the other hand, we also have 
$$\beta_{G_1}(F(\Tilde{\psi_v})(\phi))=\sum_{g\in G_1}F(\Tilde{\psi_v})(\phi)(g)g^{-1}=\sum_{g\in G_1}\phi(\Tilde{\psi}(g))g^{-1}=\sum_{g\in G_1}\sum_{\Tilde{g}\xrightarrow[]{}g}\phi(\Tilde{g})g^{-1},$$
which proves the claim. Consequently, the map $F(\Tilde{\psi'}): F(B_2)\to F(B_1)$ is given by Galois restriction in the standard bases $\mathbf e^i$ for $F(B_i)=F(\Bbb Z[G_i])^{|S'|}$, for $i=1,2$. Since the diagram of $\Bbb Z[G_2]$--modules
\[
 \begin{tikzcd}
F(B_2)\arrow{r}{F(j_2)}\arrow{d}{F(\tilde\psi')} & F(B_2^\theta)\arrow{d}{F(\tilde\psi')}  \\%
F(B_1)\arrow{r}{F(j_1)} & F(B_1^\theta),
    \end{tikzcd}
\]
commutes (here $j_i:B_i^\theta\to B_i$ are inclusions and $F(j_i)$ are surjective),
this implies that $F(\tilde\psi')$ is given by Galois restriction in the sets of $\Bbb Z[G_i]$--generators $\mathbf e^{i,\theta}:=F(j_i)(\mathbf e^i)$ of $F(B_i^\theta)$, for $i=1,2$. This means that 
$$F(\tilde\psi')(\sum_{k=1}^{|S'|}\alpha_k\cdot \mathbf e^{2,\theta}_k)=\sum_{k=1}^{|S'|}\pi(\alpha_k)\cdot \mathbf e^{1, \theta}_k, \qquad\text{ for all }\alpha_1, \dots, \alpha_{|S'|}\in \Bbb Z[G_2].$$
Now, the $\Bbb Z[G_2]$--linearity of $F(\tilde\psi'): F(B_2^\theta)\to F(B_1^\theta)$, implies that this map is given by Galois restriction in any bases of the type $\mathbf f^{2, \theta}$ and $F(\tilde\psi')(\mathbf f^{2, \theta})$ of 
$F(B_2^\theta)$ and $F(B_1^\theta)$, respectively.

Now, to conclude the proof, note that since the map $F(\tilde\psi')$ is surjective, the modules $F(B_i^\theta)$ are $\Bbb Z[G_i]$--free of common rank $|S'|-1$, for $i=1,2$, and the total ring of fractions of $\Bbb Z[G_1]$ is a finite direct sum of fields, for any $\Bbb Z[G_2]$--basis $\mathbf f^{2,\theta}$ of $F(B_2^\theta)$, the set $F(\tilde\psi')(\mathbf f^{2, \theta})$ is a set of $\Bbb Z[G_1]$--generators and therefore a $\Bbb Z[G_1]$--basis for $F(B_1^\theta)$. In the bases $\mathbf f^{2,\theta}$ and $F(\tilde\psi')(\mathbf f^{2, \theta})$, the map $F(\tilde\psi')$ is given by Galois restriction.
\end{proof}
\medskip

In what follows, if $G$ is a finite group, $p$ is a prime number and $M$ is a $\Bbb Z[G]$--module, we let
$$F_p(M):= {\rm Hom}_{\mathbb{Z}_p}(M\otimes_{\Bbb Z}\Bbb Z_p, \Bbb Z_p),$$
viewed as a $\Bbb Z_p[G]$--module with the covariant $G$--action. It is easily seen that for all free $\mathbb{Z}$-modules of finite rank $M$ one has an isomorphism of $\Bbb Z_p[G]$--modules
$$F_p(M)\simeq F(M)\otimes_{\Bbb Z}\Bbb Z_p.$$

\begin{theorem}\label{QP-compatibility} The following hold for $(K_2/K_1/F, S, T, S')$ as above and any prime $p>2$.
\begin{enumerate}
\item
If $S_{ram}(K_2/F)\subseteq S\cup T$, we have a commutative diagram of $\mathbb{Z}_p[G_2]$- modules
\[\begin{tikzcd}
F_p(B_2^{\theta})\arrow{d}{\pi'}\arrow{r} & F_p(V_2^{\theta})\arrow{d}\arrow{r} & Sel_S^T(K_2)_p\arrow{d}{res}\arrow{r} & 0 \\%
F_p(B_1^{\theta})\arrow{r} & F_p(V_1^{\theta})\arrow{r} & Sel_S^T(K_1)_p\arrow{r} & 0,
\end{tikzcd}
\]
where the left vertical map $\pi'$ is induced by Galois restriction.
\item If, in addition, $S_p(F)\cap S_{ram}(K_2/F)\subseteq S$, then we have
$$F_p(B_2^\theta)\simeq F_p(V_2^\theta)\simeq\Bbb Z_p[G_2]^k, \qquad  F_p(B_1^\theta)\simeq F_p(V_1^\theta)\simeq\Bbb Z_p[G_1]^k,$$
where $k=|S'|-1$ and appropriate bases for $F_p(B_2^\theta)$ and $F_p(B_1^\theta)$ can be chosen so that $\pi'$ is induced by the usual Galois restriction $\Bbb Z_p[G_2]\twoheadrightarrow\Bbb Z_p[G_1]$
in those bases.
\end{enumerate}
\end{theorem}
\begin{proof} Part (1)  follows directly from Theorem \ref{Z-presentations-compatible} by tensoring with $\Bbb Z_p$.
Part (2) is a consequence of part (1) and Lemma \ref{quadratic-pres-lemma}(2).
\end{proof}

The main consequence of the theorem above, which will play an important role in the the proof of the Equivariant Main Conjecture, is the following.
\begin{corollary}\label{covariance}
Let $H_0$ be an abelian CM extension of $F$ with Galois group $G_0$. Let $S_0, T_0$ be two finite, disjoint sets  places in $F$ which satisfy hypotheses ${\rm Hyp}(H_0/F)_p$. If $I$ is any subgroup of $G_0$ which does not contain the complex conjugation of $H_0$ and $R$ is a quotient of $\mathbb{Z}_p[G_0/I]$ such that $Fitt_R(Sel_{S_0}^{T_0}(H_0 ^I)_R)$ is generated by a nonzero divisor, then
$$Sel_{S_0}^{T_0}(H_0)_R\simeq  Sel_{S_0}^{T_0}(H_0 ^I)_R,$$
where the isomorphism is induced by the Selmer restriction map.
\end{corollary}
\begin{proof}
Since the complex conjugation of $H_0$ is not in $I$, $H_0 ^I$ is also a CM-field. Then, by applying Theorem \ref{QP-compatibility} and tensoring the commutative diagram in its statement with $R$, we obtain a commutative diagram of $R$--modules
with exact rows
\[\begin{tikzcd}
R^k\arrow[equal]{d}\arrow{r}{\alpha} & R^k\arrow{d}{\gamma}\arrow{r} & Sel_{S_0}^{T_0}(H_0)_R\arrow{d}{res_R}\arrow{r} & 0\\%
R^k\arrow{r}{\beta} & R^k\arrow{r} & Sel_{S_0}^{T_0}(H_0^I)_R\arrow{r} & 0,
\end{tikzcd}
\]
where $k$ is the common rank of the two quadratic presentations in Theorem \ref{QP-compatibility}. The left vertical map is the identity because in Theorem \ref{QP-compatibility} that map is induced by Galois restriction.

If $A$, $B$ and $C$ are the matrices corresponding to the maps $\alpha$, $\beta$ and $\gamma$, respectively,  in the standard bases, we have $B=CA$. So, $det(B)=det(C)\cdot det(A)$. However,  $det(B)$ is a generator of $Fitt_R(Sel_{S_0}^{T_0}(H_0 ^I)_R)$. Therefore, $det(B)$ is a non--zero divisor in $R$. Therefore, so is $det(C)$ and $det(A)$, which implies that $\alpha$,$\beta$ and $\gamma$ are injective. Therefore, the rows of the above diagram are exact. Now, Lemma \ref{surjective-res-Selmer} implies that $res_R$ is surjective. Hence, by an application of the snake lemma in the diagram above, $\gamma$ is also surjective. Therefore, $\gamma$ is an isomorphism. Consequently the map $res_R$ is an isomorphism, which concludes the proof.
\end{proof}

\begin{remark}\label{character group rings remark} The typical quotients $R$ of $\Bbb Z_p[G_0/I]$ to which the above statement applies are of the form $R:=R_\Psi$, where
$\Psi\subseteq\widehat{G_0/I}$ and
$$R_\Psi:={\rm Im}(\Bbb Z_p[G_0/I]\overset{{\rm ev}_{\Psi}}\longrightarrow \bigoplus_{\psi\in\Psi}\Bbb Z_p[\psi]), \qquad {\rm ev}_{\Psi}(x)=(\psi(x))_{\psi\in\Psi}.$$
The rings $R_\Psi$ are called character group rings in \cite{Dasgupta-Kakde}.
\end{remark}

\section{The Equivariant Main Conjecture}

Below, we fix an abelian, CM extension $H/F$ of a totally real number field $F$, of Galois group $G$. We fix a prime $p>2$ and, as usual, let $H_\infty$ and $F_\infty$ be the cyclotomic $\Bbb Z_p$--extensions of $H$ and $F$, respectively.
As before, $\mathcal G:=G(H_\infty/F)$ and $G_n:=G(H_n/F)$, for $n\geq 0$. We let $G_p$ and $G'$ denote the $p$--Sylow subgroup of $G$ and its complement, respectively. The same notational convention applies to $G_{n, p}$ and $G'_n$, for all $n$.\\

We consider two finite sets of primes $S_0$ and $T_0$ in $F$, such that
$$S_0\cap T_0=S_{ram}(H_\infty/F)\cap S_0=S_{ram}(H_\infty/F)\cap T_0=\emptyset, \qquad T_0\ne \emptyset.$$
For simplicity, we will denote $S_{p\infty}(F):=S_\infty(F)\cup S_p(F)$. Further, we let
\begin{equation}\label{Sigma-Sigma'}\Sigma:=S_{p\infty}(F)\cup S_0, \qquad \Sigma':=(S_{ram}(H_\infty/F)\setminus S_{p\infty}(F))\cup T_0.\end{equation}
Note that  $(\Sigma, \, \Sigma')$ satisfy the set of hypotheses ${\rm Hyp}(H_\infty/F)_p$ and so do $(\Sigma'\setminus J,\, \Sigma\cup J)$, for all the subsets $J\subseteq (S_{ram}(H_\infty/F)\setminus S_{p\infty}(F)).$
\\

The following is an immediate consequence of Theorem 3.3 in \cite{Dasgupta-Kakde}.
\begin{theorem}\label{empty-J-theorem}
For all $(H/F, S_0, T_0)$ as above and all $n\gg 0$, we have an equality
$$Fitt_{\mathbb{Z}_p[G_n]^-}(Sel_{\Sigma}^{\Sigma'}(H_n)_p^-)=( \Theta_{\Sigma}^{\Sigma'}(H_n/F)).$$
\end{theorem}
\begin{proof} Take any $n>0$ large enough so that $S_{ram}(H_\infty/F)=S_{ram}(H_n/F)$. For any such $n$ and $S_0=\emptyset$, the statement above is precisely Theorem 3.3 in \cite{Dasgupta-Kakde}. If $S_0\ne\emptyset$, by
Remark \ref{larger-S-remark} we have an
exact sequence of $\Bbb Z_p[G_n]^-$--modules
$$0\longrightarrow Y_{S_0}(H_n)^{\ast, -}\longrightarrow Sel_{\Sigma}^{\Sigma'}(H_n)_p^-\longrightarrow Sel_{\Sigma\setminus S_0}^{\Sigma'}(H_n)_p^-\longrightarrow 0.$$
Since the right--most non--zero module above is quadratically presented, Lemma \ref{Fitt-ses}, Lemma \ref{Fitt-Y_S-dual} and Theorem 3.3 in \cite{Dasgupta-Kakde} give equalities
$$Fitt_{\mathbb{Z}_p[G_n]^-}(Sel_{\Sigma}^{\Sigma'}(H_n)_p^-)= ( \Theta_{\Sigma\setminus S_0}^{\Sigma'}(H_n/F)\cdot\prod_{v\in S_0} (1-\sigma_v^{-1}))=( \Theta_{\Sigma}^{\Sigma'}(H_n/F)),$$
which concludes the proof.
\end{proof}

The next results aim at strengthening the Theorem above. In all the considerations which follow,  it is very important to keep in mind that, unlike \cite{Dasgupta-Kakde}, which uses the contravariant Galois action on duals of Galois modules, we use the covariant Galois action.

\begin{theorem}\label{partial}
For all $(H/F, S_0, T_0)$ as above, all $J\subseteq (S_{ram}(H_\infty/F)\setminus S_{p\infty}(F))$, and all $n\gg 0$, we have equalities of $\Bbb Z_p[G_n]^-$--ideals
$$Fitt_{\mathbb{Z}_p[G_n]^-}(Sel_{\Sigma\cup J}^{\Sigma'\setminus J}(H_n)_p^-)=( \Theta_{\Sigma\cup J}^{\Sigma'\setminus J}(H_n/F) \cdot u_n(J)),$$
where $u_n(J)\in \mathbb{Q}_p[G_n]^-$ and $\psi(u_n(J))\in\Bbb Z_p[\psi]^\times$, for all odd characters ${\psi}\in\widehat{G_n}$.
\end{theorem}
\begin{proof} We fix the prime $p>2$ and the totally real number field $F$ and prove the statement above by induction on $|J|$, for all $(H/F, S_0, T_0)$ as above.
The base case  $J=\emptyset$ holds for $u_n(\emptyset)=1$, for all $n\gg 0$, as Theorem \ref{empty-J-theorem} shows. Now let us assume that the result holds for all $(H/F, S_0, T_0, J)$, with $J$ of a fixed size. \\

For fixed data as above, let  $v\in (S_{ram}(H_\infty/F)\setminus S_{p\infty}(F))\setminus J$. We need to prove the result for the data $(H/F, S_0, T_0, J_v:= J\cup \{v\})$.\\

By \eqref{semilocal}, we know that it suffices  to prove the desired equality after tensoring with
 $$R_n^{\chi}:=\mathbb{Z}_p[\chi][G_{n,p}],$$
 for all odd characters $\chi\in\widehat{G'}$.  Let us fix such an odd character $\chi$.
Let $I_v$ be the ramification group of $v$ in $G_n$. Note that $I_v$ is independent of $n$ and can be viewed as a subgroup of $G$, as $H_\infty/H$ is unramified at $v$. Let $I_v\cong I_{v,p}\times I'_v$, where $I_{v, p}$ is the $p$--Sylow subgroup of $I_v$. Obviously, we have $I_v'\subseteq G'$. We divide the proof of the inductive step in two cases, depending on whether
$\chi$ is trivial on $I_v'$ or not.\\

{\bf $\bullet$ Case (1):} $\chi(I_v')\neq\{1\}.$  Fix $n\gg 0$.  In our context, exact sequence (40) in \cite{Dasgupta-Kakde} tensored with $R_n^\chi$ gives an ismorphism of $R_n^\chi$--modules
\begin{equation}\label{Sel-Jv-J}  Sel_{\Sigma\cup J}^{\Sigma'\setminus J_v}(H_n)_{R_n^{\chi}}\cong Sel_{\Sigma\cup J}^{\Sigma'\setminus J}(H_n)_{R_n^{\chi}}.\end{equation}
Since the right module is a quadratically presented $R_n^{\chi}$--module (see Lemma \ref{quadratic-pres-lemma}(2)), so is the left module.
Now, by \eqref{enlarge-S-Selmer-sequence}, for all $\chi\in\widehat{G'}$, we have an exact sequence of $R_n^\chi$--modules
\begin{equation}\label{case1-ses} 0\xrightarrow[]{} Y_v(H_n)^\ast_{R_n^{\chi}} \xrightarrow[]{} Sel_{\Sigma\cup J_v}^{\Sigma'\setminus J_v}(H_n)_{R_n^{\chi}}\xrightarrow[]{}Sel_{\Sigma\cup J}^{\Sigma'\setminus J_v}(H_n)_{R_n^{\chi}}\xrightarrow[]{}0. \end{equation}
Consequently, Lemma \ref{Fitt-ses} applied to the exact sequence above, combined with \eqref{Sel-Jv-J}, gives
\begin{equation}\label{case1-ses-Fitt}
{\rm Fitt}_{R_n^\chi} Sel_{\Sigma\cup J_v}^{\Sigma'\setminus J_v}(H_n)_{R_n^{\chi}} = ({\rm Fitt}_{R_n^\chi} Sel_{\Sigma\cup J}^{\Sigma'\setminus J}(H_n)_{R_n^\chi})\cdot ({\rm Fitt}_{R_n^\chi} Y_v(H_n)^\ast_{R_n^{\chi}}).
\end{equation}
However, by Lemma \ref{Fitt-Y_S-dual} combined with our assumption $\chi(I_v')\ne \{1\}$, we have
$${\rm Fitt}_{R_n^{\chi}} Y_v(H_n)_{R_n^{\chi}}^* =(1-\sigma_v^{-1}, 1-\tau \mid  \tau\in I_v)_{R_n^\chi}=R_n^\chi,$$
because for any $\tau\in I_v'$, such that $\chi(\tau)\ne 1$, we have $(1-\chi(\tau))\in (R_n^\chi)^\times.$ Consequently, the induction hypothesis combined with \eqref{case1-ses-Fitt} gives an equality
$${\rm Fitt}_{R_n^{\chi}}(Sel_{\Sigma\cup J_v}^{\Sigma'\setminus J_v}(H_n)_{R_n^{\chi}})=( \Theta_{\Sigma\cup J}^{\Sigma'\setminus J}(H_n/F) \cdot u_n(J)^\chi)_{R_n^\chi}, $$
where $u_n(J)^\chi$ is as stated in the theorem. However, we have an equality in $R_n$:
$$\Theta_{\Sigma\cup J}^{\Sigma'\setminus J}(H_n/F)\cdot (1-\sigma_v^{-1} e_v)=\Theta_{\Sigma\cup J_v}^{\Sigma'\setminus J_v}(H_n/F)\cdot (1-\sigma_v^{-1} e_v Nv).$$
On the other hand, since $\chi(I_v')\ne\{1\}$, we have $e_v R_n^\chi=0$. Consequently, when combining the last two displayed equalities, we obtain the desired result for $J_v$, with $u_n(J_v)^\chi=u_n(J)^\chi$.
\begin{remark}\label{case1-remark} Note from the considerations above that in Case 1, one can pick $u_n(J)$ and $u_n(J_v)$ such that
$$u_n(J_v)^\chi=u_n(J)^\chi, \qquad \text{ for all }n\gg 0.$$
So, in particular, if 
one can take $u_n(J)^\chi=1$ (e.g. $u_n(\emptyset)=1$), then one can take $u_n(J_v)^\chi=1$, for all $n\gg 0$.
\end{remark}
\medskip

{\bf $\bullet$ Case (2):} $\chi(I_v')=\{1\}$. As above, fix $n\gg 0$. For appropriate sets of places, the short exact sequence (40) in \cite{Dasgupta-Kakde} reads
$$0\longrightarrow Sel_{\Sigma\cup J}^{\Sigma'\setminus J_v}(H_n)_{R_n^{\chi}}\longrightarrow Sel_{\Sigma\cup J}^{\Sigma'\setminus J}(H_n)_{R_n^{\chi}}\longrightarrow R_n^{\chi}/(\tau_v-1,\sigma_v-Nv)\longrightarrow  0,$$
where $\tau_v$ is a generator of the cyclic group $I_{v,p}$ (recall that $v\nmid p$, therefore $I_{v, p}$ is cyclic) and $\sigma_v$ is a choice of Frobenius at $v$ in $G_n$. Consequently, equation (42) of \cite{Dasgupta-Kakde} reads
$$Fitt_{R_n^{\chi}}(Sel_{\Sigma\cup J}^{\Sigma'\setminus J_v}(H_n)_{R_n^{\chi}})(\sigma_v-Nv)=Fitt_{R_n^{\chi}}(Sel_{\Sigma\cup J}^{\Sigma'\setminus J}(H_n)_{R_n^{\chi}})(NI_{v},\sigma_v-Nv),$$
where $NI_v:=\sum_{\tau\in I_v}\tau$.
By the induction hypothesis  we obtain
$$Fitt_{R_n^{\chi}}(Sel_{\Sigma\cup J}^{\Sigma'\setminus J_v}(H_n)_{R_n^{\chi}})(\sigma_v-Nv)=(\Theta_{\Sigma\cup J}^{\Sigma'\setminus  J_v}(H_n/F)\cdot u_1)_{R_n^\chi}(1-e_v\sigma_v^{-1}Nv)(NI_v ,\sigma_v-Nv),$$ where $u_1:=u_n(J)^\chi$ is as stated in the theorem. Now, since $(\sigma_v-Nv)^{-1}\in {\rm Frac}(R_n)$, an easy calculation in ${\rm Frac}(R_n)$ gives equalities of $R_n^\chi$--ideals
$$(\sigma_v-Nv)^{-1}(1-e_v\sigma_v^{-1}Nv)(NI_v ,\sigma_v-Nv)R^\chi_n= (NI_v, 1-\sigma_v^{-1} e_v)R^\chi_n=(NI_v, 1-\sigma_v^{-1}e_v)R_n^\chi.$$

The last two displayed equalities combined with the surjection in \eqref{case1-ses} give
\begin{eqnarray}\label{Case 2-inclusion}
\nonumber {\rm Fitt}_{R_n^{\chi}}(Sel_{\Sigma\cup J_v}^{\Sigma'\setminus J_v}(H_n)_{R_n^{\chi}}) &\subseteq & {\rm Fitt}_{R_n^{\chi}}(Sel_{\Sigma\cup J}^{\Sigma'\setminus J_v}(H_n)_{R_n^{\chi}}=\\
\label{Fitt-Selmer-Jv-J} &=&(\Theta_{\Sigma\cup J}^{\Sigma'\setminus  J_v}(H_n/F)\cdot u_1)_{R_n^\chi}(NI_{v}, 1-\sigma_v^{-1} e_{v}).
\end{eqnarray}
However, the $R_n^\chi$--module $Sel_{\Sigma\cup J_v}^{\Sigma'\setminus J_v}(H_n)_{R_n^{\chi}}$ is quadratically presented. Therefore, its Fitting ideal is principal. Denote by $\theta$ a generator. The last inclusion shows that we have
\begin{eqnarray}\label{theta}\theta&=&\Theta_{\Sigma\cup J}^{\Sigma'\setminus J_v}(H_n/F)\cdot u_1\cdot (aNI_v+b(1-\sigma_v^{-1} e_v))\\
\nonumber &=& \Theta_{\Sigma\cup J}^{\Sigma'\setminus J_v}(H_n/F)\cdot u_1\cdot (aNI_v)+ \Theta_{\Sigma\cup J_v}^{\Sigma'\setminus J_v}(H_n/F)\cdot (bu_1)  \quad\text{  for some $a,b\in R_n^{\chi}$.}\end{eqnarray}
 \medskip

 Now let us consider the following sets of characters of $G_n$:
\begin{eqnarray}\label{Psi-Omega}\Psi'=\Psi'(G_n):=\{\psi\in\widehat{G_n/I_v}\mid \psi\vert_{G'}=\chi\}, \\
\nonumber \Omega=\Omega(G_n):=\{\psi\in\widehat{G_n} \mid \psi\vert_{G'}=\chi,  \psi(G_{v_0,n})\ne 1,\, \forall v_0\in \Sigma\cup J_v \},\\
\nonumber \Psi:=\Psi'\cap\Omega, \qquad \Omega^c:= \{\psi\in\widehat{G_n}\mid  \psi|_{G'}=\chi \}\setminus \Omega.\end{eqnarray}
Note that, when extended by $\Bbb Q_p$--linearity, the evaluation maps defined in Remark \ref{character group rings remark} give an injective $\Bbb Q_p$--algebra morphism
$${\rm ev}^\chi:={\rm ev}_\Psi\oplus{\rm ev}_{\Omega\setminus\Psi}\oplus{\rm ev}_{\Omega^c} : \Bbb Q_p (R_n^\chi) \hookrightarrow \Bbb Q_p(R_n)_\Psi\oplus \Bbb Q_p(R_n)_{\Omega\setminus\Psi}\oplus \Bbb Q_p(R_n)_{\Omega^c}.$$
Next, we will study the image of $\theta$ via this evaluation map, componentwise.\\

{\bf Step 1:} ${\rm ev}_\Psi(\theta)$. Note that since $\chi$ is odd and $\chi(I_v')=1$, $I_v$ cannot contain the complex conjugation morphism. Therefore $H^{I_v}/F$ is a CM, abelian extension where $v$ is unramified. Clearly, $H_n^{I_v}$ is the $n$--the layer in the cyclotomic $\Bbb Z_p$--extension of  $H^{I_v}$. Therefore, we may apply the induction hypothesis to the data
$$(H^{I_v}/F, S_0', T_0', J), \quad S_0'=S_0\cup\{v\}, \quad T_0'=T_0\cup (S_{ram}(H_\infty/F)\setminus S_{ram}(H^{I_v}_\infty/F).$$
As the sets $\Sigma$ and $\Sigma'$ for this new data become $\Sigma\cup\{v\}$ and $\Sigma\setminus\{v\}$, after tensoring
with $R_\Psi$ over $\Bbb Z_p[G_n/I_v]$, we obtain equalities
$${\rm Fitt}_{R_{\Psi}}(Sel_{\Sigma\cup J_v}^{\Sigma'\setminus J_v}(H_n^{I_v})_{R_{\Psi}})=(\Theta_{\Sigma\cup J_v}^{\Sigma'\setminus J_v}(H_n^{I_v}/F)\cdot u_2')_{R_\Psi},\qquad\text{ for all }n\gg0,$$
where  $u_2' \in\Bbb Q_p[G_n/I_v]^\chi$, such that $\psi(u_2')\in\Bbb Z_p[\psi]^\times$, for all $\psi\in\widehat{G_n/I_v}$ with $\psi\mid_{G'}=\chi$. Note that, for simplicity, we suppress in the notation the dependence of $u_2'$ of $n$ and $J$.

Since $\Psi\subseteq \Omega$, the element ${\rm ev}_\Psi(\Theta_{\Sigma\cup J_v}^{\Sigma'\setminus J_v}(H_n^{I_v}/F)\cdot u_2')$ is a non--zero divisor in $R_\Psi$. (See Lemma \ref{vanishing-of-Theta}.) Therefore, by Corollary \ref{covariance}, we have an $R_\Psi$--module isomorphism
$$Sel_{\Sigma\cup J_v}^{\Sigma'\setminus J_v}(H_n)_{R_{\Psi}}\cong Sel_{\Sigma\cup J_v}^{\Sigma'\setminus J_v}(H_n^{I_v})_{R_{\Psi}}.$$
Therefore, we have an equality of $R_\Psi$--ideals
$$\theta R_\Psi={\rm Fitt}_{R_{\Psi}}(Sel_{\Sigma\cup J_v}^{\Sigma'\setminus J_v}(H_n)_{R_{\Psi}})=(\Theta_{\Sigma\cup J_v}^{\Sigma'\setminus J_v}(H_n/F)\cdot u_2)_{R_\Psi},$$
for any lift $u_2$ in $\Bbb Q_p[G_n]^\chi$ of $u_2'$. Therefore, by multiplying $u_2$ and $u_2'$ with properly chosen units in the appropriate rings, we have the following equality
\begin{equation}\label{ev Psi}{\rm ev}_\Psi(\theta)={\rm ev}_{\Psi}( \Theta_{\Sigma\cup J_v}^{\Sigma'\setminus  J_v}(H_n/F)\cdot u_2),\end{equation}
where $u_2\in\Bbb Q_p[G_n]^\chi$ and $\psi(u_2)=\psi(u_2')\in\Bbb Z_p[\psi]^\times$, for all $\psi\in\Psi$.

{\bf Step 2:} ${\rm ev}_{\Omega\setminus\Psi}(\theta)$. Note that since $\psi(NI_v)=0$, for all $\psi\in\Omega\setminus\Psi$, equalities \eqref{theta} lead to
\begin{equation}\label{ev Omega minus Psi} {\rm ev}_{\Omega\setminus\Psi}(\theta)=  {\rm ev}_{\Omega\setminus\Psi}( \Theta_{\Sigma\cup J_v}^{\Sigma'\setminus J_v}(H_n/F)\cdot (bu_1) )\end{equation}

{\bf Step 3:} ${\rm ev}_{\Omega^c}(\theta)$. Lemma 3.2 in \cite{Dasgupta-Kakde} implies that, in this case, we have
$$\Theta_{\Sigma\cup J_v}^{\Sigma'\setminus J_v}(H_n/F) R_{\Omega^c}= {\rm Fitt}_{R_{\Omega^c}}(Sel_{\Sigma\cup J_v}^{\Sigma'\setminus J_v}(H_n)_{R_{\Omega^c}})=0.$$
Consequently, we have an equality in $R_{\Omega^c}$:
\begin{equation}\label{ev Omega complement}
{\rm ev}_{\Omega^c}(\theta)={\rm ev}_{\Omega^c}(\Theta_{\Sigma\cup J_v}^{\Sigma'\setminus J_v}(H_n/F))=0.
\end{equation}
\\

If we combine \eqref{ev Psi}, \eqref{ev Omega minus Psi}, \eqref{ev Omega complement} with the fact that ${\rm ev}^\chi$ is injective, we get
\begin{equation}\label{new-theta}\theta=\Theta_{\Sigma\cup J_v}^{\Sigma'\setminus J_v}(H_n/F)\cdot u_3, \qquad\text{ where } u_3=(e_v\cdot u_2+(1-e_v)\cdot bu_1) \in \Bbb Q_p[G_n]^\chi.\end{equation}
In order to complete the proof, we need to show that $\psi(u_3)\in\Bbb Z_p[\psi]^\times$, for all $\psi\in\widehat G_n$, such that $\psi\vert_{G'}=\chi$. Now, observe that for any such character $\psi$, we have
$$\psi(u_3)=\begin{cases} \psi(u_2), &\text{ if }\psi(I_v)=1;\\
\psi(b)\psi(u_1), &\text{ if }\psi(I_v)\ne 1.
\end{cases}$$
Since $\psi(u_1), \psi(u_2) \in\Bbb Z_p[\psi]^\times$, we need to show that
$$\psi(b)\in\Bbb Z_p[\psi]^\times,$$
 for all $\psi$ with $\psi\vert_{G'}=\chi$ and $\psi(I_v)\ne 1$. Note that if $|I_{v, p}|=1$, there are no such characters $\psi$, as $\psi(I_v)=\chi(I_v')=1$, by assumption. So, we assume that
$$|I_{v, p}|>1$$
and we proceed to showing that, under this hypothesis, the much stronger statement
$$b\in (R_n^\chi)^\times$$
holds.
First, note that since $R_n^\chi$ is a local ring of maximal ideal $\frak m_n^\chi=(p, 1-\sigma\mid \sigma\in G_{n, p})$, we have $b\in (R_n^\chi)^\times$ if and only if there exists $\psi\in\widehat{G_n}$, with $\psi\vert_{G'}=\chi$,  such that $\psi(b)\in\Bbb Z_p[\psi]^\times$.  Take $n\gg 0$ large enough so that $\Psi\ne\emptyset$ and Take $\psi\in\Psi$. First, note that
$$u_v^\psi:=\psi(1-\sigma_v^{-1} e_v)=(1-\psi(\sigma_v^{-1}))\in\Bbb Z_p[\psi]\setminus 0, \quad\text{ and } (u_v^\psi)^{-1}|I_{v,p}|\in\frak m_{\Bbb Z_p[\psi]},$$
where $\frak m_{\Bbb Z_p[\psi]}$ is the maximal ideal in $\Bbb Z_p[\psi]$.
Now, use \eqref{theta} and \eqref{ev Psi}, the equality
$$\psi(\Theta_{\Sigma\cup J}^{\Sigma'\setminus J_v}(H_n/F))=\psi(\Theta_{\Sigma\cup J_v}^{\Sigma'\setminus J_v}(H_n/F))\cdot (u_v^\psi)^{-1},$$
and the fact that $\psi(\Theta_{\Sigma\cup J_v}^{\Sigma'\setminus J_v}(H_n/F))$ is a non-zero divisor to obtain an equality in $\Bbb Z_p[\psi]$:
$$\psi( u_2)= \psi( u_1)\psi(a)|I_v'|\cdot (|I_{v, p}|(u_v^\psi)^{-1})+ \psi(b)\psi(u_1).$$
Since $\psi(u_1), \psi(u_2)\in\Bbb Z_p[\psi]^\times$ and $|I_{v,p}|(u_v^\psi)^{-1}\in \frak m_{\Bbb Z_p[\psi]}$,  we must have $\psi(b)\in \Bbb Z_p[\psi]^\times$, which concludes the proof of the theorem.
\end{proof}

\medskip

Next, under some additional assumption on the extension $H/F$, we remove the term $u$ in the above theorem. We let $F_\infty$ denote the cyclotomic $\Bbb Z_p$--extension of $F$.
\begin{theorem}\label{full}
Assume that data $(H/F, S_0, T_0)$ is as above and $H\cap F_\infty=F$. Then, for all sets $J\subseteq (S_{ram}(H_\infty/F)\setminus S_{p\infty}(F))$, we have equalities of $\Bbb Z_p[G_n]^-$--ideals
$$Fitt_{\mathbb{Z}_p[G_n]^-}(Sel_{\Sigma\cup J}^{\Sigma'\setminus J}(H_n)_p^-)=( \Theta_{\Sigma\cup J}^{\Sigma'\setminus J}(H_n/F)), \qquad \text{for all }n\gg0.$$
\end{theorem}
\begin{proof}  In what follows, we let $\Gamma=G(F_\infty/F)$. Since, $F_{\infty}/F$ and $H/F$ are disjoint, Galois restriction gives isomorphisms  $G(H_\infty/H)\simeq\Gamma$ and $\mathcal{G}=Gal(H_{\infty}/F)\cong G\times \Gamma$. Therefore,  we have isomorphisms of topological, compact $\Bbb Z_p$--algebras
$$\mathbb{Z}_p[[\mathcal{G}]]\simeq \mathbb{Z}_p[G][[\Gamma]]=\mathbb{Z}_p[G'][G_p][[\Gamma]]\simeq \bigoplus_{[\chi]\in[\widehat{G'}]}R_{\infty}^{\chi},$$
 where the last isomorphism given by the direct sum of the usual $\chi$--evaluation maps and  and $R_{\infty}^{\chi}=\mathbb{Z}_p[\chi][G_p][[\Gamma]]$.  Obviously, we have
 $$R_\infty^\chi\simeq \varprojlim_n R_n^\chi, $$
 where the rings $R_n^\chi=\Bbb Z_p[\chi][G_{n,p}]$, as before.\\


The strategy of proof is the same as that of Theorem \ref{partial}. That is, we fix the field $F$ and proceed by induction on $|J|$, for all data $(H/F, S_0, T_0)$ as above. The base case $J=\emptyset$ is Theorem \ref{empty-J-theorem} above. The proof of the inductive step, going from a given $J$ to $J_v=J\cup\{v\}$, with $v\in (S_{ram}(H_\infty/F)\setminus S_{p\infty}(F))\setminus J$, is again done on a character--by--character basis, with respect to all the odd characters $\chi\in\widehat{G'}$. For a given such  character $\chi$, the proof is again divided into two cases, {\bf Case 1:} $\chi\vert_{I_v'}\ne 1$ and {\bf Case 2:} $\chi\vert_{I_v'}= 1$. According to Remark \ref{case1-remark}, in {\bf Case 1}, the implication in the inductive step holds. So, we need to prove this implication
in {\bf Case 2}, under the induction hypothesis that $u_n(J)=1$, for $n\gg 0$, with notations as in the proof of Theorem \ref{partial}. \\

From the proof of {\bf Case 2} in Theorem \ref{partial}, since under our new induction hypothesis $u_1=u_2=1$, we have  $u_3=e_{v,p}+(1-e_{v,p})b$. Therefore, for all $n\gg 0$, we have an equality
$${\rm Fitt}_{R_n^{\chi}}(Sel_{\Sigma\cup J_v}^{\Sigma'\setminus J_v}(H_n)_{R_n^{\chi}}^-)=( \Theta_{\Sigma\cup J_v}^{\Sigma'\setminus J_v}(H_n/F)(e_{v,p}+(1-e_{v,p})b_n))_{R_n^\chi},$$
with $b_n\in R_n^\chi$, such that
\begin{equation}\label{vn-psi-eval}\psi(e_{v,p}+(1-e_{v,p})b_n))\in\Bbb Z_p[\psi]^\times,\qquad\text{  for all $\psi\in\widehat{G_n}$.}\end{equation}
 Further, from \eqref{theta} and \eqref{new-theta}, there exist $a_n\in R_n^\chi$, such that
$$
\Theta_{\Sigma\cup J_v}^{\Sigma'\setminus J_v}(H_n/F)(e_{v,p}+(1-e_{v,p})b_n)= \Theta_{\Sigma\cup J}^{\Sigma'\setminus J_v}(H_n/F)\cdot  (a_nNI_{v, p}) + \Theta_{\Sigma\cup J_v}^{\Sigma'\setminus J_v}(H_n/F)\cdot b_n.
$$
Since ${\rm ev}_{\Omega(G_n)}( \Theta_{\Sigma\cup J}^{\Sigma'\setminus J_v}(H_n/F))$ is not a zero-divisor in $R_{\Omega(G_n)}$, this implies that
\begin{equation}\label{an-bn} {\rm ev}_{\Omega(G_n)}((1-\sigma_v^{-1}e_{v,p})(e_{v,p}+(1-e_{v,p})b_n))={\rm ev}_{\Omega(G_n)}(a_nNI_{v,p}+ (1-\sigma_v^{-1}e_{v,p})b_n)\end{equation}
in $\frac{1}{|I_{v,p}|}R_{\Omega(G_n)}$, for all $n\gg 0$.\\

Next, we look at the above data in the Iwasawa tower $H_\infty/H$. Since the restriction maps
$$Sel_{\Sigma\cup J_v}^{\Sigma'\setminus J_v}(H_{n+1})_p^-\to Sel_{\Sigma\cup J_v}^{\Sigma'\setminus J_v}(H_{n})_p^-$$
are surjective (see lemma \ref{surjective-res-Selmer} ),  if $\pi_n^{n+1}:R_{n+1}^{\chi}\xrightarrow[]{}R_n^{\chi}$
denote the ring morphisms induced by Galois restriction, we have inclusions of $R_n^\chi$--ideals
$$\pi_n^{n+1}({\rm Fitt}_{R_{n+1}^{\chi}}(Sel_{\Sigma\cup J_v}^{\Sigma'\setminus J_v}(H_{n+1})_{R_{n+1}^{\chi}}^-)\subseteq {\rm Fitt}_{R_n^{\chi}}(Sel_{\Sigma\cup J_v}^{\Sigma'\setminus J_v}(H_n)_{R_n^{\chi}}^-), \qquad\text{ for all }n\gg0.$$
Therefore, for all $n\gg0$,  there exists $x_n\in R_n^{\chi}$,  such that
$$\Theta_{\Sigma\cup J_v}^{\Sigma'\setminus J_v
}(H_n/F)(e_{v,p}+(1-e_{v,p})\pi_n^{n+1}(b_{n+1}))=x_n\cdot \Theta_{\Sigma\cup J_v}^{\Sigma'\setminus  J_v}(H_n/F)(e_{v,p}+(1-e_{v,p})b_n).$$
Now, if one evaluates the equality above against a character $\psi\in\Omega(G_n)$, as in \eqref{Psi-Omega}, one sees that $\psi(x_n)\in\Bbb Z_p[\psi]^\times$. Consequently, $x_n\in (R_n^\chi)^\times$. Consequently, for any lift $x_{n+1}\in (R_{n+1}^\chi)^\times$ of $x_n$, one has the following equalities in $R_n^\chi$:
$$\pi_n^{n+1}(\Theta_{\Sigma\cup J_v}^{\Sigma'\setminus J_v}(H_{n+1}/F)(e_{v,p}+(1-e_{v,p})b_{n+1})x_{n+1}^{-1}) =\Theta_{\Sigma\cup J_v}^{\Sigma'\setminus J_v}(H_n/F)(e_{v,p}+(1-e_{v,p})b_n).$$
\\

Now, let us fix $n\gg 0$. The above considerations allow us to construct inductively elements $v_m\in \frac{1}{|I_{v,p}|}R_m^\chi$, for all $m\geq n$, such that
\begin{eqnarray}\label{the-vm}  v_n&=&(e_{v,p}+ (1-e_{v, p})b_n),\\
\nonumber (\Theta_{\Sigma\cup J_v}^{\Sigma'\setminus J_v}(H_m/F)\cdot v_m)_{R_m^\chi}&=&{\rm Fitt}_{R_m^{\chi}}(Sel_{\Sigma\cup J_v}^{\Sigma'\cup J_v}(H_m)_{R_m^{\chi}}^-),\quad \text{ for all }m\geq n,\\
\nonumber \pi^{m+1}_m(\Theta_{\Sigma\cup J_v}^{\Sigma'\setminus J_v}(H_{m+1}/F)\cdot v_{m+1})&=&(\Theta_{\Sigma\cup J_v}^{\Sigma'\setminus J_v}(H_m/F)\cdot v_m), \qquad \text{ for all }m\geq n.\end{eqnarray}
Consequently, from Proposition \ref{limit-of-Fitt}  and Lemma \ref{proj-lim-ideals}, we have
\begin{eqnarray}\label{Fitt-Sel-Jv-infty}{\rm Fitt}_{R_{\infty}^{\chi}}(Sel_{\Sigma\cup J_v}^{\Sigma'\setminus J_v}(H_{\infty})_{R_{\infty}^{\chi}}^-)&=&\varprojlim_{m\geq n}\text{ }( \Theta_{\Sigma\cup J_v}^{\Sigma'\setminus J_v}(H_m/F)v_m)\\
\nonumber &=&(\varprojlim_{m\geq n}\text{ } \Theta_{\Sigma\cup J_v}^{\Sigma'\setminus J_v}(H_m/F)v_m )\\
\nonumber &=& (\Theta_{\Sigma\cup J_v}^{\Sigma'\setminus J_v}(H_{\infty}/F)\cdot u_{\infty}),\end{eqnarray}
where $u_{\infty}\in\frac{1}{|I_{v,p}|}R_{\infty}^{\chi}$ can be chosen such that, if $\pi_m:R_\infty^\chi\to R_m^\chi$ is the $\Bbb Z_p[\chi]$--algebra morphism induced by Galois restriction and $u_m:=\pi_m(u_\infty)$, then
\begin{eqnarray}\label{ev-Omega-uinfty}
\pi_m( \Theta_{\Sigma\cup J_v}^{\Sigma'\setminus J_v}(H_{\infty}/F)\cdot u_{\infty} )&=& \Theta_{\Sigma\cup J_v}^{\Sigma'\setminus J_v}(H_m/F)v_m,\qquad\text{ for all }m\geq n,\\
\nonumber{\rm ev}_{\Omega(G_n)}(u_n)&=&{\rm ev}_{\Omega(G_n)}(e_{v,p}+(1-e_{v,p})b_n), \qquad \text{ in }\frac{1}{|I_{v,p}|}R_{\Omega(G_n)},\end{eqnarray}
where the last equality is a direct consequence of \eqref{the-vm} and the fact that ${\rm ev}_{\Omega(G_n)}(\Theta_{\Sigma\cup J_v}^{\Sigma'\setminus J_v}(H_n/F))$ is not a zero--divisor in $R_{\Omega(G_n)}$.\\

Now,  pass to the limit in \eqref{Fitt-Selmer-Jv-J} (keep in mind that $u_1=1$), apply Proposition \ref{limit-of-Fitt} and Lemma \ref{proj-lim-ideals} to obtain the following relations among  $R_\infty^\chi$--ideals
\begin{eqnarray}\label{Fitt-Sel-Jv-infty-inclusion} Fitt_{R_\infty^\chi}( Sel_{\Sigma\cup J_v}^{\Sigma'\setminus  J_v}(H_{\infty})_{R_{\infty}^{\chi}}^-) &\subseteq&  Fitt_{R_{\infty}^{\chi}}(Sel_{\Sigma\cup J}^{\Sigma'\setminus J_v}(H_{\infty})_{R_{\infty}^{\chi}}^-)=\\
\nonumber &=&(\Theta_{\Sigma\cup J}^{\Sigma'\setminus J_v}(H_{\infty}/F)(NI_{v,p},1-\sigma_v^{-1} e_{v,p}))_{R_\infty^\chi},\end{eqnarray}
where $\sigma_v$ is a choice of Frobenius of $v$ in $\mathcal G$.

Equalities \eqref{Fitt-Sel-Jv-infty} combined with \eqref{Fitt-Sel-Jv-infty-inclusion} and the fact that $\Theta_{\Sigma\cup J}^{\Sigma'\setminus J_v}(H_{\infty}/F)$ is a nonzero divisor in $R_\infty^\chi$ (see Proposition \ref{Theta-non-zero-divisor}) imply that there exist $\alpha, \beta\in R_\infty^\chi$, such that
\begin{equation}\label{alpha-beta} (1-\sigma_v^{-1} e_{v,p})\cdot u_{\infty}=\alpha NI_{v,p}+\beta (1-\sigma_v^{-1} e_{v,p}),\end{equation}
viewed as an equality inside $\frac{1}{|I_{v,p}|}R_\infty^\chi$. Now, we let $\alpha_n:=\pi_n(\alpha)$, $\beta_n:=\pi_n(\beta)$. If we apply ${\rm ev}_{\Omega(G_n)}\circ\pi_n$ to the equality above and combine it with \eqref{an-bn}, we obtain
\begin{eqnarray}\label{alphan-betan} {\rm ev}_{\Omega(G_n)}((1-\sigma_v^{-1} e_{v,p})\cdot u_{n})&=&{\rm ev}_{\Omega(G_n)}(a_n NI_{v,p}+b_n (1-\sigma_v^{-1} e_{v,p}))\\
\nonumber &=&{\rm ev}_{\Omega(G_n)}(\alpha_n NI_{v,p}+\beta_n (1-\sigma_v^{-1} e_{v,p})),\end{eqnarray}
viewed as equalities inside $\frac{1}{|I_{v,p}|}R_{\Omega(G_n)}$.\\

Assume that the isomorphism $\mathcal G\simeq G'\times G_p\times \Gamma$ maps $\sigma_v\mapsto (\sigma_v',\,\sigma_{v,p},\, \gamma_v)$. Then $\gamma_v$ is the Frobenius automorphism associated to $v$ in $G(F_\infty/F)$ (recall that $v\not\in S_p$, therefore $v$ is unramified in $F_\infty/F$.) Consequently, the closed subgroup $\overline{\langle\gamma_v^{-1}\rangle}$ generated by $\gamma_v^{-1}$ in $\Gamma$ has finite index and it is therefore equal to $\Gamma^{p^k}$, for some $k\geq 0$.
Consequently, we can pick a topological generator $\gamma$ of $\Gamma$, such that $\gamma_v^{-1}=\gamma^{p^k}$. Once $\gamma$ is chosen this way, we have an isomorphism of compact $\Bbb Z_p[\chi][G_p]$--algebras
$$R_\infty^\chi=\Bbb Z_p[\chi][G_p][[\Gamma]]\simeq \Bbb Z_p[\chi][G_p][[T]], \qquad \gamma\mapsto (T+1),$$
which sends the image $(\chi(\sigma_v')\sigma_{v,p}\cdot\gamma_v)^{-1}$ of $\sigma_v^{-1}$ in $R_\infty^\chi$ to $\chi(\sigma_v')^{-1}\sigma_{v,p}^{-1}\cdot(1+T)^{p^k}$. Consequently, in the following calculations which happen inside $R_\infty^\chi$, we can and will  identify
$$R_\infty^\chi=\Lambda[G_p], \qquad \Lambda=\Bbb Z_p[\chi][[T]], \quad \sigma_v^{-1}=c\cdot h\cdot (1+T)^{p^k},$$
where $h:=\sigma_{v,p}^{-1}\in G_p$ and $c:=\chi'(\sigma_v')^{-1}\in\Bbb Z_p[\chi]^\times$. \\

Now, equality \eqref{alpha-beta} projected via Galois restriction onto  $\frac{1}{|I_{v,p}|}\Lambda[G_p/I_{v,p}]$ reads
\begin{equation}\label{equation u_infty-beta}(\overline{u_{\infty}}-\overline\beta)(1-c\cdot \overline h\cdot (1+T)^{p^k})=\overline \alpha|I_{v,p}|, \qquad \text{ in } \frac{1}{|I_{v,p}|}\Lambda[G_p/I_{v,p}].\end{equation}
 Next, we solve this equation for $(\overline{u_{\infty}}-\overline\beta)$ in $\frac{1}{|I_{v,p}|}\Lambda[G_p/I_{v,p}]$. To that end, assume that
$$\overline{u_{\infty}}-\overline \beta=\sum_{g\in G_p/I_{v,p}}a_gg, \qquad \overline\alpha=\sum_{g\in G_p/I_{v,p}}b_gg,$$
where $a_g\in\frac{1}{|I_{v,p}|}\Lambda$ and $b_g\in\Lambda$, for all $g\in G_p/I_{v,p}$. Now, fix $g\in G_p/I_{v,p}$ and compare the coefficients of $g\overline h^{-i}$ on both sides of equality\eqref{equation u_infty-beta} to obtain
$$a_{g\overline h^{-i}}-a_{g\overline h^{-i-1}}c(1+T)^{p^k}=|I_{v,p}|b_{g\overline h^{-i}}, \qquad\text{ for all }i=0,1,2,\dots, {\rm ord}(\overline h)-1,$$
where ${\rm ord}(\overline h)$ is the order of $\overline h$ in $G_p/I_{v, p}$.
Now, by multiplying the equalities above with $c^i(1+T)^{p^k\cdot i}$ and taking the sum for all $i$ as above, we obtain an equality in $\Lambda$
 $$(|I_{v,p}|a_g)(1-c^{ord(\overline h)}(1+T)^{p^k ord(\overline h)})=|I_{v,p}|^2\sum_{1=0}^{ord(\overline h)-1}b_{g\overline h^{-i}}c^i(1+T)^{p^k\cdot i}.$$
 However, $\Lambda$ is a unique factorization domain, and $p\in\Lambda$ is prime. Therefore, by comparing the power of $p$ dividing both sides in the above equality, one can conclude right away
 that $a_g\in |I_{v,p}|\Lambda$. Therefore, we have
 $$\overline{u_{\infty}}-\overline \beta=|I_{v,p}|\overline{\alpha'}, \qquad \text{ with }\alpha'\in \Lambda[G_p].$$
 Therefore, if $\tau $ is a generator of the cyclic group $I_{v,p}$, we have
 $$u_{\infty}-\beta=|I_{v,p}|\cdot\alpha' +(1-\tau)\alpha'', \qquad \text { with }\alpha'\in \Lambda[G_p], \quad \alpha''\in\frac{1}{|I_{v,p}|}\Lambda[G_p].$$
After reidentifying $\Lambda[G_p]$ with $R_\infty^\chi$ and applying $\pi_n$ to the above equality, we obtain the following equality in $\frac{1}{|I_{v,p}|}R_n^\chi$.
\begin{equation}\label{un-minus-betan}u_n-\beta_n=|I_{v,p}|\cdot\alpha_n' +(1-\tau)\alpha_n'', \qquad \text { with }\alpha'_n\in R_n^\chi, \quad \alpha_n''\in\frac{1}{|I_{v,p}|}R_n^\chi.\end{equation}
\\

Now, the second equality in \eqref{alphan-betan} gives us the following
\begin{equation}\label{ev(1-bn)} {\rm ev}_{\Omega(G_n)}(e_{v,p}(1-b_n))={\rm ev}_{\Omega(G_n)}([(\alpha_n-a_n)|I_{v,p}|+\beta_n-b_n]\sigma+(1-\beta_n)e_{v,p}).\end{equation}
Now, observe that by \eqref{ev-Omega-uinfty}, we have
$${\rm ev}_{\Omega(G_n)}(u_n e_{v,p})={\rm ev}_{\Omega(G_n)}((e_{v,p}(1-b_n)+b_n)e_{v,p})={\rm ev}_{\Omega(G_n)}(e_{v,p}).$$
Therefore, since  $(1-\tau)e_{v,p}=0$, by \eqref{un-minus-betan} we obtain
\begin{eqnarray} \nonumber {\rm ev}_{\Omega(G_n)}((1-\beta_n)e_{v,p})={\rm ev}_{\Omega(G_n)}((u_{n}-\beta_n)e_{v,p})&=&{\rm ev}_{\Omega(G_n)}((|I_{v,p}|{\alpha'_n}+(1-\tau)\alpha''_n)e_{v,p})\\
\nonumber &=&{\rm ev}_{\Omega(G_n)}(|I_{v,p}|{\alpha'_n})\in R_{\Omega(G_n)}.\end{eqnarray}
When combined with \eqref{ev(1-bn)}, this leads to
$${\rm ev}_{\Omega(G_n)}(b_n+ e_{v,p}(1-b_n))\in R_{\Omega(G_n)}.$$
When combined with \eqref{vn-psi-eval}, this leads to
$${\rm ev}_{\Omega(G_n)}(b_n+ e_{v,p}(1-b_n))\in R_{\Omega(G_n)}^\times.$$
Now, the map ${\rm ev}_{\Omega(G_n)}: R_n^\chi \to R_{\Omega(G_n)}$ clearly maps the unit group $(R_n^\chi)^\times$ surjectively onto the unit group $R_{\Omega(G_n)}^\times$. Therefore, there exists $\xi_n\in (R_n^\chi)^\times$ and $\rho_n\in\ker({\rm ev}_{\Omega(G_n)})$, such that
$$ b_n+ e_{v,p}(1-b_n)=\xi_n+\rho_n, \qquad \text{ in }R_n^\chi.$$
Consequently, since $\Theta_{\Sigma\cup J_v}^{\Sigma'\setminus J_v}(H_n/F)\cdot\ker({\rm ev}_{\Omega(G_n)})=0$, we have equalities of $R_n^\chi$--ideals
\begin{eqnarray} \nonumber {\rm Fitt}_{R_n^\chi}(Sel_{\Sigma\cup J_v}^{\Sigma'\setminus J_v}(H_n)_{R_n^\chi})&=&(\Theta_{\Sigma\cup J_v}^{\Sigma'\setminus J_v}(H_n/F)\cdot(b_n+ e_{v,p}(1-b_n))\\
\nonumber &=&(\Theta_{\Sigma\cup J_v}^{\Sigma'\setminus J_v}(H_n/F)\cdot\xi_n)= (\Theta_{\Sigma\cup J_v}^{\Sigma'\setminus J_v}(H_n/F))_{R_n^\chi},\end{eqnarray}
which completes the inductive step and hence the proof of the theorem.
\end{proof}
A consequence of the above theorem is the following.
\begin{theorem}\label{Fitt-at-n}
Let $F_{\infty}$ be the cyclotomic $\mathbb{Z}_p$-extension of $F$. Assume that $F_\infty\cap H=F$. Let $S$ and $T$ be two finite, disjoint sets of places in $F$ satisfying hypotheses ${\rm Hyp}(H_\infty/F)_p$  in Lemma \ref{Deligne-Ribet-Kurihara}. Then,   we have equalities of $\Bbb Z_p[G_n]^-$--ideals
$$Fitt_{\mathbb{Z}_p[G_n]^-}(Sel_S^T(H_n)_p^-)=(\Theta_S^T(H_n/F)), \qquad\text{ for all }n\gg0.$$
\end{theorem}
\begin{proof}
 Apply Theorem \ref{full} to the set of data $(H/F, S_0, T_0, J)$, where
$$S_0:=S\setminus S_{ram}(H_\infty/F), \quad T_0:=T\setminus S_{ram}(H_\infty/F), \quad J:=(S_{ram}(H_\infty/F)\setminus S_{p\infty}(F))\setminus T, $$
and observe that,  for $\Sigma$ and $\Sigma'$ as defined in \eqref{Sigma-Sigma'}, we have set equalities
$$\Sigma\cup J=S, \qquad \Sigma'\setminus J= T,$$
which concludes the proof.
\end{proof}

Further, we obtain the following strengthening of Theorem 3.3 in \cite{Dasgupta-Kakde} (the keystone result in loc.cit.), in the case where $S_p(F)\subseteq S$. This will imply our main result. The theorem below will be considerably strengthened again in \S7.4. (See Theorem \ref{BKS-Conj} in \S7.4.)
\begin{theorem}\label{keystone-improved}
Let $H/F$ be a finite, abelian extension of number fields, with $G=Gal(H/F)$, $F$ totally real and $H$ a CM field.  Let $p>2$ be prime and $S$ and  $T$ be disjoint,  finite, nonempty sets of places in $F$, satisfying hypotheses ${\rm Hyp}(H_\infty/F)_p$. Then, we have
$$Fitt_{\mathbb{Z}_p[G]^-}(Sel_S^T(H)_p^-)=(\Theta_S^T(H/F))$$
\end{theorem}
\begin{proof} Throughout, if $k$ is a number field, $k_\infty$ denotes its cyclotomic $\Bbb Z_p$--extension. 
As noted before, there exists a CM field $H'$ such that $H'\cap F_{\infty}=F$ and $H'_\infty=H_\infty$. Now, let $H_n'$ denote the $n$-th layer of $H_{\infty}/H'$. Since, $\bigcup_{n=1}^{\infty}H_n'=H_{\infty}$, there exists $n\gg 0$ such that $H\subseteq H_n'$. Let us fix an $n$ with this property.

As in the proof of Theorem \ref{full}, we consider the set of characters of $p$--adic characters of $G$
$$\Omega:=\Omega(G):=\{\psi\in\widehat G\,\vert\, \psi(\Theta_S^T(H/F))\ne 0\}.$$  
Then, for all $n\geq 0$, we have canonical surjective ring morhpism 
$$\pi_n: \Bbb Z_p[G'_n]^-\overset{res_n}\longrightarrow\Bbb Z_p[G]^-\overset{{\rm ev}_\Omega}\longrightarrow R_{\Omega},$$
such that $\pi_n(\Theta_S^T(H_n'/F))=\pi_0(\Theta_S^T(H/F))$ is not a zero divisor in $R:=R_\Omega$. We have the following surjective restriction morphism at the level of Selmer modules
\begin{equation}\label{Selmer-surjection}(Sel_S^T(H_n')_p^-)_R\xrightarrow[]{} (Sel_S^T(H)_p^-)_R\end{equation}
Therefore, by the properties of Fitting ideals and Theorem \ref{Fitt-at-n},  we have 
$$(\pi_n(\Theta_S^T(H_n'/F)))=Fitt_R((Sel_S^T(H_n')_p^-)_R)\subseteq Fitt_R((Sel_S^T(H)_p^-)_R)=(\pi_0(\theta)), \qquad\text{ for } n\gg 0.$$
Here,  $\theta$ is any generator of $Fitt_{\Bbb Z_p[G]^-}(Sel_S^T(H)_p^-)$, which we know to be principal since $Sel_S^T(H)_p^-$ is quadratically presented over $\Bbb Z_p[G]^-$. (See Lemma \ref{quadratic-pres-lemma}(2).) Now, since we know that in 
$\pi_n(\Theta_S^T(H_n'/F))$ is a nonzero divisor in $R$, so is $\pi_0(\theta)$. Consequently, by Corollary \ref{covariance}, morphism \eqref{Selmer-surjection} is an isomorphism. Therefore, we have an equality of Fitting ideals
$$Fitt_R((Sel_S^T(H)_p^-)_R)=Fitt_R((Sel_S^T(H_n')_p^-)_R)=(\pi_0(\Theta_S^T(H/F)))$$
However, if $R'=R_{\Omega^c}$ and $\pi'_0:\Bbb Z_p[G]^-\to R'$ is the corresponding surjection, by Lemma 3.2 in \cite{Dasgupta-Kakde}, we have
$$Fitt_{R'}((Sel_S^T(H)_p^-)_{R'})=(0)=(\pi_0'(\Theta_S^T(H/F)))$$
The he last two equalities and the injectivity of the map $\pi_0\oplus\pi_0':\Bbb Z_p[G]^-\to R\oplus R'$ give the desired result.
\end{proof}

Now, we are ready to state and prove the main  theorem of this paper.
\begin{theorem}[EMC] \label{EMC} Let $H/F$ be an abelian, CM extension of a totally real number field. Let $H_\infty$ be the cyclotomic $\Bbb Z_p$--extension of $H$, for a prime $p>2$.
Let $\mathcal G:=Gal(H_{\infty}/F)$. Suppose that $S$ and $T$ are finite, disjoint sets of places of $F$, satisfying hypotheses ${\rm Hyp}(H_\infty/F)_p$ in Lemma \ref{Deligne-Ribet-Kurihara}. Then, the following hold.
\begin{enumerate}
\item We have an equality of $\Bbb Z_p[[\mathcal G]]^-$--ideals
$${\rm Fitt}_{\mathbb{Z}_p[[\mathcal{G}]]^-}(Sel_S^T(H_{\infty})_p^-)=(\Theta_S^T(H_\infty/F)).$$
\item The $\Bbb Z_p[[\mathcal G]]^-$--module $Sel_S^T(H_{\infty})_p^-$ sits in a short exact sequence
$$0\longrightarrow(\Bbb Z_p[[\mathcal G]]^-)^k\longrightarrow (\Bbb Z_p[[\mathcal G]]^-)^k\longrightarrow Sel_S^T(H_{\infty})_p^-\longrightarrow 0,$$
for some $k>0$. In particular, ${\rm pd}_{\Bbb Z_p[[\mathcal G]]^-} Sel_S^T(H_{\infty})_p^- =1.$
\end{enumerate}
\end{theorem}
\begin{proof}  (1) Theorem \ref{keystone-improved} applied to $H_n/F$, for all $n\geq 0$, and Proposition \ref{limit-of-Fitt} imply that we have an equality of $\Bbb Z_p[[\mathcal G]]^-$--ideals
 $${\rm Fitt}_{\mathbb{Z}_p[[\mathcal{G}]]^-}(Sel_S^T(H_{\infty})_p^-)=\varprojlim_n \, (\Bbb Z_p[G_n]^-\cdot \Theta_S^T(H_n/F)).$$
When combined with Lemma \ref{proj-lim-ideals}, this gives a further equality of $\Bbb Z_p[[\mathcal G]]^-$--ideals
 $${\rm Fitt}_{\mathbb{Z}_p[[\mathcal{G}]]^-}(Sel_S^T(H_{\infty})_p^-)=\Bbb Z_p[[\mathcal G]]^-\cdot \Theta_S^T(H_\infty/F),$$
 which concludes the proof of part (1).
 \\

(2) Since the ring $R:=\Bbb Z_p[[\mathcal G]]^-$ is semilocal (see Remark \ref{semi-local-remark}) and $\Theta_S^T(H_\infty/F)$ is not a zero divisor in $R$  (see Proposition \ref{Theta-non-zero-divisor}), we can apply Proposition \ref{Cornacchia-Greither} to the finitely presented, torsion  $R$--module $Sel_S^T(H_{\infty})_p^-$ to conclude that, in light of part (1),  part (2) holds as well.
\end{proof}

We conclude this section with a Corollary of the above theorem which gives information on a more classical Iwasawa module, obtained by taking projective limits of ray class--groups in the Iwasawa tower $H_\infty/H$.
\begin{corollary}
With notations as above, assume that $$S:=S_{p\infty}(F), \qquad T=(S_{ram}(H_\infty/F)\setminus S_{p\infty}(F))\cup T_0,$$ for a finite, non-emtpy set of primes $T_0$ in $F$ such that $T_0\cap S_{ram}(H_\infty/F)=\emptyset.$
Then
\begin{enumerate}
\item We have an equality of $\Bbb Z_p[[\mathcal G]]^-$--ideals
$${\rm Fitt}_{\mathbb{Z}_p[[\mathcal{G}]]^-}(A^{T,\vee, -})=(\Theta_S^T(H_\infty/F)).$$
\item The $\Bbb Z_p[[\mathcal G]]^-$--module $A^{T,\vee, -}$ sits in a short exact sequence
$$0\longrightarrow(\Bbb Z_p[[\mathcal G]]^-)^k\longrightarrow (\Bbb Z_p[[\mathcal G]]^-)^k\longrightarrow A^{T,\vee, -}\longrightarrow 0,$$
for some $k>0$. In particular, ${\rm pd}_{\Bbb Z_p[[\mathcal G]]^-} A^{T,\vee, -} =1.$
\end{enumerate}
\end{corollary}
\begin{proof}
Under our hypotheses, Proposition \ref{SES-infinity} gives an isomorphism of $\Bbb Z_p[[\mathcal G]]^-$--modules
$$Sel_{S}^{T}(H_{\infty})_p^-\cong A^{T,\vee,-}.$$
Consequently, the statement above follows from the previous Theorem.
\end{proof}


\section{Applications}\label{Applications}

In this section, we give several applications of our main result, Theorem \ref{EMC}. The notations are as above. In addition, we let $\mu:=\mu_{H,p}$ denote the classical Iwasawa $\mu$--invariant associated to the $\Bbb Z_p$--extension $H_\infty/H$.
In \S7.1 we establish a perfect, $\Bbb Z_p[[\mathcal G]]$--bilinear duality pairing between the Selmer module $Sel_S^T(H_\infty)^-$ and a certain module $\mathcal M_S^T(H_\infty)^-$, which was essentially introduced in \cite{Greither-Popescu}. As a consequence, in \S7.2 we derive (as a corollary to Theorem \ref{EMC}) the main result of \cite{Greither-Popescu}, namely the Equivariant Main Conjecture for the Tate module $T_p( \mathcal M_S^T(H_\infty)^-)$, proved in loc.cit. under the hypothesis  that $\mu=0$ (a classical conjecture of Iwasawa.) Further, in \S7.3 we give an unconditional proof of a refinement of the Coates--Sinnott Conjecture, proved in \cite{Greither-Popescu} under the hypothesis that $\mu=0$. Finally, in \S7.4 we prove a considerable strengthening of our Theorem \ref{keystone-improved}, by replacing hypotheses ${\rm Hyp}(H_\infty/F)_p$ with the weaker ${\rm Hyp}(H/F)_p$ and, as a consequence, prove a conjecture of Burns--Kurihara--Sano.
\\

Before we begin, recall that Iwasawa's $\mu=0$ conjecture (currently known to hold if $H/\Bbb Q$ is abelian) is equivalent to the $p$--divisibility of the group  $A_\infty^{T, -}$, for any finite set $T$ (including the empty set) of non--archimedean primes in $F$.  In general, one has group isomorphisms
$$A_\infty^{T}\simeq (\Bbb Q_p/\Bbb Z_p)^{\lambda_T}\bigoplus (A_\infty^{T})_{nd}, \qquad A_\infty^{T,\pm}\simeq (\Bbb Q_p/\Bbb Z_p)^{\lambda_T^\pm}\bigoplus (A_\infty^{T,\pm})_{nd},$$
where $\lambda_T, \lambda_T^\pm\in\Bbb Z_{\geq 0}$ and  $(A_\infty^{T})_{nd}$ and $(A_\infty^{T, \pm})_{nd}$  are torsion $p$--groups of exponents $p^{\mu}$ and $p^{\mu^\pm}$, respectively, where $\mu$ and $\mu^{\pm}$ are independent of $T$.  It turns out that $\mu=0$ iff $\mu^-=0$, and that happens iff $A_\infty^{T,-}$ is divisible, as the isomorphisms above show. (See \S2 in \cite{Greither-Popescu} for more details.)

 From now on, for notational simplicity we let $\cH:=H_\infty$ and $A^T:=A_\infty^T$, for all finite sets $T$ of non-archimedean primes in $F$.

\subsection{Duality pairings} Let $S$ and $T$ be two finite, disjoint  sets of primes in $F$, with $S_{p\infty}(F)\subseteq S$. The set $T$ can be empty. For all sets $X$ of non--archimedian primes in $F$, such that
$S_p\cap X=\emptyset$ (so $\cH/H$ is unramified at $X$), we have injective transition maps
\begin{equation}\label{maps-between-divisors} Y_X(H_n)\longrightarrow  Y_X(H_{n+1}), \qquad v\mapsto \sum_{w|v} w,\end{equation}
 which are compatible in the obvious way with the divisor maps $div_X^{H_i}:H_i^\times\to Y_X(H_i)$, for $i=n, n+1$ and the inclusions $H_n^\times\subseteq  H_{n+1}^\times$. Consequently, we obtain a divisor map
 $$div_X^{\cH}:=\varinjlim_n\, div_X^{H_n}: \,\, \cH^\times\to Y_X(\cH):=\varinjlim_n\, Y_X(H_n).$$


Following \cite{Greither-Popescu}, for all $m\geq 1$, we define the following subgroups of $\cH_T^\times$:
$$\cH_{S,T}^{(p^m)}=\{f\in \cH_T^{\times}\, \vert\, div_{\overline{S_{p\infty}}}^{\cH}(f)=p^m\cdot D+y,\, \text{ with } y\in Y_{S\setminus S_{p\infty}}(\cH), \, D\in Y_{\overline{S}}(\cH)\}.$$

\begin{definition} For all $m>0$,  we define the $\Bbb Z_p[[\mathcal G]]$--linear morphisms
$$\iota_m: {\cH_{S,T}^{(p^m)}}/{(\cH_T^\times)^{p^m}} \longrightarrow { \cH_{S,T}^{(p^{m+1})}}/{(\cH_T^\times)^{p^{m+1}}} , \qquad \iota_m(\widehat f)=\widehat{f^p},\text { for all $f\in \cH_{S,T}^{(p^m)}$},$$
and let  $\mathcal M_S^T(\cH):=\varinjlim_m {\cH_{S,T}^{(p^m)}}/{(\cH_T^\times)^{p^m}}$, where the injective limit is taken with respect to the transition maps $\iota_m$ defined above.
\end{definition}

One can define in a similar manner the $\Bbb Z_p[[\mathcal G]]$--modules $\mathcal M_S^T(H_n)$, for all $n\geq 1$ and one can check without dificulty that the inclucions $H_{n}^\times\subseteq H_{n+1}^\times\subseteq \cH^\times$ induce injective morphisms
$$\mathcal M_S^T(H_n)\overset{j_n}\longrightarrow\mathcal M_S^T(H_{n+1})\longrightarrow \mathcal M_S^T(\cH).$$
Further, it is easily checked that we have $\Bbb Z_p[[\mathcal G]]$--module isomorphisms
\begin{equation}\label{inj-limit-motives}
\mathcal M_S^T(\cH)\simeq \varinjlim_n\mathcal M_S^T(H_n), \qquad \mathcal M_S^T(\cH)^\vee\simeq \varprojlim_n\mathcal M_S^T(H_n)^\vee,
\end{equation}
where the injective limit is taken with respect to the injective morphisms $j_n$, while the projective limit is taken with respect to their surjective duals $j_n^\vee$.

\begin{remark}\label{maps-pi_m} It is easily checked that the maps $\iota_m$ are injective, which allows us to view $\cH_{S,T}^{(p^m)}/(\cH_T^\times)^{p^m}$  as a submodule of $\mathcal M_S^T(\mathcal H)$.  In fact, $\iota_m$ induces an isomorphism
$$\cH_{S,T}^{(p^m)}/(\cH_T^\times)^{p^m}\simeq \mathcal M_S^T(\cH)[p^m],$$
for all $m>0$, as one can check. (The group of roots of unity in the cyclotomic $\Bbb Z_p$--extension $\cH$ is always $p$--divisible!) Further, it is obvious that these isomorphism fit into commutative diagrams
\[\begin{tikzcd}
\cH_{S,T}^{(p^{m+1})}/(\cH_T^\times)^{p^{m+1}}\arrow{r}{\sim}\arrow{d}{\pi_m} & \mathcal M_S^T(\cH)[p^{m+1}]\arrow{d}{\times p}\\
\cH_{S,T}^{(p^m)}/(\cH_T^\times)^{p^m}\arrow{r}{\sim} & \mathcal M_S^T(\cH)[p^m],
\end{tikzcd}\]
where the left--vertical map is induced by the inclusion $\cH_{S,T}^{(p^{m+1})}\subseteq \cH_{S,T}^{(p^{m})}$ and the right--vertical map is multiplication by $p$. This shows that there is an isomorphism of $\Bbb Z_p[[\mathcal G]]$--modules
$$ \varprojlim_m \cH_{S,T}^{(p^m)}/(\cH_T^\times)^{p^m}\simeq T_p(\mathcal M_S^T(\cH)),$$
where the inverse limit is taken with respect to the maps $\pi_m$ and  $T_p(\mathcal M_S^T(\cH))$ denotes the $p$--adic Tate module of the $p$--group $\mathcal M_S^T(\cH)$.
\end{remark}

The description of $\mathcal M_S^T(\cH)[p^m]$ in Remark \ref{maps-pi_m} gives rise to canonical exact sequences of $\Bbb Z_p[[\mathcal G]]^-$--modules
\begin{equation}\label{p^m-sequence} 0\longrightarrow  A^T[p^m]^-\overset{\tau_m}\longrightarrow {\mathcal M}_S^T(\cH)[p^m]^-\overset{d_m}\longrightarrow Y_{S\setminus S_p}(\cH)_p^-/p^m,\end{equation}
where $\tau_m$ takes a divisor class $\widehat D$, such that $p^mD={\rm div}(f)$, for some $f\in\cH_T^\times$, to the class $\widehat f$ in
${\mathcal M}_S^T(\cH)[p^m]$, and $d_m(\widehat g)= (div_{S\setminus S_{p\infty}}^\cH(g)\mod p^mY_{S\setminus S_{p}}(\cH)^-_p)$, for all $g\in \cH_{S,T}^{(p^m)}$. We leave it to the reader to check that the maps $\tau_m$ and $d_m$ are well--defined and that the sequence above is exact. Further, we have commutative diagrams of $\Bbb Z_p[[\mathcal G]]^-$--modules
 \[\begin{tikzcd}
 0\arrow{r} &A^T[p^{m+1}]^-\arrow{r}{\tau_{m+1}} & {\mathcal M}_S^T(\cH)[p^{m+1}]^-\arrow{r}{d_{m+1}} &Y_{S\setminus S_p}(\cH)_p^-/p^{m+1}\\
0\arrow{r} &A^T[p^m]^-\arrow{r}{\tau_m}\arrow{u} & {\mathcal M}_S^T(\cH)[p^m]^-\arrow{r}{d_m}\arrow{u}{\iota_m} &Y_{S\setminus S_p}(\cH)_p^-/p^m\arrow{u}{\times p}
\end{tikzcd}\]
where the left--vertical map inclusion and the right--vertical maps is given by multiplication by $p$. In particular, all the vertical arrows are injective.

\begin{lemma}
If we take an injective limit with respect to $m$ in the commutative diagrams above, we obtain an exact sequence of $\Bbb Z_p[[\mathcal G]]^-$--modules
\begin{equation}\label{p^infty-sequence} 0\longrightarrow A^{T, -}\overset{\tau}\longrightarrow {\mathcal M}_S^T(\cH)^-\overset{d}\longrightarrow Y_{S\setminus S_p}(\cH)_p^-\otimes\Bbb Q_p/\Bbb Z_p\longrightarrow 0,\end{equation}
where
$\tau:=(\tau_m)_m$, and $d:=(d_m)_m$.
\end{lemma}
\begin{proof} It suffices to show that the map $d$ is surjective (although the maps $d_m$ may not be surjective, in general.) This is equivalent to showing that for all $m$ and all $y\in Y_{S\setminus S_{p}}(\cH)^-_p$, there exists an
$n\geq 0$ such that $$p^n y\mod  p^{m+n}Y_{S\setminus S_{p}}(\cH)^-_p \, \in\, {\rm Im}(d_{m+n}).$$
Fix $m$ and $y$ as above. Since $p^\mu A^{T,-}$ is a $p$--divisible group, there exists a divisor $D$ in $Y_{\overline {S_{p\infty}}}(\cH)_p^-$, such that we have the following equality of divisor classes in $A^{T, -}.$
$$p^\mu\widehat y=p^{m+\mu}\widehat D.$$
Equivalently, there is an element $g\in \cH_T^\times$, such that
$$div^{\cH}(g)=p^\mu y-p^{m+\mu} D.$$
Note that this implies that $g\in \cH_{S,T}^{(p^{m+\mu})}$ and that $$d_{m+\mu}(\widehat g)=p^\mu y \mod p^{m+n}Y_{S\setminus S_{p}}(\cH)^-,$$ which means that $(p^\mu y\mod  p^{m+\mu}Y_{S\setminus S_{p}}(\cH)^-_p) \in {\rm Im}(d_{m+\mu})$, concluding the proof.
\end{proof}
\begin{remark} One can easily show by following the steps in the proof of the above Lemma that if $\mu=0$, then the maps $d_m$ are surjective, for all $m$.
\end{remark}
After taking Pontrjagin duals in \eqref{p^infty-sequence} and noting that since $Y_{S\setminus S_p}(\mathcal H)_p$ is $\Bbb Z_p$--free we have a natural isomorphism
$(Y_{S\setminus S_p}(\cH)_p^-\otimes\Bbb Q_p/\Bbb Z_p)^\vee\simeq Y_{S\setminus S_p}(\cH)_p^{\ast, -}$, we obtain an exact sequence of $\Bbb Z_p[[\mathcal G]]^-$--modules
\begin{equation}\label{dual-p^infty-sequence}
0\longrightarrow Y_{S\setminus S_p}(\cH)_p^{\ast, -}\overset{d^\vee}\longrightarrow {\mathcal M}_S^T(\cH)^{\vee, -}\overset{\tau^\vee}\longrightarrow A^{T, \vee, -}\longrightarrow 0.
\end{equation}
When combined with Proposition \ref{SES-infinity}, the exact sequence above suggests
that there should be a close link between the Selmer module $Sel_S^T(\cH)_p^-$ and
the module $\mathcal M_S^T(\cH)^{\vee, -}$.

Indeed, one can check that for all $n\geq 1$ we have well--defined,  $\Bbb Z_p[[\mathcal G]]^-$--bilinear pairings
$$\langle\cdot\, ,\, \cdot\rangle_n \, :\,  Sel_S^T(H_n)_p^-\times \mathcal M_S^T(H_n)^-\longrightarrow \Bbb Q_p/\Bbb Z_p,  \qquad \langle\widehat\phi, \widehat{f_m}\rangle:= \frac{\phi(f_m)}{p^m},$$
for all $\phi\in (H_{n,T}^\times)_p^\ast$ and all  $f_m\in H_{n, S, T}^{(p^m)}.$ The above bilinear pairings induce natural $\Bbb Z_p[[\mathcal G]]^-$--linear morphisms
$$\Phi_{n,S}^T:  Sel_S^T(H_n)_p^-\longrightarrow \mathcal M_S^T(H_n)^{\vee, -}, \qquad \Phi_{n, S}^T(\widehat\phi)(\widehat{f_m}):=\frac{\phi(f_m)}{p^m},$$
for all $n\geq 1$, all $\phi\in (H_n^\times)_p^\ast$ and all  $f_m\in H_{n, S, T}^{(p^m)}.$ For all $n\geq 1$, we have obvious commutative diagrams
\[\begin{tikzcd}
Sel_S^T(H_{n+1})_p^-\arrow{r}{\Phi_{n+1,S}^T}\ar{d}{res_n} &  \mathcal M_S^T(H_{n+1})^{\vee, -}\ar{d}{j_n^\vee}\\
Sel_S^T(H_{n})_p^-\arrow{r}{\Phi_{n,S}^T} &  \mathcal M_S^T(H_{n})^{\vee, -}
\end{tikzcd}\]
Now, if we take the projective limit with respect to $n$ in the diagrams above and take into account isomorphisms \eqref{inj-limit-motives}, we obtain a $\Bbb Z_p[[\mathcal G]]^-$--linear morphism
$$\Phi_S^T:\, Sel_S^T(\cH)_p^-\longrightarrow \mathcal M_S^T(\cH)_p^{\vee, -}, \qquad \Phi_S^T:=(\Phi_{n,S}^T)_n,$$
whose properties are listed in the following propositon.
\medskip

\begin{proposition}\label{link-sel-mst}
The following hold, for all finite, disjoint sets $S$ and $T$, with $S_{p\infty}(F)\subseteq S$.
\begin{enumerate}
\item
We have a commutative diagram of $\Bbb Z_p[[\mathcal G]]^-$--modules with exact rows
\[ \begin{tikzcd}
0\arrow{r} & Y_{S\setminus S_p}(\cH)_p^{\ast, -}\arrow{r}\arrow{d}{{\rm id}}[swap]{=} & Sel_S^T(\cH)_p^-\arrow{r}\arrow{d}{\Phi_S^T} & A^{T,\vee,-}\arrow{r}\arrow{d}{{\rm id}}[swap]{=} & 0\\
0\arrow{r}  &Y_{S\setminus S_p}(\cH)_p^{\ast, -}\arrow{r}{d^\vee} &{\mathcal M}_S^T(\cH)^{\vee, -}\arrow{r} {\tau^\vee}& A^{T, \vee, -}\arrow{r} & 0,
\end{tikzcd}
\]
where the top row is the exact sequence in Proposition \ref{SES-infinity} and the bottom row is \eqref{dual-p^infty-sequence} above.
\item The map $\Phi_S^T$ induces a natural isomorphism of $\Bbb Z_p[[\mathcal G]]^-$--modules
$$
\Phi_S^T: Sel_S^T(\cH)_p^- \simeq  {\mathcal M}_S^T(\cH)^{\vee, -}.
$$
\item If $\mu:=\mu_{H,p}=0$, then we have a natural isomorhpsim of $\Bbb Z_p[[\mathcal G]]^-$--modules
$$ Sel_S^T(\cH)_p^- \simeq T_p(\mathcal M_S^T(\cH)^-)^\ast.$$
\end{enumerate}
\end{proposition}
\begin{proof} Part (1) can be easily checked by using the explicit definitions of all the maps involved. We leave it to the reader, as an exercise. Part (2) is an immediate consequence of part (1) via the snake lemma.

Now, assume that $\mu=0$. Since $A^{T,-}$ is $p$--divisible, torsion, of finite corank, \eqref{p^infty-sequence} implies that so is  $\mathcal M_S^T(\cH)^{-}$ (i.e. it is isomorphic to a direct sum of finitely many copies of $\Bbb Q_p/\Bbb Z_p$.) Consequently, there is a natural isomorphism of $\Bbb Z_p[[\mathcal G]]^-$--modules
$$\mathcal M_S^T(\cH)^{-,\vee}\simeq T_p(\mathcal M_S^T(\cH)^-)^\ast,\qquad \phi\mapsto \overline\phi,\quad \overline\phi((x_m)_m)=(p^m\widetilde{\phi(x_m)}\mod p^m)_m,$$
for all $\phi\in\mathcal  M_S^T(\cH)^{-,\vee}$ and all $(x_m)_m\in \varprojlim_m\mathcal M_S^{T}(\cH)^-[p^m]=T_p(\mathcal M_S^T(\cH)^{-})$, where $\widetilde{\phi(x_m)}$ is a lift of $\phi(x_m)$ to $\frac{1}{p^m}\Bbb Z_p$. Now, part (3) follows from part (2).
\end{proof}
\begin{corollary}\label{true-Greither-Popescu}
If $S_{p\infty}(F)\subseteq S$, $S_{ram}(\cH/F)\subseteq S\cup T$ and $T\not\subseteq S_{ram}(\cH/F)$, then the following hold.
\begin{enumerate}
\item The $\Bbb Z_p[[\mathcal G]]^-$--module ${\mathcal M}_S^T(\cH)^{\vee, -}$ is quadratically presented and of projective dimension $1$.
\item ${\rm Fitt}_{\mathbb{Z}_p[[\mathcal{G}]]^-}({\mathcal M}_S^T(\cH)^{\vee, -})=(\Theta_S^T(H_\infty/F)).$
\end{enumerate}
\end{corollary}
\begin{proof}
This is a direct consequence of Theorem \ref{EMC} and part (2) of Proposition \ref{link-sel-mst}.
\end{proof}
\medskip

\subsection{The Equivariant Main Conjecture of Greither--Popescu}
In \cite{Greither-Popescu}, {\it working under the hypothesis that $\mu=0$,} and assuming that $S_{ram}(\cH/F)\subseteq S$ and $T\neq\emptyset$, the authors introduce the abstract $p$--adic $1$--motive
$$\mathcal M_{S,T}^{\cH, -}:=[Y_{S\setminus S_p}(\cH)_p^-\overset{\widehat{\rm div}}\longrightarrow A^{T, -}],$$
(See \S3 in loc.cit. for the relevant definitions.) and prove an equivariant main conjecture for its $p$--adic Tate module $T_p(\mathcal M_{S,T}^{\cH, -})$, viewed as a $\Bbb Z_p[[\mathcal G]]^-$--module. (See Theorems 4.6 and  5.6 in loc.cit.) In this section, we establish a link between the modules $T_p(\mathcal M_{S,T}^{\cH, -})$ and $Sel_S^T(\cH)_p^-$  and show that Theorems 4.6 and 5.6 in loc.cit. are direct consequences of our
Theorem \ref{EMC}. Throughout this section, we assume that $\mu=0$, $S_{ram}(\cH/F)\subseteq S$ and $T\neq\emptyset$.\\

\begin{proposition}\label{link-Selmer-Picard}
If $\mu=0$, $S_{ram}(\cH/F)\subseteq S$ and $T\neq\emptyset$, we have a canonical isomorphism
$$Sel_S^T(\cH)_p^-\, \simeq\, T_p(\mathcal M_{S,T}^{\cH, -})^\ast$$
in the category  of $\Bbb Z_p[[\mathcal G]]^-$--modules.
\end{proposition}
\begin{proof}
By Corollary 3.4 in \cite{Greither-Popescu}, we have commutative diagrams of $\Bbb Z_p[[\mathcal G]]^-$--modules
\[\begin{tikzcd}
(\cH_{S,T}^{(p^{m+1})}/(\cH_T^\times)^{p^{m+1}})^-\arrow{r}{\sim}\arrow{d}{\pi_m} & \mathcal M_{S,T}^{\cH, -}[p^{m+1}]\arrow{d}{\times p}\\
(\cH_{S,T}^{(p^m)}/(\cH_T^\times)^{p^m})^-\arrow{r}{\sim} & \mathcal M_{S,T}^{\cH,-}[p^m],
\end{tikzcd}\]
where the maps $\pi_m$ are described in Remark \ref{maps-pi_m} and the (inverses of the) horizontal isomorphisms are explicitly constructed in Proposition 3.2 of loc.cit. When combined with the last isomorphism in Remark \ref{maps-pi_m}, these lead to canonical isomorphisms of $\Bbb Z_p[[\mathcal G]]^-$--modules
$$ \varprojlim_m (\cH_{S,T}^{(p^m)}/(\cH_T^\times)^{p^m})^-\simeq T_p(\mathcal M_S^T(\cH)^-)\simeq T_p(\mathcal M_{S,T}^{\cH, -}).$$
When combined with Proposition \ref{link-sel-mst}(3), the last isomorphism above gives
$$Sel_S^T(\cH)_p^-\, \simeq\, T_p(\mathcal M_S^T(\cH)^-)^\ast \,\simeq\, T_p(\mathcal M_{S,T}^{\cH, -})^\ast,$$
which concludes the proof.
\end{proof}

Now, we are ready to derive the main result of \cite{Greither-Popescu} from our Theorem \ref{EMC}.

\begin{theorem}[Greither-Popescu, Theorems 4.6, 5.6 in \cite{Greither-Popescu}]\label{Greither-Popescu-Main} If $\mu=0$, $S_{ram}(\cH/F)\subseteq S$ and $T\neq\emptyset$,  then the following hold.
\begin{enumerate}
\item ${\rm pd}_{\Bbb Z_p[[\mathcal G]]^-} T_p(\mathcal M_{S, T}^{\cH, -}) =1.$
\item ${\rm Fitt}_{\Bbb Z_p[[\mathcal G]]^-} T_p(\mathcal M_{S, T}^{\cH, -}) =(\Theta_S^T(\cH/F)).$
\end{enumerate}
\end{theorem}
\begin{proof}
From Theorem \ref{EMC}(1), Proposition \ref{link-Selmer-Picard}, and Proposition 7.2(3) in \cite{Greither-Popescu}, we have
$${\rm Fitt}_{\Bbb Z_p[[\mathcal G]]^-} T_p(\mathcal M_{S, T}^{\cH, -}) = {\rm Fitt}_{\Bbb Z_p[[\mathcal G]]^-} T_p(\mathcal M_{S, T}^{\cH, -})^\ast= {\rm Fitt}_{\Bbb Z_p[[\mathcal G]]^-} Sel_S^T(\cH)^-_p=(\Theta_S^T(\cH/F)).$$
Further, since $\Theta_S^T(\cH/F)$ is not a zero--divisor in $\Bbb Z_p[[\mathcal G]]^-$ (see Proposition \ref{Theta-non-zero-divisor}), Proposition \ref{Cornacchia-Greither} implies that
$${\rm pd}_{\Bbb Z_p[[\mathcal G]]^-} T_p(\mathcal M_{S, T}^{\cH, -}) =1,$$
which concludes the proof.
\end{proof}
\smallskip

\subsection{The refined Coates--Sinnott Conjecture}
Let $K/k$ be an abelian extension of number fields, of Galois group $G$. We fix a finite set $S$  of places in $k$, satisfying the usual  $S_{ram}(K/k)\cup S_\infty(k)\subseteq S$ and let $K_i(\mathcal O_{K,S})$ denote the Quillen $K$--groups of the ring of $S$--integers $\mathcal O_{K,S}$, for all $i\geq 0$. We remind the reader (see \S6 in \cite{Greither-Popescu} and references therein) that $K_i(\mathcal O_{K,S})$ is a finitely generated $\Bbb Z[G]$--module for all $i$, and it is finite for all even $i\geq 2$. We let $K_i(\mathcal O_{K,S})^{tors}$ denote the torsion subgroup of $K_i(\mathcal O_{K,S})$, for all $i\geq 0$.
\\

Further, we follow \S6 in \cite{Greither-Popescu} and, for every $n\in\Bbb Z_{\geq 2}$, we define the idempotent element
$$e_n(K/k):=
\begin{cases}

\prod_{v\in S_\infty(k)}\frac{1+(-1)^n\sigma_v}{2}, & \text{ if $k$ is totally real},\\
0, &\text{ otherwise}
\end{cases}
$$
viewed in $\Bbb Z[1/2][G]$, where $\sigma_v$ is the Frobenius element of $v$ in $G$.  Recall that
$$e_n(K/k)\cdot\Theta_{S,K/k}(1-n)=\Theta_{S,K/k}(1-n),$$
and that $\Theta_{S, K/k}(1-n)$ is not a zero--divisor in $e_n(K/k)\Bbb Q[G]$, for all $n\geq 2$. (See equalities (35) in \S6 of \cite{Greither-Popescu}.) Consequently, if $T$ is a non-empty, finite set of primes in $k$, which is disjoint from $S\cup S_p(F)$, then
$\Theta^T_{S,K/k}(1-n)$ is also a non zero--divisor in $e_n(K/k)\Bbb Z_p[G]$, as
$$\Theta^T_{S, K/k}(1-n)=\delta_T(1-n)\Theta_{S, K/k}(1-n), \qquad\text{ for }\delta_T(1-n):=\prod_{v\in T}(1-\sigma_v^{-1}\cdot Nv^n)$$
and $\delta_T(1-n)$ is clearly a non zero--divisor in $\Bbb Z_p[G]$, for all $n\geq 1$.
\\

 In this section, we follow the strategy developed in \S6 of \cite{Greither-Popescu} and use our Theorem \ref{EMC} to give an unconditional proof of the following refinement of the Coates--Sinnott Conjecture, away from its $2$--primary part.
\begin{theorem}[The refined Coates--Sinnott Conjecture]\label{CS-K}
For $(K/k, S)$ as above, we have
$$(Fitt_{\Bbb Z_p[G]} K_{2n-1}(\mathcal O_{K,S})_p^{tors})\cdot\Theta_{S,K/k}(1-n)=e_n(K/k)\cdot Fitt_{\Bbb Z_p[G]} K_{2n-2}(\mathcal O_{K,S})_p$$
for all $n\in\Bbb Z_{\geq 2}$ and all primes $p>2$.
\end{theorem}
\begin{remark} The original Coates--Sinnott conjecture states that the left--hand side of the equality in the theorem above is included in $Ann_{\Bbb Z_p[G]} K_{2n-2}(\mathcal O_{K,S})_p$, for all primes $p$. The theorem above obviously implies this original statement, at primes $p>2$.
\end{remark}
As discussed in detail in \S6 of \cite{Greither-Popescu}, work of Soul\'e  and Dwyer-Friedlander (late 1970s, early 1980s) lead to the construction of canonical surjective $\Bbb Z_p[[\mathcal G]]$--linear morphisms
$$ch_{p,n}^i:K_{2n-i}(O_{K,S})\twoheadrightarrow H^i_{et}(\mathcal O_{K,S}[1/p], \Bbb Z_p(n)), \qquad\text{ for all } i=1,2, \quad n\geq 2,\quad p>2,$$
called $p$--adic Chern character maps, which were conjectured by Quillen--Lichtenbaum and Dwyer--Friedlander to be isomorphisms, and this was confirmed to hold true by more recent work of Rost and Vojevodsky. Consequently, Theorem \ref{CS-K} above is equivalent to the following \'etale--cohomological version.
\begin{theorem}\label{CS-H}
For $(K/k, S, n, p)$ as above, we have an equality
$$(Fitt_{\Bbb Z_p[G]} H^1_K(S,p,n)^{tors})\cdot\Theta_{S,K/k}(1-n)=e_n(K/k)\cdot Fitt_{\Bbb Z_p[G]} H^2_K(S,p,n),$$
where $H^i_K(S,p,n):=H^i_{et}(\mathcal O_{K,S}[1/p], \Bbb Z_p(n))$, for all $i=1,2$.
\end{theorem}
In \S6 of \cite{Greither-Popescu} it is proved (see Lemmas 6.12, 6.13, 6.14, 6.15 in loc.cit.) that it suffices to prove the theorem above under the following hypotheses.
\begin{itemize}
\item $K$ is CM and $k$ is totally real.
\item $S_p(F)\subseteq S$.
\item $\boldsymbol\mu_p\subseteq K$, where $\boldsymbol\mu_p$ is the set of $p$--roots of unity. (Equivalently, $\mathcal K:=K(\boldsymbol\mu_{p^\infty})$ is the cyclotomic $\Bbb Z_p$--extension of $K$.)
\end{itemize}
\medskip

Under these hypotheses, {\it and assuming further that the Iwasawa $\mu$--invariant $\mu_{K,p}$ vanishes}, the authors of \cite{Greither-Popescu} develop a technical machinery which allows them to link the modules $T_p(\mathcal M^{\mathcal K,-}_{S,T})$ to the cohomology groups $H^i_K(S,p, n)$ and use the information provided by Theorem \ref{Greither-Popescu-Main} above for the data $(K/k,S, p)$ and varying sets $T$ to derive Theorem \ref{CS-H}.
\medskip

In this section, {\it we remove the hypothesis $\mu_{K,p}=0$} and, instead of using the modules $T_p(\mathcal M_{S,T}^{\mathcal K, -})$ and Theorem \ref{Greither-Popescu-Main}, we feed the modules $Sel_S^T(\mathcal K)_p^-$ and our Theorem \ref{EMC} through the technical machinery developed in \S6 of \cite{Greither-Popescu} and prove Theorem \ref{CS-H} unconditionally. To that end, from now on we fix data $(K/k, S, p, n)$ as above, assume the three hypotheses listed above and, for further notational simplicity, we set:
$$e_n:=\widetilde{e_n(K/k)}:=\frac{1+(-1)^nj}{2}, \qquad H^i_K:=H^i_K(S,p,n), \text{ for }i=1,2,$$
where $j$ is the complex conjugation in $\mathcal G$. Note that we have ring isomorphisms
$$ e_n\Bbb Z_p[[\mathcal G]] \simeq\Bbb Z_p[[\mathcal G]]/(1+(-1)^n j), \qquad t_{1-n}:\Bbb Z_p[[\mathcal G]]^-\simeq e_n\Bbb Z_p[[\mathcal G]],$$
 where $t_m:\Bbb Z_p[[\mathcal G]]\simeq \Bbb Z_p[[\mathcal G]]$ is the continuous $\Bbb Z_p$--algebra isomorphism such that
 $$t_m(\sigma)=\kappa(\sigma)^m\cdot\sigma,\qquad \text{ for all }\sigma\in\mathcal G,\, m\in\Bbb Z,$$
 and $\kappa:\mathcal G\to\Bbb Z_p^\times$ is the $p$--adic cyclotomic character.
Throughout, the reader should keep in mind that all our duals are endowed with the co-variant Galois action, as before.
\\

The first step in this process consists of establishing a link between the \'etale cohomology group $H^2_K$ and the Selmer group $Sel_S^\emptyset(\mathcal K)_p^-$.  This is given by part (1) of the next Proposition.

\begin{proposition}\label{link-etale-selmer}
Under the above hypotheses, if we let $\Gamma:=Gal(\mathcal K/K)$, we have the following canonical $\Bbb Z_p[[\mathcal G]]$--module isomorphisms.
\begin{enumerate}
\item $e_n\cdot H^2_K\simeq (Sel_S^\emptyset(\mathcal K)_p^-(n-1)_\Gamma)^\vee$.
\item $e_n\cdot (H^1_K)^{tors}=(H^1_K)^{tors}\simeq \Bbb Z_p(n)_\Gamma$.
\end{enumerate}
\end{proposition}
\begin{proof} Part (2) is Proposition 6.17(2) in \cite{Greither-Popescu}.

In order to prove part (1), we make the notational convention that if $M$ is a $\Bbb Z_p[[\mathcal G]]$--module, then $\overline M$ is the $\Bbb Z_p[[\mathcal G]]$--module which is identical to $M$ as a $\Bbb Z_p$--module, but it is endowed with
the {\bf contravariant} $\mathcal G$--action
$$(\sigma, m)\mapsto \sigma^{-1}\cdot m, \qquad \text{ for all }\sigma\in\mathcal G,\,\, m\in M.$$
The identity map induces isomorphisms of $\Bbb Z_p[[\mathcal G]]$--modules $\overline{\overline M}\simeq M$, $\overline{M^\vee}\simeq (\overline M)^\vee$,  $\overline{M^\ast}\simeq (\overline M)^\ast$, $(\overline M)_{\Gamma}\simeq \overline {M_\Gamma}$, and $\overline{M(m)}\simeq \overline M(-m)$, for all $m\in\Bbb Z$.

Let $\mathfrak X_S$ denote the Galois group of the maximal abelian, pro-$p$ extension $\mathcal K_S/\mathcal K$ of $\mathcal K$, which is unramified away from $S$.
Since $\boldsymbol\mu_{p^\infty}\subseteq \mathcal K$, $\mathcal K_S/\mathcal K$ is a Kummer extension, and we have a perfect, $\Bbb Z_p[[\mathcal G]]$--bilinear Kummer pairing (see \S3.3. in \cite{Greither-Popescu} as well)
$$\mathfrak X_S\times \mathcal M_S^\emptyset(\mathcal K)\longrightarrow \Bbb Q_p/\Bbb Z_p(1)=\boldsymbol\mu_{p^\infty}, \qquad (\sigma, \widehat x)\mapsto \frac{\sigma(x^{1/p^m})}{x^{1/p^m}},$$
for all $m\geq 1$, all $\widehat x\in \mathcal M_S^\emptyset(\mathcal K)[p^m]$ (i.e. $x\in \mathcal K_{S,\emptyset}^{(p^m)}$ and $\widehat x\in K_{S,\emptyset}^{(p^m)}/\mathcal K^{\times p^m}$), and all $\sigma\in\mathfrak X_S$. The above perfect pairing gives a $\Bbb Z_p[[\mathcal G]]$--module isomorphism
$$\overline{\mathfrak X_S^+(-1)}\simeq \mathcal M_S^\emptyset(\mathcal K)^{-,\vee}.$$
When combined with Proposition \ref{link-sel-mst}(2), this leads to the $\Bbb Z_p[[\mathcal G]]$--module isomorphisms
$$ \overline{\mathfrak X_S^+(-1)}\simeq Sel_S^\emptyset(\mathcal K)_p^-, \qquad \overline{\mathfrak X_S^+(-n)}\simeq Sel_S^\emptyset(\mathcal K)_p^-(n-1).$$
When combined with the first isomorphism in Proposition 6.17(1) of \cite{Greither-Popescu}, this gives
$$e_n\cdot H^2_K\simeq \overline{(\mathfrak X_S^+(-n)_\Gamma)^\vee}\simeq \left(\left(\overline{\mathfrak X_S^+(-n)}\right )_\Gamma\right )^\vee\simeq (Sel_S^\emptyset(\mathcal K)_p^-(n-1)_\Gamma)^\vee,$$
which concludes the proof of part (1).

Note that the first isomorphism in Proposition 6.17(1) of loc.cit is written there as
$$e_n\cdot H^2_K\simeq (\mathfrak X_S^+(-n)_\Gamma)^\vee.$$
That is because in \S6 of loc.cit. the duals are endowed with the contravariant Galois action, whereas here we are viewing them as endowed with the covariant Galois action. Also note that
in loc.cit. it is also stated that $(\mathfrak X^+(-n)_\Gamma)^\vee\simeq (\mathfrak X^+(-n)^\ast)_\Gamma$. This last isomorphism only holds if $\mu=0$ (see proof in loc.cit.), which we are not assuming here.
\end{proof}
\medskip

Next, we fix a finite, non-empty set $T$ of non--archimedean primes in $F$, disjoint from $S$. For simplicity, we assume that $T=\{v\}$, although the arguments below work for sets $T$ of arbitrary size.
 As in \S2 of \cite{Greither-Popescu}, we consider the $\Bbb Z[\mathcal G]$--module (denoted in loc.cit. $\Delta_{\mathcal K, \mathcal T}$)
$$\Delta_T(\mathcal K):=\bigoplus_{w\in T(\mathcal K)}\kappa(w)^\times\simeq \left(\kappa(\overline v)^\times\otimes_{\Bbb Z[\mathcal G_v]}\Bbb Z[\mathcal G]\right)$$
where $\kappa(w)$ is the residue field of the prime $w$ of $\mathcal K$, $\overline v$ is a fixed prime of $\mathcal K$ sitting above $v$, and $\mathcal G_v$ is the decomposition group of $v$ in $\mathcal G$.
Note that since $S_p\cap T=\emptyset$ and $\boldsymbol \mu_{p^\infty}\subseteq \mathcal K$, $\boldsymbol\mu_{p^\infty}$ injects in $\kappa(w)^\times$ via the reduction mod $w$ map and this injection identifies the $\Bbb Z_p[[\mathcal G_v]]$--module  $\boldsymbol\mu_{p^\infty}$ with the $p$--primary part $\kappa(w)^\times_p$ of $\kappa(w)^\times$, for all $w\in T(\mathcal K)$. (See the proof of Lemma 2.8 in loc.cit.) Consequently, we have an isomorphism of $\Bbb Z_p[[\mathcal G]]$--modules
$$\Delta_T(\mathcal K)_p\simeq \left(\Bbb Q_p/\Bbb Z_p(1)\otimes_{\Bbb Z_p[[\mathcal G_v]]}\Bbb Z_p[[\mathcal G]]\right).$$
So, $\Delta_T(\mathcal K)_p$ is a divisible, $p$--torsion module of corank $\vert T(\mathcal K)\vert$, which contains $\mu_{p^\infty}=\Bbb Q_p/\Bbb Z_p(1)$ (via the diagonal reduction mod $w$ embedding, for all $w\in T(\mathcal K)$). Consequently, we have an isomorphism of abstract groups
$$\Delta_T(\mathcal K)_p^-\simeq (\Bbb Q_p/\Bbb Z_p)^{\delta_{T,\mathcal K}^-}, \qquad\text{ where } 1\leq \delta_{T,\mathcal K}^-\leq \vert T(\mathcal K)\vert.$$
As a consequence, we have an  isomorphism of $\Bbb Z_p[[\mathcal G]]$--modules
\begin{equation}\label{Delta-dual} \Delta_T(\mathcal K)_p^{\vee, -}\simeq T_p(\Delta_T(\mathcal K))^{\ast, -}.
\end{equation}
Now, we have a canonical exact sequence of $\Bbb Z_p[[\mathcal G]]$--modules (see \S2 of loc.cit.)
$$
0\longrightarrow \Bbb Q_p/\Bbb Z_p(1)\overset{d_T}\longrightarrow \Delta_T(\mathcal K)_p^-\longrightarrow A^{T,-}\overset{\pi_T}\longrightarrow A^{\emptyset, -}\longrightarrow 0,
$$
where $\pi_T$ maps the class $\widehat D$ in $A^{T, -}$ of a divisor $D$ to the class of $D$ in $A^{\emptyset, -}$, and $d_T$ is the diagonal mod $w$ embedding, for $w\in T(\cH)$.  After taking Pontrjagin duals and using \eqref{Delta-dual}, we get an exact sequence of $\Bbb Z_p[[\mathcal G]]$--modules
$$
0\longrightarrow A^{\emptyset, \vee,-}\overset{\pi_T^\vee}\longrightarrow A^{T, \vee, -}\longrightarrow  T_p(\Delta_T(\mathcal K))^{\ast, -}\overset{d_T^\vee} \longrightarrow \Bbb Z_p(1)^\ast\longrightarrow 0.
$$
Now, we have an obvious commutative diagram of $\Bbb Z_p[[\mathcal G]]$--modules

\[ \begin{tikzcd}
0\arrow{r} & Y_{S\setminus S_p}(\mathcal K)_p^{\ast, -}\arrow{r}\arrow{d}{{\rm id}}[swap]{=} & Sel_S^\emptyset(\mathcal K)_p^-\arrow{r}\arrow{d} & A^{\emptyset,\vee,-}\arrow{r}\arrow{d}{\pi_T^\vee} & 0\\
0\arrow{r}  &Y_{S\setminus S_p}(\mathcal K)_p^{\ast, -}\arrow{r} &Sel_S^T(\mathcal K)_p^-\arrow{r} & A^{T, \vee, -}\arrow{r} & 0,
\end{tikzcd}
\]
where the middle vertical map is induced by the inclusion $\mathcal K_T^\times\subseteq \mathcal K^\times.$ Now, when applying the snake lemma to the diagram above, combining it with the exact sequence just above it, and noting that $\Bbb Z_p(1)^\ast\simeq\Bbb Z_p(1)$ (as duals are endowed with the covariant $\mathcal G$--action), we  obtain an exact sequence of $\Bbb Z_p[[\mathcal G]]$--modules
\begin{equation}\label{long-seq-sel}
0\longrightarrow  Sel_S^\emptyset(\mathcal K)_p^-\longrightarrow Sel_S^T(\mathcal K)_p^-\longrightarrow T_p(\Delta_T(\mathcal K))^{\ast, -}\overset{d_T^\vee} \longrightarrow \Bbb Z_p(1)\longrightarrow 0.
\end{equation}
This is the unconditional analogue of the dual of exact sequence (15) in loc.cit. If we tensor the above sequence with $\Bbb Z_p(n-1)$, we obtain the following exact sequence
\begin{equation}\label{long-seq-sel-twist}
0\longrightarrow  Sel_S^\emptyset(\mathcal K)_p^-(n-1)\longrightarrow Sel_S^T(\mathcal K)_p^-(n-1)\longrightarrow (T_p(\Delta_T(\mathcal K))^-(n-1))^{\ast} \overset{d_T^\vee(n-1)}\longrightarrow \Bbb Z_p(n)\longrightarrow 0,
\end{equation}
which is the unconditional analogue of the dual of exact sequence (39) in loc.cit. Now, we are ready to prove Theorem \ref{CS-H}, by feeding our exact sequence \eqref{long-seq-sel-twist} through the machinery developed in \cite{Greither-Popescu}  (in the paragraphs right after exact sequence (39) in loc.cit.)\\

\begin{proof}[Proof of Theorem \ref{CS-H}] As in loc.cit.,  the idea is to take $\Gamma$--coinvariants in exact sequence \eqref{long-seq-sel-twist}, show that the resulting sequence of $\Gamma$--coinvariants  is exact, that all its terms are finite, and the middle two terms are of
${\rm pd}_{\Bbb Z_p[G]}\leq 1$ and then apply Lemma \ref{Burns-Greither} to derive the desired equality of Fitting ideals.
We start by noting that for all finitely generated $\Bbb Z_p[[\mathcal G]]$--modules $X$, we have equalities
\begin{equation}\label{module-twist}
e_n(X^-(n-1))=X^-(n-1), \qquad Fitt_{e_n\Bbb Z_p[[\mathcal G]]}(X^-(n-1))=t_{1-n}( Fitt_{\Bbb Z_p[[\mathcal G]]^-}X^-),
\end{equation}
where the first is obvious and the second is a direct consequence of Lemma 7.4(2) in loc.cit.
\begin{lemma}\label{sel-twist-coinv} For $(K/k,S,T,p, n)$ as above, let $N:=Sel_S^T(\mathcal K)_p^-(n-1)$. Then:
\begin{enumerate}
\item  $Fitt_{e_n\Bbb Z_p[G]}(N_\Gamma)=(\Theta_{S,K/k}^T(1-n))$.
\item ${\rm pd}_{e_n\Bbb Z_p[G]} N_{\Gamma}\leq 1.$
\item $N_{\Gamma}$ is finite and quadratically presented over $e_n\Bbb Z_p[G]$.
\end{enumerate}
\end{lemma}
\begin{proof} Our Theorem \ref{EMC}(1) combined with \eqref{module-twist} above give an equality
$$Fitt_{e_n\Bbb Z_p[[\mathcal G]]} N=(t_{1-n}(\Theta_S^T(\mathcal K/k))).$$
However, if $\pi:\Bbb Z_p[[\mathcal G]] \to\Bbb Z_p[G]$ is the $\Bbb Z_p$--algebra surjection given by Galois restriction, the equality above and Lemma 5.14(2) in loc.cit gives us
$$
 Fitt_{e_n\Bbb Z_p[G]}(N_{\Gamma})=\pi(Fitt_{e_n\Bbb Z_p[[\mathcal G]]}N)=(\pi\circ t_{1-n}(\Theta_S^T(\mathcal K/k)))=(\Theta_{S,K/k}^T(1-n)),
$$
which concludes the proof of (1).

 Since $\Theta_S^T(1-n)$ is not a zero--divisor in $e_n\Bbb Z_p[G]$, $N_\Gamma$ is finitely generated over $e_n\Bbb Z_p[G]$ (as  $N$ is finitely generated over $\Bbb Z_p[[\mathcal G]]^-$),  and   $e_n\Bbb Z_p[G]$ is a semilocal ring, Proposition \ref{Cornacchia-Greither} and part (1) imply part (2) and the quadratic presentation statement in part (3).

 From part (1), we conclude that $\Theta_{S, K/k}^T(1-n)$ annihilates $N_\Gamma$. Therefore, $N_\Gamma$ is a isomorphic to a quotient of a direct sum of finitely many copies of
 $$e_n\Bbb Z_p[G]/(\Theta_{S, K/k}^T(1-n)),$$
 which is easily seen to be finite, as $\Theta_{S,K/k}^T(1-n)$ is not a zero divisor in $e_n\Bbb Z_p[G]$ and $e_n\Bbb Z_p[G]$ is a finite, integral extension of $\Bbb Z_p$. This concludes the proof of part (3).
\end{proof}
\begin{lemma}\label{delta-twist-coinv}For $(K/k,S,T,p, n)$ as above, let $M:= T_p(\Delta_T(\mathcal K))^-(n-1)$. Then:
\begin{enumerate}
\item  $Fitt_{e_n\Bbb Z_p[G]} (M^{\ast})_\Gamma =(e_n\delta_T(1-n))$.
\item ${\rm pd}_{e_n\Bbb Z_p[G]}(M^\ast)_\Gamma \leq 1.$
\item $(M^\ast)_\Gamma$ is finite and quadratically presented over $e_n\Bbb Z_p[G]$.
\end{enumerate}
\end{lemma}
\begin{proof} First, note that (1), (2) and (3) are all true if we replace $M^\ast$ by $M$ in the above statements. These are all immediate consequences of the isomorphism of $e_n\Bbb Z_p[G]$--modules
\begin{equation}\label{iso1-delta_T}M_{\Gamma}\simeq e_n\Bbb Z_p[G]/e_n\delta_T(1-n).\end{equation}
and the fact that $e_n\delta_T(1-n)$ is not a zero--divisor in $e_n\Bbb Z_p[G]$. See the displayed isomorphism right below (39) in \S6 of \cite{Greither-Popescu} and the resulting arguments therein.) Now, since $M$ is $\Bbb Z_p$--free and $M_\Gamma$ is finite, Lemma 6.16 in loc.cit gives an isomorphism of $e_n\Bbb Z_p[G]$--modules
$$\label{iso2-delta_T}(M^\ast)_{\Gamma}\simeq (M_\Gamma)^\vee.$$
On the other hand, the ring $R:=e_n\Bbb Z_p[G]/e_n\delta_T(1-n)$ is a finite, semi-local Gorenstein $\Bbb Z_p$--algebra (as a finite quotient of $\Bbb Z_p[G]$.) Therefore, we have isomorphisms of $R$--modules
$$(M_\Gamma)^\vee  \simeq {\rm Hom}_R(M_\Gamma, R)\simeq R,$$
where the first follows from Proposition 4 in the Appendix of \cite{Mazur-Wiles} and the second follows from \eqref{iso1-delta_T}. The last three displayed isomorphisms show that
$$(M^\ast)_{\Gamma}\simeq M_\Gamma$$
as $R$--modules and, consequently, as $e_n\Bbb Z_p[G]$--modules, which concludes the proof.
\end{proof}
Next, we will take $\Gamma$--coinvariants in exact sequence \eqref{long-seq-sel-twist}. With notations as in the previous two Lemmas, note that since $M^\ast$ and $\Bbb Z_p(n)$ are $\Bbb Z_p$--free of finite rank and $(M^\ast)_\Gamma$ and $\Bbb Z_p(n)_\Gamma$ are finite, we have
$$(M^\ast)^{\Gamma}=0, \qquad \Bbb Z_p(n)^{\Gamma}=0, \qquad (\ker d_T^\vee(n-1))^\Gamma=0,$$
where the first two equalities follow from Lemma 6.16 of \cite{Greither-Popescu} and the third follows from the first. This fact combined with Proposition \ref{link-etale-selmer} implies right away that when taking $\Gamma$--coinvariants in \eqref{long-seq-sel-twist}, we get a long exact sequence of
$e_n\Bbb Z_p[G]$--modules
\begin{equation}\label{coinvariant-seq} 0\to (e_nH_K^2)^\vee \to N_\Gamma\to (M^\ast)_\Gamma\to (H^1_K)^{tors}\to 0,\end{equation}
whose terms are all finite and whose middle two terms are of ${\rm pd}_{e_n\Bbb Z_p[G]}\leq 1$ (see Lemma \ref{sel-twist-coinv} and Lemma \ref{delta-twist-coinv}.) Now, since the ring $R:=e_n\Bbb Z_p[G]$ is a $\Bbb Z_p$--algebra which is $\Bbb Z_p$--free of finite rank and relatively Gorensteing over $\Bbb Z_p$, we can apply Lemma \ref{Burns-Greither} to the above exact sequence and combine that with Lemma \ref{sel-twist-coinv}(1) and Lemma \ref{delta-twist-coinv}(1) to  obtain
$${\rm Fitt}_{e_n\Bbb Z_p[G]}(e_nH^2_K)\cdot (e_n\delta_T(1-n))= ({\rm Fitt}_{e_n\Bbb Z_p[G]}(H^1_K)^{tors})\cdot\Theta_{S, K/k}(1-n)\cdot (e_n\delta_T(1-n)).$$
Since $e_n\delta_T(1-n)$ is not a zero--divisor in $e_n\Bbb Z_p[G]$, we can cancel it out to obtain
$$e_n{\rm Fitt}_{\Bbb Z_p[G]}(H^2_K)=({\rm Fitt}_{\Bbb Z_p[G]}(H^1_K)^{tors})\cdot \Theta_{S,K/k}(1-n),$$
recalling that $e_n\Theta_{S,K/k}(1-n)=\Theta_{S,K/k}(1-n)$. This concludes the proof of Theorem \ref{CS-H}.
\end{proof}
\smallskip

\subsection{A conjecture of Burns--Kurihara--Sano} In this section, we prove a stronger version of our Theorem \ref{keystone-improved}, which implies a conjecture of Burns--Kurihara--Sano. In what follows, $(H/F, S, T, p)$ are as in Theorem 
\ref{keystone-improved}. The main goal is to replace hypotheses ${\rm Hyp}(H_\infty/F)_p$ with the weaker hypotheses ${\rm Hyp}(H/F)_p$, the difference being that $S$ will only be required to contain the set $S_p^{ram}(H/F):=S_p(F)\cap S_{ram}(H/F)$ of $p$--adic primes which ramify in $H/F$ as opposed to the entire set of $p$--adic primes $S_p(F)=S_p(F)\cap S_{ram}(H_\infty/F)$. More precisely, our goal is to prove the following.

\begin{theorem}\label{BKS-Conj}The following hold.
\begin{enumerate}
\item If $p$ is an odd prime and $(H/F, S, T, p)$ satisfy ${\rm Hyp}(H/F)_p$, then we have
$${\rm Fitt}_{\mathbb{Z}_p[G]^-}(Sel_S^T(H)_p^-)=(\Theta_S^T(H/F)).$$
\item {\bf (The Burns--Kurihara--Sano Conjecture.)}  If $(H/F, S, T)$ satisfy ${\rm Hyp}(H/F)_p$, for all odd primes $p$ (e.g. if $S_{ram}(H/F)\subseteq S$ and $H_T^\times$ has no coprime--to--$2$ torsion), then we have
$${\rm Fitt}_{\mathbb{Z}[G]^-}(Sel_S^T(H)^-)=(\Theta_S^T(H/F)).$$
\end{enumerate}

\end{theorem}
\begin{proof}
Obviously, part (2) is a consequence of part (1). In order to prove part (1), let us fix $(K/k, S, T, p)$ satisfying hypotheses ${\rm Hyp}(H/F)_p$. Let $S_p^{ur}(H/F):=S_p(F)\setminus S_{ram}(H/F)$. Note that it suffices to prove the required equality of ideals under the assumption that $T\cap S_p^{ur}(H/F)=\emptyset.$ Indeed, if this is not the case, then we let $T':=T\setminus S_p^{ur}(H/F)$ and observe that $(K/k, S, T', p)$ still satisfies ${\rm Hyp}(H/F)_p$. Now, Lemma 4.1 in \cite{Dasgupta-Kakde} shows that we have equalities of ideals 
$${\rm Fitt}_{\mathbb{Z}_p[G]^-}(Sel_S^T(H)_p^-)={\rm Fitt}_{\mathbb{Z}_p[G]^-}(Sel_S^{T'}(H)_p^-), \qquad (\Theta_S^T(H/F))=(\Theta_S^{T'}(H/F)),$$
showing that the equality in (1) for $(H/F,S,T',p)$ is equivalent to the same equality for $(H/F,S,T,p)$. So, from now on, we assume that $T\cap S_p^{ur}(F)=\emptyset$.

Let $S_\ast:=S\cup S_p^{ur}(H/F)$ and $u_p(S, H/F):=\prod_{v\in S_\ast\setminus S}(1-\sigma_v^{-1})$.  Since $(H/F, S_\ast, T, p)$ satisfies ${\rm Hyp}_p(H_\infty/F)$, Theorem \ref{keystone-improved} gives us an equality of $\Bbb Z_p[G]^-$--ideals
\begin{equation}\label{Fitt-S_ast}{\rm Fitt}_{\mathbb{Z}_p[G]^-}(Sel_{S_\ast}^T(H)_p^-)=(\Theta_{S_\ast}^T(H/F))=u_p(S,H/F)\cdot (\Theta_{S}^T(H/F)).\end{equation}
On the other hand, we have an exact sequence of $\Bbb Z_p[G]^-$--modules
$$0\to (Y_{S_\ast\setminus S}(H)_p^-)^\ast\to Sel_{S_\ast}^T(H)_p^-\to Sel_{S}^T(H)_p^-\to 0,$$
(see \eqref{enlarge-S-Selmer-sequence}), where the right--most nonzero module is quadratically presented. (See Lemma \ref{quadratic-pres-lemma}.) Consequently, Lemma \ref{Fitt-ses} and Lemma \ref{Fitt-Y_S-dual}(2) give us an equality 
\begin{equation}\label{Fitt-S_ast-S}{\rm Fitt}_{\mathbb{Z}_p[G]^-}(Sel_{S_\ast}^T(H)_p^-)= u_p(S, H/F)\cdot {\rm Fitt}_{\mathbb{Z}_p[G]^-}(Sel_{S}^T(H)_p^-).\end{equation}
The equalities  of $\Bbb Z_p[G]^-$--ideals \eqref{Fitt-S_ast} and \eqref{Fitt-S_ast-S} show that we have
\begin{equation}\label{u-times-equality}u_p(S, H/F)\cdot {\rm Fitt}_{\mathbb{Z}_p[G]^-}(Sel_{S}^T(H)_p^-) = u_p(S,H/F)\cdot (\Theta_{S}^T(H/F)).\end{equation}
Therefore, proving the equality in part  (1) would mean cancelling the $u_p(S, H/F)$--factor in the equality above. Unfortunately, in general, this factor is a zero--divisor in $\Bbb Z_p[G]^-$, which makes this problem quite difficult. We will eliminate this factor by applying what is now called ``the method of Taylor--Wiles primes'', introduced first by Wiles in \cite{Wiles-Brumer} and refined further by Greither in \S4 of \cite{Greither-Brumer}. We proceed in two steps.
\\

{\bf Step 1.} First, we intend to prove that, under the above hypotheses, we always have an inclusion 
\begin{equation}\label{inclusion1} (\Theta_S^T(H/F))\subseteq {\rm Fitt}_{\mathbb{Z}_p[G]^-}(Sel_S^T(H)_p^-).\end{equation}
Following \S4 in \cite{Greither-Brumer}, we give the following definition.
\begin{definition}
An odd prime $p$ is called $H/F$--bad if there exists $v\in S_p(F)$ such that $j\not\in G_v(H/F)$ or, equivalently, if there exists $v\in S_p(F)$ which splits completely in $H/H^+$, where $H^+:=H^{j=1}$. 
An odd prime which is not $H/F$--bad will be called $H/F$--good.
\end{definition}
\begin{remark}\label{remark-bad} Note that if the prime $p>2$ is $H/F$--good, then $(1-\sigma_v^{-1})$ is a divisor of $(1-j)$ in $\Bbb Z_p[G]^-$, for all $v\in S_p^{ur}(H/F)$. Since $(1-j)$ is invertible in $\Bbb Z_p[G]^-$,  $u_p(S, H/F)$ is also invertible in $\Bbb Z_p[G]^-$. Therefore, equality \eqref{u-times-equality} implies the desired result. 
\end{remark}
As a conseqeunce of the above Remark, we may and will assume that $p$ is $H/F$--bad. 
The following is proved in \cite{Greither-Brumer}, as an application of the Tchebotarev density theorem. (See Proposition 4.1 in loc.cit. and note that the assumption that $H/F$ ($K/F$ in in loc.cit.) is a ``nice'' extension is never used in the proof.)
\begin{proposition}[Greither \cite{Greither-Brumer}]\label{Tchebotarev}  Let $N\in\Bbb Z_{>0}$. Then, there exist infinitely many primes $r>2$, such that:
\begin{enumerate}
\item $r\equiv 1\mod p^N$.
\item $r$ is $H/F$--good.
\item The Frobenius morphism $\sigma_p(E/\Bbb Q)$ generates $G(E/\Bbb Q)$, where $E\subseteq\Bbb Q(\mu_r)$ such that $[E:\Bbb Q]=p^N$. 
\end{enumerate}
\end{proposition} 
Next, we fix two integers $N, M\in\Bbb Z_{>0}$, such that 
\begin{equation}\label{N-M}N-M\geq {\rm ord}_p(|G|\cdot f),\end{equation}
where $f$ is the least common multiple of the residual degrees $f(v/p):=[\mathcal O_F/v:\Bbb F_p]$, for all $v\in S_p(F)$. Further, for the chosen integer $N$, we pick an odd prime $r$  as in Proposition \ref{Tchebotarev}, which satisfies the additional hypotheses (excluding only finitely many choices):
$$r\not\in S_{ram}(H/\Bbb Q), \qquad S_r(F)\cap (S\cup T)=\emptyset.$$
Once these choices are made, we take the compositum $H':=H\cdot E$, where $E$ is the unique (totally real) field described in Proposition \ref{Tchebotarev}(3) and observe that $H'/F$ is an abelian, CM--extension of $F$, of Galois group
$$G(H'/F)\simeq G\times \Delta, \qquad \text{ where } \Delta:=G(E/\Bbb Q)=\langle \sigma_p(E/\Bbb Q)\rangle \simeq \Bbb Z/p^N\Bbb Z.$$
It is also easy to see that $S_p^{ur}(H'/F)=S_p^{ur}(H/F)$ and that, if we set 
$$S':=S\cupdot S_r(F) \text{ (disjoint union) }, \qquad  S'_\ast:=S'\cup S_p^{ur}(F)=S_\ast\cupdot S_r(F),$$ 
then $(H'/F, S', T, p)$ satisfies hypotheses ${\rm Hyp}(H'/F)_p$ and  $(H'/F, S_\ast', T, p)$ satisfies hypotheses ${\rm Hyp}(H_\infty'/F)_p$. Consequently, we have an equality of $\Bbb Z_p[G(H'/F)]^-$--ideals similar to \eqref{u-times-equality}:
\begin{equation}\label{u'-times-equality}u_p(S', H'/F)\cdot {\rm Fitt}_{\mathbb{Z}_p[G']^-}(Sel_{S'}^T(H')_p^-) = u_p(S',H'/F)\cdot (\Theta_{S'}^T(H'/F)), \end{equation}
where $u_p(S', H'/F):=\prod_{v\in S'_\ast\setminus S'}(1-\sigma_v(H'/F)^{-1}).$ 

 If  $\nu\in\Bbb Z_p[\Delta]\subseteq  \Bbb Z_p[G(H'/F)]^-$ is given by
$$\nu:=1+\sigma_p(E/\Bbb Q)^{p^{N-M}}+ (\sigma_p(E/\Bbb Q)^{p^{N-M}})^2+\dots + (\sigma_p(E/\Bbb Q)^{p^{N-M}})^{p^M-1},$$
then the following is proved in \cite{Greither-Brumer} (see Proposition 4.6 in loc.cit.), as an immediate consequence of \eqref{N-M}.
\begin{proposition}[Greither \cite{Greither-Brumer}]\label{Greither-nzd} The element $u_p(S', H'/F)\in \Bbb Z_p[G(H'/F)]^-$ becomes a non--zero--divisor when projected to $\Bbb Z_p[G(H'/F)]^-/\nu$.
\end{proposition}
For simplicity, let us set $R:=\Bbb Z_p[G]^-$, $R':=\Bbb Z_p[G(H'/F)]^-$. Note that, as a direct consequence of the definition of $\nu$,  the canonical surjective ring morphism  $\rho: R'\twoheadrightarrow R$ (induced by Galois restriction, taking every element in $\Delta$ to $1$) induces a surjective ring morphism $\widehat{\rho} :R'/\nu\twoheadrightarrow R/p^M$. This is ilustrated in the following commutative diagram of surjective ring morphisms, where the vertical maps are the usual projections.
\[\begin{tikzcd}
R'\arrow{r}{\rho}\arrow{d}{\pi_{R'}} & R\arrow{d}{\pi_R}\\
R'/\nu\arrow{r}{\widehat{\rho}} & R/p^M
\end{tikzcd}\]
Now, we apply $\pi_{R'}$ to equality \eqref{u'-times-equality}, and use Proposition \ref{Greither-nzd} above ($\pi_{R'}(u_p(S', H'/F))$ is a non--zero--divisor) together with property \eqref{base-change-Fitt} of Fitting ideals, to conclude that we have an equality of $R'/\nu$--ideals
$${\rm Fitt}_{R'/\nu}\,(Sel_{S'}^T(H')_{R'/\nu}) = (\Theta_{S'}^T(H'/F))_{R'/\nu}.$$
Next, we apply $\widehat{\rho}$ to the equality above and use \eqref{base-change-Fitt} again to obtain an equality of $R/p^M$--ideals
$${\rm Fitt}_{R/p^M}\,(Sel_{S'}^T(H')_{R/p^M}) = (\Theta_{S'}^T(H'/F))_{R/p^M}.$$
If we take the preimage via $\pi_R$ of the equality above and use \eqref{base-change-Fitt} again, we obtain equalities of $R$--ideals
\begin{equation}\label{R-Fitting-equality}{\rm Fitt}_{R}\,(Sel_{S'}^T(H')_{R}) +p^M R = (\Theta_{S'}^T(H'/F))_{R}+p^MR= (\Theta_{S'}^T(H/F)) +p^MR,\end{equation}
where the last equality is the obvious restriction property of Artin $L$--functions $\rho(\Theta_{S'}^T(H'/F))=\Theta_S^T(H/F)$.
\\

Corollary \ref{maps-Fitt-Selmer} combined with \eqref{base-change-Fitt} shows that we have an inclusion of $R$--ideals
$${\rm Fitt}_{R}\,(Sel_{S'}^T(H')_{R}) =\rho({\rm Fitt}_{R'}\,(Sel_{S'}^T(H')))\subseteq {\rm Fitt}_{R}\,(Sel_{S'}^T(H)_p^-) .$$
The argument leading to the proof of \eqref{Fitt-S_ast-S}, applied to $S$ and $S'=S\cupdot S_r(F)$, gives an equality of $R$--ideals
$${\rm Fitt}_{R}\,(Sel_{S'}^T(H)_p^-)=u_r(H/F)\cdot {\rm Fitt}_{R}\,(Sel_{S}^T(H)_p^-),$$
where $u_r(H/F):=\prod_{v\in S_r(F)}(1-\sigma_v(H/F)^{-1}).$ Obviously, we also have an equality
$$\Theta_{S'}^T(H/F)=u_r(H/F)\cdot\Theta_S^T(H/F).$$
When combining \eqref{R-Fitting-equality} with the last three displayed relations, we obtain an inclusion of $R$--ideals
$$u_r(H/F)\cdot(\Theta_S^T(H/F))\subseteq u_r(H/F)\cdot {\rm Fitt}_{R}\,(Sel_{S}^T(H)_p^-) + p^M R.$$
However, since the auxiliary prime $r$ is $H/F$--good (see Proposition \ref{Tchebotarev}(2)), Remark \ref{remark-bad} shows that $u_r(H/F)$ is a unit in the ring $R$. Therefore, we can cancel it to obtain an inclusion
$$(\Theta_S^T(H/F))\subseteq {\rm Fitt}_{R}\,(Sel_{S}^T(H)_p^-) + p^M R.$$
However, since the integer $M$ can be taken to be arbitrarily large (see \eqref{N-M}) and the Noetherian $\Bbb Z_p$--algebra  $R$ is $p$--adically complete and all its (finitely generated) ideals are closed in the $p$--adic topology, the last inclusion implies that, in fact, we have an inclusion of $R$--ideals.
 $$(\Theta_S^T(H/F))\subseteq \bigcap_{M\gg 0} ({\rm Fitt}_{R}\,(Sel_{S}^T(H)_p^-) + p^M R)={\rm Fitt}_{R}\,(Sel_{S}^T(H)_p^-),$$
 which concludes our {\bf Step 1}.
\medskip

{\bf Step 2.} Next, we prove that the inclusion \eqref{inclusion1} is an equality, as stated in Theorem \ref{BKS-Conj}(1).  Let $G=G'\times G_p$, where $G_p$ is the $p$--Sylow subgroup of $G$, and $R:=\Bbb Z_p[G]^-$. For any odd character $\chi\in\widehat{G'}$ and any character $\psi\in\widehat G$, such that $\psi\vert_{G'}=\chi$ (i.e. $\psi$ ``belongs'' to $\chi$), we let $R^\chi:=\Bbb Z_p[\chi][G_p]$ and $R^\psi=\Bbb Z_p[\psi]$, viewed as $R$--algebras via the usual surjective $\Bbb Z_p$--algebra morphisms $\chi$ and $\psi$ (abuse of notation!), induced by the characters $\chi$ and $\psi$, respectively:
$$\psi:\, R\overset{\chi}\twoheadrightarrow R^\chi\twoheadrightarrow  R^\psi.$$
Since we have an obvious $\Bbb Z_p$--algebra isomorphism $R\simeq \bigoplus_{[\chi]}R^\chi$ (sum taken over all $G(\overline{\Bbb Q_p}/\Bbb Q_p)$--conjugacy classes of odd characters $\chi\in \widehat{G'}$), we need to show that we have an equality
$$(\Theta_{S, H/F}^T)_{R^\chi}= {\rm Fitt}_{R^\chi}(Sel_S^T(H)_{R^\chi}),$$
for all characters $\chi$ as above. 

Let us fix $\chi$ as above. Lemma 4.1 in \cite{Dasgupta-Kakde} implies that if we have  $\psi(\Theta_{S, H/F}^T)=0$, for all $\psi\in\widehat G$ which belong to $\chi$, then we have a double equality
$$(\Theta_{S, H/F}^T)_{R^\chi}= (0)= {\rm Fitt}_{R^\chi}(Sel_S^T(H)_{R^\chi}).$$
Consequently, we may assume that there exists a character $\psi$ which belongs to $\chi$, such that 
$$\psi(\Theta_{S, H/F}^T)\ne 0.$$ We fix such a character $\psi$.
Since the $R$--module $Sel_S^T(H)_p^-$ is quadratically presented (see Lemma \ref{quadratic-pres-lemma}(2)), its Fitting ideal is principal. Hence, inclusion \eqref{inclusion1} shows that there exists $x\in R$, such that 
$$(\Theta_{S,H/F}^T)=(x)\cdot {\rm Fitt}_{R}(Sel_S^T(H)_p^-), \qquad (\Theta_S^T(H/F))_{R^\chi}=({\chi}(x))\cdot {\rm Fitt}_{R^\chi}(Sel_S^T(H)_{R^\chi}).$$
Therefore, if we show that $\chi(x)$ is a unit in $R^\chi$, we are done. However, the rings $R^\chi$ and $R^\psi$ are local, of maximal ideals $\mathfrak m_\chi:=(p, I_{G_p})$ (here $I_{G_p}$ is the augmentatation ideal, generated by $\{\sigma-1\vert \sigma\in G_p\}$), and $\mathfrak m_\psi=\psi(\mathfrak m_\chi)$,  respectively, we have equivalences
\begin{equation}\label{double-equivalence}\chi(x)\in (R^\chi)^\times \iff \psi(x)\in R_\psi^\times \iff (\psi(\Theta_S^T(H/F)))={\rm Fitt}_{R^\psi}(Sel_S^T(H)_{R^\psi}).\end{equation}

In what follows, we will use ideas developed in the proof of Theorem \ref{partial} above to prove the rightmost statement in the above double equivalence. This will conclude {\bf Step 2} and the proof of Theorem \ref{BKS-Conj}. More precisely, we prove the following.

\begin{proposition}\label{final-prop} For all $(H/F, S, T, p)$ as in Theorem \ref{BKS-Conj}(1),  we let 
$$\Sigma_0:=S\setminus ((S_{ram}(H/F)\setminus S_{p\infty}(F))\cap S), \qquad \Sigma_0':=T\cupdot ((S_{ram}(H/F)\setminus S_{p\infty}(F))\cap S).$$
Then, the following hold, for all odd characters $\psi\in\widehat G$.
\begin{enumerate}
\item For all sets  $J\subseteq (S_{ram}(H/F)\setminus S_{p\infty}(F))\cap S$, we have an equality of $R^\psi$--ideals
$${\rm Fitt}_{R^\psi}(Sel_{\Sigma_0\cup J}^{\Sigma_0'\setminus J}(H/F)_{R^\psi})=(\psi(\Theta_{\Sigma_0\cup J}^{\Sigma_0'\setminus J}(H/F))).$$
\item We have an equality of $R^\psi$--ideals 
$${\rm Fitt}_{R^\psi}(Sel_S^T(H)_{R^\psi})=(\psi(\Theta_{S, H/F}^T)).$$
\end{enumerate}
\end{proposition}
\begin{proof}
First, note that part (2) (which is of direct interest to us, as it implies the rightmost statement in  the double equivalence \eqref{double-equivalence}) follows from part (1), if we set $J:=(S_{ram}(H/F)\setminus S_{p\infty}(F))\cap S$.
\\

We prove the statement in part (1), for all $(H/F,\Sigma_0,\Sigma_0',J)$, by induction on $|J|$, as we proceeded in the proof of Theorem \ref{partial}. 
In the case  $J=\emptyset$, the proof of Theorem \ref{empty-J-theorem} (applied to $(H/F, \Sigma_0, \Sigma_0')$ instead of $(H_n/F, \Sigma, \Sigma')$) 
produces an equality of $R$--ideals
$${\rm Fitt}_R(Sel_{\Sigma_0}^{\Sigma_0'}(H)_p^-)=(\Theta_{\Sigma_0}^{\Sigma_0'}(H/F)).$$
This implies the desired equality of $R^\psi$--ideals by applying $\psi$ on both sides, for all odd characters $\psi\in\widehat G$.
\\

Now, let us assume that the result holds for all $(H/F,\Sigma_0,\Sigma_0',J)$, with $J$ of a fixed size. 
For fixed data as above, let  $v\in ((S_{ram}(H/F)\setminus S_{p\infty}(F))\cap S)\setminus J$. We need to prove the result for the data $(H/F, \Sigma_0, \Sigma_0', J_v:= J\cup \{v\})$.
Let $I_v$ be the ramification group of $v$ and let $I_v':=I_v\cap G'$. Let $\chi\in\widehat{G'}$, such that  $\psi|_{G'}=\chi$. As in the proof of Theorem \ref{partial}, we have two cases.
\\

{\bf Case (1):} $\chi(I_v')\neq\{1\}.$ The proof of {\bf Case 1} in Theorem \ref{partial} (applied to the data $(H/F, \Sigma_0, \Sigma_0')$ instead of $(H_n/F, \Sigma, \Sigma')$) gives equalities of $R^\chi$--ideals 
$${\rm Fitt}_{R^\chi}(Sel_{\Sigma_0\cup J}^{\Sigma_0'\setminus J}(H/F)_{R^\chi})={\rm Fitt}_{R^\chi}(Sel_{\Sigma_0\cup J_v}^{\Sigma_0'\setminus J_v}(H/F)_{R^\chi}),  \qquad (\Theta_{\Sigma_0\cup J}^{\Sigma_0'\setminus J}(H/F))_{R^\chi}=
(\Theta_{\Sigma_0\cup J_v}^{\Sigma_0'\setminus J_v}(H/F))_{R^\chi}.$$
The desired equality of $R^\psi$--ideals is obtained by applying $\psi$ to the two equalities above and using the induction hypothesis.
\\

{\bf $\bullet$ Case (2):} $\chi(I_v')=\{1\}.$ We follow the initial steps of the proof of {\bf Case 2} in Theorem \ref{partial}, applied to the data $(H/F, \Sigma_0, \Sigma_0')$ instead of $(H_n/F, \Sigma, \Sigma')$. After evaluation against $\psi$, inclusion \eqref{Case 2-inclusion} reads as follows

\begin{eqnarray}
\nonumber {\rm Fitt}_{{R^\psi}}(Sel_{\Sigma_0\cup J_v}^{\Sigma_0'\setminus J_v}(H)_{R^\psi}) &\subseteq & {\rm Fitt}_{R^\psi}(Sel_{\Sigma_0\cup J}^{\Sigma_0'\setminus J_v}(H)_{R^\psi})=\\
\nonumber &=&(\psi(\Theta_{\Sigma_0\cup J}^{\Sigma_0'\setminus  J_v}(H/F))(\psi(NI_{v}), \psi(1-\sigma_v^{-1} e_{v})),
\end{eqnarray}
where the notations are as in loc.cit.

Now, if $\psi$ is ramified at $v$, we have $\psi(NI_v)=0$ and then the above inclusion implies,
$${\rm Fitt}_{R^\psi}(Sel_{\Sigma_0\cup J_v}^{\Sigma_0'\setminus J_v}(H)_{\psi}) \subseteq (\psi(\Theta_{\Sigma_0\cup J_v}^{\Sigma_0'\setminus  J_v}(H/F))$$
Since we have already proved the reverse inclusion, we are done proving the desired equality of $R^\psi$--ideals.

Next, assume that $\psi$ is not ramified at $v$. Therefore, we have $\psi(I_v)=\{1\}$. We let $H_{\psi}:=H^{ker(\psi)}$. By Lemma 5.4 of \cite{Dasgupta-Kakde}, we have an isomorphism of $R^\psi$--modules 
$$Sel_{\Sigma_0\cup J_v}^{\Sigma_0'\setminus J_v}(H)_{R^\psi}\cong Sel_{\Sigma_0\cup J_v}^{\Sigma_0'\setminus J_v}(H_{\psi})_{R^\psi}$$ 
induced by the natural restriction map. Observe that $v$ is unramified in $H_{\psi}/F$. Therefore, we may apply the induction hypothesis to the data
$(H_{\psi}/F, \Sigma_0\cup \{v\}, \Sigma_0'\setminus \{v\}, J)$ to obtain equalities of $R^\psi$--ideals
$${\rm Fitt}_{R^\psi}(Sel_{\Sigma_0\cup J_v}^{\Sigma_0'\setminus J_v}(H)_{R^\psi})={\rm Fitt}_{R^\psi}(Sel_{\Sigma_0\cup J_v}^{\Sigma_0'\setminus J_v}(H_{\psi})_{R^\psi})=(\psi(\Theta_{\Sigma_0\cup J_v}^{\Sigma_0'\setminus  J_v}(H_{\psi}/F)))=(\psi(\Theta_{\Sigma_0\cup J_v}^{\Sigma_0'\setminus  J_v}(H/F))).$$
Above, the first equality is due to the module isomorphism mentioned above, the second comes from the induction hypothesis, and the third is due to the inflation property of Artin L-functions. This completes the proof of the inductive step and that of Proposition \ref{final-prop}.
\end{proof}
As explained above, Theorem \ref{BKS-Conj}(1) follows from Proposition \ref{final-prop}(2) and equivalences \eqref{double-equivalence}.
\end{proof}

\bibliographystyle{amsplain}
\bibliography{EMCIIBib}
\end{document}